\documentclass[final, leqno]{siamltex}
\usepackage{amsmath,amssymb}
\usepackage[pdftex]{graphicx,color}
\usepackage{epstopdf}
\usepackage{multicol}
\usepackage{enumerate}
\usepackage{float}
\usepackage{hyperref} % must be last package
\hypersetup{
 colorlinks=true,
 linktocpage=true,
 bookmarksnumbered=true
}
\usepackage{algorithm,algorithmic}
\newtheorem{thm}{Theorem}[section]
\newtheorem{rem}[thm]{Remark}

\newtheorem{ex}[thm]{Example}

\newcommand*{\rom}[1]{\expandafter\@slowromancap\romannumeral #1@}
\newcommand{\seqof}[3]{(#1)_{#2}^{#3}}

\newcommand{\MLE}{\text{MLE}}

\newcommand{\avg}[2]{\mathcal{A}(#1;#2)}
\newcommand{\avgsub}[3]{\mathcal{A}_{#3}(#1;#2)}
\newcommand{\sdev}[2]{\mathcal{S}(#1;#2)}

\newcommand{\dbar}[1]{\Bar{\Bar{#1}}}

\newcommand{\ie}{\emph{i.e.}}

\newcommand{\ud}{\mathrm{d}}
\newcommand{\prob}[1]{\mathrm{P}\left(#1\right)}
\newcommand{\expt}[1]{\mathrm{E}\left[#1\right]}

\newcommand{\norm}[1]{\left\|#1\right\|}
\newcommand{\abs}[1]{\left|#1\right|}
\newcommand{\rset}{\mathbb{R}}

\newcommand{\zset}{\mathbb{Z}}

\newcommand{\ordo}[1]{{o}\left(#1\right)}

\newcommand{\thx}[1]{\thanks{#1}}

\newcommand{\latt}{\mbox{$\zset_+^d$}}

\newcommand{\indicator}[1]{\mathbf{1}_{\left\{#1\right\}}} 

\newcommand{\PERIOD}{.}
\newcommand{\COMMA}{,}

\newcommand{\LP}{\left(}
\newcommand{\RP}{\right)}

\newcommand{\SEP}{\, \big| \,}

\title{An Efficient Forward-Reverse Expectation-Maximization Algorithm for Statistical Inference in Stochastic Reaction Networks}

\author{Christian Bayer\thanks{Weierstrass Institute for Applied Analysis and Stochastics,
 Berlin, Germany ({\tt christian.bayer@wias-berlin.de}).}\and
Alvaro Moraes\thanks{Computer, Electrical and Mathematical Sciences and Engineering,
 King Abdullah University of Science and Technology (KAUST),
 Thuwal, Saudi Arabia ({\tt alvaro.moraesgutierrez@kaust.edu.sa}).}
 \and Raul Tempone\thanks{Mathematical and Computer Sciences and Engineering Division,
 King Abdullah University of Science and Technology (KAUST),
 Thuwal, Saudi Arabia ({\tt raul.tempone@kaust.edu.sa}).}
\and Pedro Vilanova\thanks{Mathematical and Computer Sciences and Engineering Division,
 King Abdullah University of Science and Technology (KAUST),
 Thuwal, Saudi Arabia ({\tt pedro.guerra@kaust.edu.sa}).}}

\begin{document}

\maketitle

\begin{abstract}
% flatex input: [abstract.tex]
In this work, we present an extension to the context of Stochastic Reaction Networks (SRNs) of the forward-reverse representation introduced in ``Simulation of forward-reverse stochastic representations for  conditional diffusions'', a 2014 paper by Bayer and Schoenmakers. 
We apply this stochastic representation in the computation of efficient  approximations of expected values of functionals of SNR bridges, \ie, SRNs conditioned to its values in the extremes of given time-intervals.
We then employ this SNR bridge-generation technique to the statistical inference problem of approximating the reaction propensities based on discretely observed data. To this end, we introduce a two-phase iterative inference method in which, during phase I, we  solve a set of deterministic optimization problems where the SRNs are replaced by their reaction-rate Ordinary Differential Equations (ODEs) approximation; then, during phase II, we apply the Monte Carlo version of the Expectation-Maximization (EM) algorithm starting from the phase I output. 
By selecting a set of over dispersed seeds as initial points for  phase I, the output of parallel runs from our two-phase method is a cluster of approximate maximum likelihood estimates. 
Our results are illustrated by numerical examples.

% flatex input end: [abstract.tex]
  
%\title{An Efficient Forward-Reverse Algorithm for Stochastic Reaction Networks with Applications to Statistical Inference}
\end{abstract}



\pagestyle{myheadings}
\thispagestyle{plain}

\setcounter{tocdepth}{1}

\section{Introduction}
% flatex input: [introduction.tex]

Stochastic Reaction Networks (SRNs) are a class of continuous-time Markov
chains, $X{\equiv}\{X(t)\}_{t\in[0,T]}$, that take values in $\latt$, \ie, the
lattice of $d$-tuples of non-negative integers.  SRNs are mathematical models
employed to describe the time evolution of many natural and {artificial}
systems. Among them we find biochemical reactions, spread of epidemic
diseases, communication networks, social networks, transcription and
translation in genomics, and virus kinetics.

For historical reasons, the jargon from chemical kinetics is used to describe
the elements of SRNs.  The integer $d{\geq}1$ is the number of chemical species
reacting in our system. The coordinates of the Markov chain,
$X(t){=}(X_1(t),\ldots,X_d(t))$, account for the number of molecules or
individuals of each species present in the system at time $t$.  The
transitions in our system are given by a finite number $J$ of 
\emph{reaction channels}, $\seqof{\mathcal{R}_j}{j=1}{J}$.  Each reaction channel
$\mathcal{R}_j$ is a pair formed by a vector $\nu_j$ of $d$ integer components and a non-negative function $a_j(x)$ of the state of the system. Usually, $\nu_j$ and $a_j$ are named \emph{stoichiometric vector} and
\emph{propensity function}, respectively.  Because our state space is a lattice, our
system evolves in time by jumping from one state to the next, and for that reason
$X$ is a pure jump process.

The propensity functions, $a_j$, are usually derived through \emph{the mass action 
principle} also known as \emph{the law of mass action}, see for instance Section 3.2.1 in \cite{Holmes}. 
For that reason,
we assume that $a_j(x) {=} c_j\, g_j(x)$, where $c_j$ is a non negative
coefficient and $g_j(x)$ is a given monomial in the coordinates of the
process, $X$.  However, our results can be easily extended to polynomial
propensities.

In this work, we address the statistical inference problem of estimating the
coefficients $\theta {=} (c_1,\ldots,c_J)$ from \emph{discretely observed data}, \ie,
data collected by observing one or more paths of the process $X$ at
a certain finite number of \emph{observational times} or epochs.  It means that our data, $\mathcal{D}$, is a finite collection  
$\{(t_{n,m},x(t_{n,m}))\}$, where  $m{=}1,2,\ldots,M$ indicates the
observed path,  $n{=}1,2,\ldots,N(m)$ indicates the $n$-th observational 
time corresponding to the $m$-th path, and the datum $x(t_{n,m})$ can be considered as 
an observation of the $m$-path of the process $X$ at time time $t_{n,m}$. 
The observational times, $t_{n,m}$,
are either deterministic or random but independent from the state of the
process $X$. 
In what follows, we denote
with $X_{i,n,m}$ the $i$-th coordinate of $X(t_{n,m},\omega_m)$, 
with $X_{\cdot,n,m}$ the vector $X(t_{n,m},\omega_m)$, where
$\omega_m$ is the $m$-th path of the process $X$.

{Let us remark that we observe all the coordinates of $X$ and not only a fixed subset at each observational time $t_{n,m}$. In that sense, we are not treating the case of \emph{partially observed data} where only a fixed proper subset of coordinates of $X$ is observed.}
\begin{rem}
The partially observed case can in principle also be treated by a variant of the FREM algorithm based on \cite{Bayer}  (Corollary 3.8).
\end{rem}

For further convenience, we organize the information in our data set,
$\mathcal{D}$, as a finite collection,
\begin{align}\label{def:data}
\mathcal{D} = \seqof{[s_k,t_k],x(s_k),x(t_k)}{k=1}{K}\COMMA
\end{align} 
such that for each $k$, $I_k:=[s_k,t_k]$ is the time interval determined by
two consecutive observational points $s_k$ and $t_k$, where the states
$x(s_k)$ and $x(t_k)$ have been observed.
{Notice that the set 
$\mathcal{D}$  collects all the data corresponding to the $M$ observed paths of the process $X$.
For that reason, it is possible to have $[s_k,t_k]{=}[s_{k'},t_{k'}]$ for $k{\neq} k'$, for instance, in the case of repeated measurements.}  

For technical reasons,
we need to define a sequence of \emph{intermediate times}, $\seqof{t_k^*}{k=1}{K}$;
for instance, $t_k^*$ could be the midpoint of $[s_k,t_k]$.

It turns out that the likelihood function, $\text{lik}^c(\theta)$,
corresponding to data obtained from continuously observed paths of $X$ is
relatively easy to derive (see Section \ref{sec:contobservedpaths}). It depends on the total
number of times that each reaction channel fires over the time interval
$[0,T]$ and the values of the monomials $g_j$ evaluated at the jump times of
$X$.  
Since the observational times, $t_{n,m}$, are not necessarily equal to the jump times of the process $X$,
we can not directly deal  with  the likelihood $\text{lik}^c(\theta)$. 
For that reason, we consider the Monte Carlo version of the expectation-maximization (EM) algorithm \cite{Dempster77,Casella, WatanabeYamaguchi, McLachlanEM} in which we treat the jump times of $X$ and their corresponding reactions as missing data.  
The ``missing data'' can be gathered by simulating \emph{SRN bridges} of the process  $X$  conditional on  
$\mathcal{D}$,  \ie, $X(s_k){=}x(s_k)$ and $X(t_k){=}x(t_k)$ for all intervals $[s_k,t_k]$.
To simulate SRN bridges, we extend the \emph{forward-reverse} technique developed by Bayer and Schoenmakers \cite{Bayer} for It\^o diffusions to the case of SRNs.
As explained in Section \ref{sec:forwardreverse}, the forward-reverse algorithm generates forward paths from $s_k$ to $t_k^*$ and backward paths from $t_k$ to $t_k^*$. An exact SRN bridge is formed when forward and backward paths meet at 
$t_k^*$. Observe that the probability of producing SRN bridges strongly depends on the approximation of $\theta$ that we use to generate the forward and backward paths.  In addition to exact bridges, in this work we also relax this meeting condition by using a kernel $\kappa$.  

{Here, we present a two-phase algorithm that approximates the Maximum Likelihood Estimator, $\hat{\theta}_{\MLE}$,  of the vector $\theta$ using the collected data, $\mathcal{D}$.

Phase I is the result of a deterministic procedure while  phase II is the result of a stochastic one.
The purpose of  phase I is to generate an estimate of  $\theta$ that will be used as initial point for  phase II.
To this end, in the phase I  we solve a deterministic global optimization problem obtained by substituting at each time interval, $[s_k,t_k]$, 
the ODE approximations to the mean of the forward and reverse stochastic paths and minimizing a weighted sum of the squares of the Euclidean distances of the ODE approximations at the times $t^*_k$. Using this value as a starting point for phase II, we hope to simulate an acceptable number of SRN bridges in the interval $[s_k,t_k]$ without too much computational effort. 
Phase I starts at ${\theta^{(0)}_{I}}$ and provides 
$\theta^{(0)}_{I\!I}$. 
In phase II  we run a Monte Carlo EM stochastic sequence $\seqof{\hat{\theta}^{(p)}_{I\!I}}{p=1}{+\infty}$ until a certain convergence criterion is fulfilled. Here we have a schematic representation of the two-phase method:
%See Figure \ref{fig:scheme} for a schematic representation of our two-phase approach.
\begin{equation*}
\theta^{(0)}_{I} \rightarrow \theta^{(0)}_{II}\rightarrow
\hat{\theta}^{(1)}_{II} \rightarrow \cdots \,\, 
\hat{\theta}^{(p)}_{II}\rightarrow \cdots \rightarrow \hat \theta\PERIOD
\end{equation*}

%\begin{figure}[h!]
%\centering
%$
%\theta^{(0)}_{I} \rightarrow \theta^{(0)}_{II}\rightarrow
%\hat{\theta}^{(1)}_{II} \rightarrow \cdots \,\, 
%\hat{\theta}^{(p)}_{II}\rightarrow \cdots \rightarrow \hat \theta
%$
%\caption{The two-phase estimation process. In the first step, we obtain $\theta^{(0)}_{II}$ from $\theta^{(0)}_{I}$ by solving an optimization problem \eqref{eq:seedI}). 
%In the subsequent steps, we generate the stochastic sequence 
%$\seqof{\theta^{(p)}_{II}}{p=1}{+\infty}$ using Monte Carlo EM \eqref{eq:MCEMiteration}.}
%\label{fig:scheme}
%\end{figure}

During phase II, we intensively use a computationally efficient implementation of the SRN-bridge simulation algorithm for simulating the ``missing data'' that feeds the Monte Carlo EM algorithm. Details are provided in Section \ref{sec:FREM}.
Our two-phase algorithm is named FREM as the acronym for Forward-Reverse Expectation Maximization. } 

Although our FREM algorithm has certain similarity with the estimation methodology proposed in \cite{daigle2012accelerated}, there are also notable differences.
In terms of the similarity, in \cite{daigle2012accelerated} the authors propose a two-phase method where the first phase is intended to select a seed for the second phase, which is an implementation of the Monte Carlo EM algorithm. 
While our first phase is deterministic and uses the reaction-rate ODEs as approximations of the SRN paths,
theirs is stochastic and a number of parameters should be chosen to determine the amount of computational work and the accuracy of the estimates.
There is also a main difference is the implementation of the second phase: 
while the FREM algorithm is focused in efficiently generating kernel-based SRN bridges using the novel forward-reverse technology introduced by Bayer and Schoenmakers in \cite{Bayer}, the authors of \cite{daigle2012accelerated} propose a trial-and-error shooting method for sampling SRN bridges. This  shooting method can be viewed as a particular case of the FREM algorithm by systematically choosing the intermediate point $t^*_k$ as the right extreme point $t_k$, giving no place for  backward paths. 
To quantify the uncertainty in our estimates, we prefer to have the outputs of our algorithm starting from a set of over dispersed initial points without assuming Gaussianity in its distribution (see \cite{Casella}). 
The variance of our estimators can be easily assessed by bootstrap calculations. In our numerical experiments, we observe that the outputs lie on a low-dimensional manifold in parameter space; this is a motivation against the use of the Gaussiantiy assumption.
Regarding the stopping criterion proposed in \cite{daigle2012accelerated}, we found that the condition imposed there, of obtaining three consecutive iterations close to each other up to a certain tolerance, could be considered as a rare event in some examples and it may lead to the generation of an excessive number of Monte Carlo EM iterations. We refer to \cite{daigle2012accelerated} for comparisons against other existing related statistical inference methods for SRNs.

In \cite{wang} the authors propose a method based on maximum likelihood for parameter inference. It is based on first estimating the gradient of the likelihood function with respect to the parameters by using reversible-jump Markov chain Monte Carlo sampling (RJMCMC) \cite{green95,BoysEtAl2008} and then applying a gradient descent method to obtain the maximum likelihood estimation of the  parameter values. The authors provide a formula for the gradient of the likelihood function given the observations.
The idea of the RJMCMC method is to generate an initial reaction path and then generate  new samples by adding or deleting a set of reactions from the path using an acceptance method. The authors propose a general method for obtaining a sampler that can work for any reaction system. This sampler can be inefficient in the case of large observation intervals. At this point, we would like to observe that their approach can be combined with ours if, instead of using the RJMCMC method for computing the gradient of the likelihood function, we use our forward-reverse method. 
We think that this combination may be useful in cases in which many iterations of our method are needed (see Section \ref{ex:bd} for such an example). This is left as future work.

In the remainder of this section, we formally introduce SRNs and their reaction-rate ODE approximations, the stochastic simulation algorithm and the forward-reverse method. In Section \ref{sec:forwardreverse}, we develop the main result of this article: the extension of the forward-reverse technique to the context of SRNs. The EM algorithm for SRNs is introduced in Section \ref{sec:EM}. Next, in Section \ref{sec:FREM}, we introduce the main application of this article: the forward-reverse EM (FREM) algorithm for SRNs. In Section \ref{sec:compdetails}, we provide computational details for the practical implementation of the FREM algorithm. Later, in Section \ref{sec:numex}, we present numerical examples to illustrate the FREM algorithm and finally, we present our conclusions in Section \ref{conclusions}. 
Appendix \ref{sec:algorithms} contains the pseudo-code for the implementation of the FREM algorithm.

\subsection{Stochastic Reaction Networks}
Stochastic Reaction Networks are continuous time Markov chains, $X:[0,T]\times \Omega \to\latt$, that describe the stochastic evolution of a system of $d$ interacting species.
%For the sake of brevity, we write $X(t,\omega) \equiv X(t)$.  
In this context, the $i$-th coordinate of the process $X$, $X_i(t)$, can be interpreted as the
number of individuals of species $i$ present in the system at time $t$. 

The system evolves randomly through $J$ different reaction channels $\mathcal{R}_j:=(\nu_j,a_j)$. 
Each stoichiometric vector $\nu_j{\in}\zset^d$ represents a possible jump of the system, $x \rightarrow x{+}\nu_j$. 
The probability that the reaction $j$ occurs during an infinitesimal interval
$(t,t+\ud t)$ is given by
\begin{equation}\label{eq:infdefX}
  \prob{\text{reaction } j \text{ fires during } (t,t+\ud t) \bigm|
    X(t) = x} =
  a_j(x) \ud t + \ordo{\ud t},
\end{equation}
where $a_j:\rset^d \to [0,\infty)$ are known as propensity functions. 
We set $a_j(x){=}0$ for those $x$ such that $x{+}\nu_j\notin \latt$.
We assume that the initial condition of $X$, $X(0)=x_0\in\latt$ is deterministic and known.
The \emph{stoichiometric matrix} $\nu$ is defined as the matrix whose $j$-column is $\nu_j$ ($\nu^T$ denotes its transpose). The \emph{propensity vector} $a(x) \in \rset^J$ has $a_j(x)$ as components.

\begin{ex}[Simple decay model]\label{ex:AB}
 Consider the reaction $X \xrightarrow{c} \emptyset$ where one particle is
  consumed. In this case,
  the state vector $X(t)$ is in $\zset_+$ where $X$  denotes the number
  of particles in the system. The vector for this
  reaction is $\nu = -1$. 
  The propensity functions in this case could be, for example, $a(X)= c\,X$, where $c>0$.
\end{ex}
Section \ref{sec:numex} contains more examples of stochastic reaction networks.

\subsection{ Deterministic Approximations of SRNs}
The infinitesimal generator $\mathcal{L}_X$ of the process $X$
is a linear operator defined on the set of bounded functions \cite{kurtzmp} .
In the case of SRN, it is given by 
\begin{equation}\label{eq:genX}
\mathcal{L}_X(f)(x) := \sum_j a_j(x) ({f(x+\nu_j)-f(x)})\PERIOD
\end{equation}
The Dynkin formula, (see \cite{Klebaner})
\begin{equation} \label{eq:Dynkin}
\expt{f(X(t))} = f(X(0)) + \int_0^t \expt{\mathcal{L}_X(f) (s)}\ud s\COMMA
\end{equation}
can be used to obtain integral equations describing the time evolution of any observable of the process $X$. 
In particular, taking the canonical projections $f_i(x)=x_i$, we obtain a system of equations for $\expt{X_i(t)}$,
\begin{align*}
\expt{X_i(t)} = x_0 + \int_0^t \sum_j \expt{a_j(X(s))} \nu_{j,i} \ud s.%\\[5.9pt]
\end{align*}
%{\color{green}
If all the propensity functions, $a_j$, are affine functions of the state, then this system of equations leads to a closed system of ODEs.
%}
In general, some propensity functions may not depend on their
coordinates $x$ in an affine way, and for that reason, the integral equations for
$\expt{X_i(t)}$ obtained from the Dynkin formula depend on higher moments of
$X$. This can be treated  using moment  closure techniques \cite{MomentClousure, MomentClousure2} or by
taking a different approach: using a formal first-order Taylor expansion of
$f$ in \eqref{eq:genX}, we obtain the generator
\begin{equation*}
\mathcal{L}_Z(f)(x) := \sum_j a_j(x) {\partial_x f(x) \nu_j }\COMMA
\end{equation*} 
which corresponds to the reaction-rate ODEs (also known as the {mean field}
ODEs) 
%\todo{CB: I guess this is helpful when $f$ is non-linear. But how is it
%necessary when $f$ is already linear like here?. AM: The problem arises when %some $a_j$ is not affine, for instance $a_j(x)=x_1\times x_2$ like in %epidemics. I will send you a picture by email.}
\begin{align}\label{eq:ODE}
\left\{ \!
\begin{array}{l@{\;}c@{\;}l}
dZ(t) &=& \nu a(Z(t)) dt , \,\, t \in \rset_+, \\
Z(0) &=& x_0,
\end{array}
\right.
\end{align}
where the $j$-column of the matrix $\nu$ is $\nu_j$ and $a$ is a column vector with components $a_j$.

This derivation motivates the use of $Z(t)$ as an approximation of $\expt{X(t)}$
in phase I of our FREM algorithm.

\subsection{The Stochastic Simulation Algorithm}
\label{sec:SSA_num_approx}
To simulate paths of the process $X$, we employ the stochastic simulation algorithm (SSA) by Gillespie \cite{Gillespie1976}.
The SSA  simulates statistically exact paths of $X$, \ie, the probability law of any path generated by the SSA 
satisfies (\ref{eq:infdefX}). 
It requires one to sample two independent uniform random variables per time step: one is used to find
the time of the next reaction and the other to determine which is the
reaction that fires at that time.  
Concretely, given the current state of the system, $x:= X(t)$, we simulate 
two independent uniform random numbers, $U_1,U_2 \sim \mathcal{U}(0,1)$
and compute:
\begin{equation*}
 j =  \min \Big \{ k\in \{1,\ldots,J\}: \sum_{i=1}^{k} {a_i(x)} {>}
    U_1\, {a_0(x)}\Big\}
    \COMMA \, \, \tau_{\min} = -\LP a_0(x)\RP ^{-1} \ln \left( U_2 \right)\COMMA
\end{equation*}
where $a_0(x):=\sum_{j=1}^J a_j(x)$. 
The system remains in the state $x$ until the time $t+\tau_{\min}$ when it jumps, $X(t+\tau_{\min})= x+\nu_j$.
In this way, we can simulate a full path of the process $X$.

{Exact paths can be generated using more efficient algorithms like the modified next reaction method by Anderson \cite{Anderson2007}, where only one uniform variate is needed at each step. However, in regimes where the total propensity, $a_0(x)$, is high, approximate path-simulation methods like the hybrid Chernoff tau-leap \cite{ourSL} or its multilevel versions \cite{ourML,ourMixed} may be required.}

\subsection{Bridge Simulation for SDEs}\label{sec:orders}
In \cite{Bayer}, Bayer and Schoenmakers introduced the
so-called forward-reverse algorithm for computing conditional expectations of
path-dependent functionals of a diffusion process conditioned on the values of
the diffusion process at the end-points of the time interval. More precisely,
let $X = X(t)$, $0 \le t \le T$, denote the solution of a $d$-dimensional
stochastic differential equation (SDE) driven by standard Brownian motion. Under
mild regularity conditions, a \emph{stochastic representation} is provided for
conditional expectations of the form,
\begin{equation*}
  \mathcal{H} \equiv \expt{\left. g(X) \ \right| \ X_0 = x, \, X_T = y }\COMMA
\end{equation*}
for fixed values $x, y \in \mathbb{R}^d$ and a (sufficiently regular)
functional $g$ on the path-space.\footnote{In fact, Bayer and Schoenmakers
  \cite{Bayer} require $g$ to be a smooth function of the values $X_{t_i}$
of the process $X$ along a grid $t_i$, but a closer look at the paper reveals
that more general, truly path-dependent functionals can be allowed.} More
precisely, they prove an limiting equality of the form
\begin{equation}
  \label{eq:forrev-sde}
  \mathcal{H} = \frac{\lim_{\epsilon \to 0} \expt{ g( X^{(f)} \circ X^{(b)})
      \kappa_\epsilon(X^{(f)}(t^\ast) - X^{(b)}(t^\ast)) \mathcal{Y} }}{
    \lim_{\epsilon \to 0} \expt{ \kappa_\epsilon(X^{(f)}(t^\ast) -
      X^{(b)}(t^\ast)) \mathcal{Y} } }\PERIOD
\end{equation}
Here, $X^{(f)}$ is the solution of the original SDE (i.e., is a copy of $X$)
started at $X^{(f)}(0) = x$ and solved until some time $0 < t^\ast <
T$. $X^{(b)}$ is the time-reversal of another diffusion process $Y$ whose
dynamics are again given by an SDE (with coefficients \emph{explicitly} given
in terms of the coefficients of the original SDEs) started at $Y(t^\ast) = y$
and run until time $T$. Hence, $X^{(b)}$ starts at $t^\ast$ and ends at
$X^{(b)}(T) = y$. We then evaluate the functional $g$ on the ``concatenation''
$X^{(f)} \circ X^{(b)}$ of the paths $X^{(f)}$ and $X^{(b)}$, which is a path
defined on the full interval $[0,T]$ defined by
\begin{equation*}
  X^{(f)} \circ X^{(b)} (s) \equiv
  \begin{cases}
    X^{(f)}(s), & 0 \le s \le t^\ast, \\
    X^{(b)}(s), & t^\ast < s \le T.
  \end{cases}
\end{equation*}
In particular, we remark that $X^{(f)} \circ X^{(b)}$ may exhibit a jump at
$t^\ast$. Here, $\mathcal{Y}$ is an exponential weighting term of the
form $\mathcal{Y} = \exp\left( \int_{t^\ast}^T c(Y_s) ds \right)$. At last,
$\kappa_\epsilon$ denotes a \emph{kernel} with bandwidth $\epsilon >
0$. Notice that the processes $X^{(f)}$ and the pair $\left( X^{(b)},
  \mathcal{Y} \right)$ are chosen to be independent.

Let us roughly explain the structure of the representation
(\ref{eq:forrev-sde}). First note that the term on the right-hand side only
contains standard (unconditional) expectations, implying that the right-hand
side (unlike the left-hand side) is amenable to standard Monte
Carlo simulation  which is why we call (\ref{eq:forrev-sde}) a ``stochastic
representation''. The denominator of (\ref{eq:forrev-sde}) actually equals the transition density $p(0,x,T,y)$ of the solution $X$, and its presence
directly follows from the same term in the (analytical) definition of the
conditional expectation in terms of densities. In fact, it was precisely in
this context (i.e., in the context of density estimation) that Milstein,
Schoenmakers and Spokoiny introduced the general idea for the first time
\cite{Milstein2004}. 
In essence, the \emph{reverse} process $Y$ can be thought as
an ``adjoint'' process to $X$, as its infinitesimal generator is essentially
the adjoint operator of the infinitesimal generator of $X$ (see below for a
more detailed discussion in the SRN setting).

In a nutshell, the idea is that the law of the diffusion bridge admits a
Radon-Nikodym density with respect to the law of the concatenated process
$X^{(f)} \circ X^{(b)}$ with density given by $\mathcal{Y}$, \emph{provided}
that the trajectories meet at time $t^\ast$, i.e., provided that
$X^{(f)}(t^\ast) = X^{(b)}(t^\ast)$. Of course, this happens only with zero
probability\footnote{In the SRN setting, the probability is
  positive, since the state space is discrete.}, so we relax the above
equality with the help of a kernel with a positive bandwidth
$\epsilon$. Furthermore, note that by the independence of $X^{(f)}$ and
$X^{(b)}$, we can independently sample many trajectories of $X^{(f)}$ and many
trajectories of $X^{(b)}$ and then identify all pairs of trajectories
satisfying the approximate identity $X^{(f)}(t^\ast) \approx X^{(b)}(t^\ast)$
as determined by the kernel $\kappa_\epsilon$. This results in a Monte Carlo
algorithm, which, in principle, requires the calculation of a huge double sum
by summing over all pairs of $N$ samples from $X^{(f)}$ and $M$ samples from
$X^{(b)}$. A naive implementation of that algorithm would require a
prohibitive computational cost of order $O(M^2)$ operations, but
fortunately there are more efficient implementation relying on the structure
of the kernel and often reducing the complexity to $O(M \log(M))$
(see \cite{Bayer, BayerMC}). In this way, the
forward-reverse algorithm can nearly achieve the optimal Monte Carlo
convergence rate  of $1/2$. More precisely, assuming enough regularity on
the density of $X$ and assuming the use of a kernel of sufficiently high order
(depending on the dimension), the root-mean-squared error of the estimator is
$O(M^{-1/2})$ with a complexity $O(M\log(M))$ and a
bandwidth of $\epsilon = O(M^{-1/d})$. These statements
assume that we can exactly solve the SDEs driving the forward and the
reverse processes. Otherwise, the error induced by, say, the Euler scheme,
will be added.

The structure of the construction of the forward-reverse representation
(\ref{eq:forrev-sde}) and later of the corresponding Monte Carlo estimator in
\cite{Bayer} strongly suggests that the forward-reverse approach does not
rely on the continuity of diffusion processes, but merely on the Markov
property. Hence, the approach was generalized to discrete time Markov chains
in \cite{BayerMC} and is generalized to the case of continuous time
Markov chains with discrete state space in the this work. 

For a literature review on computational algorithms for computing conditional
expectations of functionals of diffusion processes we refer to \cite{Bayer}.
% flatex input end: [introduction.tex]

%}

%\clearpage
\section{Expectations SRN-Bridge Functionals}
% flatex input: [ReversedPJP.tex]
\label{sec:forwardreverse}
In this section, we derive the dynamics of the reverse paths and the expectation formula for  SRN-brige functionals.
The derivation follows the same scheme used in \cite{Milstein2004} , that is, 
i) write the master equation, ii) manipulate the master equation to obtain a backward Kolmogorov equation and, iii) derive the infinitesimal generator of the reverse process. 

\subsection{The Master Equation}
Let $X$ be a SRN defined by the intensity-reaction pairs $\seqof{(\nu_j,
  a_j(x))}{j=1}{J}$. Let $p(t,x,s,y)$ be its transition probability function,
\ie, $p(t,x,s,y){:=}\prob{X(s){=}y\SEP X(t){=}x}$ where $x,y\in\latt$ and
$0{<}t{<}s{<}T$.  The function $p$ satisfies the following linear system of
ODEs known as the master equation \cite{Gardiner,Risken,Kampen}:
%\begin{align}
%\partial_s p(t,x,s,y) &= \sum_{j=1}^J \LP a_j(y-\nu_j)p(t,x,s,y-\nu_j) - a_j(y)p(t,x,s,y)\RP\\
%p(t,x,t,y) &= \delta_{x=y} 
%\end{align}
\begin{align}\label{eq:ME}
\left\{ 
\begin{array}{rl}
\partial_s p(t,x,s,y) &= \sum_{j=1}^J \LP a_j(y-\nu_j)p(t,x,s,y-\nu_j) - a_j(y)p(t,x,s,y)\RP\COMMA\\
p(t,x,t,y) &= \delta_{x=y}\COMMA
\end{array}
\right. 
\end{align}
where $\delta$ is the Kronecker delta function.

A general analytic solution of \eqref{eq:ME} is in general computationally infeasible. 
Even numerical solutions are infeasible for systems with infinite or large number of states. 
For continuous state spaces, \eqref{eq:ME} becomes a parabolic PDE known as the Fokker-Planck Equation.
Next, we derive the generator of the reverse process in the SRN setting.

\subsection{Derivation of the Reverse Process}\label{sec:reverse}
Let us consider a fixed time interval $[t,T]$. For $s\in[t,T]$ and
$x,y\in\latt$, let us define $v(s,y):= \sum_x g(x) p(t,x,s,y)$ { provided
  that the sum converges. We remark here that $v$ cannot in general be
  interpreted as an expectation of $g$. Indeed, while $\sum_y p(t,x,s,y) = 1$,
  the sum over {$x$} could, in principle, even diverge. Hence, it is not
  a priori clear that $v$ admits a stochastic representation. However, }
multiplying both sides of the master equation \eqref{eq:ME} by $g(x)$ and
summing over $x$, we obtain:
\begin{align}\label{eq:sumoverx}
\left\{ 
\begin{array}{rl}
\partial_s v(s,y) &= \sum_{j=1}^J \LP a_j(y-\nu_j)v(s,y-\nu_j) - a_j(y) v(s,y)\RP\COMMA\\
v(t,y) &= g(y)\PERIOD
\end{array}
\right. 
\end{align} 

Now, let us consider a time reversal induced by a change of variables 
$\tilde s = T+t-s$  with  $\tilde v(\tilde s, y) := v(T+t-\tilde s,y) = v(s,y)$ leading to the following backward equation:
\begin{align}\label{eq:sumoverxback}
\left\{ 
\begin{array}{rl}
-\partial_{\tilde s} \tilde v(\tilde s, y) &= \sum_{j=1}^J \LP a_j(y-\nu_j) \tilde v (\tilde s, y-\nu_j) - a_j(y) \tilde v(\tilde s, y) \RP, \,\,t<\tilde s <T \COMMA\\
\tilde v(T,y) &= v(t,y) = g(y) \PERIOD
\end{array}
\right.
\end{align} 

Let $\tilde \nu_j := - \nu_j$. By adding and subtracting the term $a_j(y+\tilde \nu_j)\tilde v(\tilde s, y)$, we can write the first equation of \eqref{eq:sumoverxback} as
\begin{align*}
\partial_{\tilde s}\tilde v(\tilde s, y) + \sum_{j=1}^J \LP a_j(y+\tilde \nu_j)\LP \tilde v(\tilde s, y +\tilde \nu_j) -  \tilde v(\tilde s, y)\RP + \LP a_j(y+\tilde \nu_j)- a_j(y)\RP \tilde v(\tilde s, y) \RP=0\PERIOD
\end{align*}
As a consequence, the system \eqref{eq:sumoverxback} can be written as
\begin{align}\label{eq:kbeprevious}
\left\{ 
\begin{array}{ll}
\partial_{\tilde s} \tilde v(\tilde s, y) + \sum_{j=1}^J a_j(y+\tilde \nu_j)\LP \tilde v(\tilde s, y+\tilde \nu_j) - \tilde v(\tilde s, y)  \RP + c(y) \tilde v(\tilde s, y) = 0\COMMA\\
\tilde v(T,y) = g(y)\COMMA
\end{array}
\right. 
\end{align} 
where $c(y):= \sum_{j=1}^J a_j(y+\tilde \nu_j) {-} a_j(y)$.

Let us now define $\tilde a_j(y) := a_j(y+\tilde \nu_j)$ and substitute it into \eqref{eq:kbeprevious}. We have arrived at
the following backward Kolmogorov equation \cite{RogersWilliams} for the cost-to-go function $v(\tilde s, y)$, 
\begin{equation}\label{eq:reverse}
\left\{ 
\begin{array}{ll}
\partial_{\tilde s} \tilde v(\tilde s, y) + \sum_{j=1}^J \tilde a_j(y)\LP \tilde v(\tilde s, y+\tilde \nu_j) - \tilde v(\tilde s, y)  \RP + c(y) \tilde v(\tilde s, y) = 0\COMMA\\
\tilde v(T,y) = g(y) \PERIOD
\end{array}
\right.
\end{equation}
 
We recognize in \eqref{eq:reverse} the generator 
$\mathcal{L}_Y (\tilde v)(\tilde s, y):= \sum_{j=1}^J \tilde a_j(y)\LP \tilde v(\tilde s, y+\tilde \nu_j) - \tilde v(\tilde s, y)\RP$ that defines the so-called reverse process $Y\equiv \{Y(\tilde s,\omega)\}_{t \leq \tilde s \leq T}$ by 
\begin{align}
\label{eq:reverse-dynamics}
\prob{Y(\tilde s+d\tilde s) = y+\tilde \nu_j \SEP Y(\tilde s)=y} = \tilde a_j(y) d\tilde s
\end{align}
or equivalently by, 
\begin{align}
\prob{Y(\tilde s+d\tilde s) = y-\nu_j \SEP Y(\tilde s)=y} =  a_j(y-\nu_j) d\tilde s\PERIOD
\end{align}

The Feynman-Kac formula \cite{RogersWilliams} provides a stochastic representation of the solution of \eqref{eq:reverse},
\begin{align}
\tilde v(\tilde s, y) = \expt{g(Y(T)) \exp\LP \int_{\tilde s}^T c(Y(s))ds\RP\SEP Y(\tilde s)= y}\PERIOD
\end{align}
Notice that $Y$ is a SRN in its own right. 
We note in passing that stochastic representations based on shifted evaluations of the propensities have been derived independently in \cite{kt,kkst} to estimate 
variations  and differences of the cost to go function.

\subsection{The Forward-Reverse Formula for SRN}\label{sec:SRNformula}
%\red{I have taken the liberty to re-formulate this subsection. Please note
%  that the original formula (2.9) was incorrect, as the additive structure
%  $\Phi(X^{(f)},[s,t^\ast]) + \Phi(X^{(b)},[t^\ast,T])$ is only correct if the
%functional $\Phi$ is additive in terms of segments of the path --- as is the
%case in \emph{our} EM algorithm.}
Let us consider a time interval $[s,t]$ and assume that we only observe the
process $X$ on the end points, i.e., that we have $X(s) = x$ and $X(t) = y$
for some observed values $x,y \in \zset^d_+$. Fix an intermediate time 
$s {<}t^\ast {<} t$, which will be considered a numerical input parameter later
on. Denote by $X^{(f)}$ the process $X$ conditioned on \emph{starting} at
$X^{(f)}(s) = x$ and restricted to the time domain $[s, t^\ast]$.

Furthermore, let $Y$ denote the reverse process constructed
in~(\ref{eq:reverse-dynamics}) on the time domain $[t^\ast, t]$ (i.e.,
inserting $t^\ast$ for $t$ and $t$ for $T$ in the above subsection) started at
$Y(t^\ast) = y$. As noted above, $Y$ is again an SRN with reaction channels
$\seqof{(-\nu_j,\tilde a_j)}{j{=}1}{J}$. For convenience, we also introduce
the notation $X^{(b)}$ for the process $Y$ run backward in time, i.e., we
define $X^{(b)}(u){:=}Y(t^\ast{+}t{-}u)$ for $u{\in}[t^*,t]$, and notice that
$X^{(b)}(t) = y$. 

Recall that we aim to provide a \emph{stochastic representation},
\ie, a representation containing standard expectations only, for
conditional expectations of the form,
\begin{equation}
  \label{eq:H-general}
  \mathcal{H}(x,y) \equiv \expt{\left. \Phi\left(X, [s,t] \right) \, \right| \, X(s) =
    x,\, X(t) = y}\COMMA
\end{equation}
for $\Phi$ mapping $\zset^d_+$-valued paths to real numbers. Obviously, $\Phi$
needs to be integrable in order for $\mathcal{H}$ to be well defined, and we shall also
assume polynomial growth conditions on $\Phi$ and its derivatives %\blue{(we assume here that it is possible to extend $\Phi$ to the real domain in a sensible manner)} 
with respect
to jump-times of the underlying path. Moreover, we assume that $p(s,x,t,y) >
0$. 
%\red{How detailed should we be here?} 
Once again, the fundamental idea of
the forward-reverse algorithm of Bayer and Schoenmakers \cite{Bayer} is to
simulate trajectories of $X^{(f)}$ and (independently) of $X^{(b)}$ and then
look for any pairs that are ``linked''. Since the state space is now
discrete, we may, in principle, require exact linkage in the sense that we may
only consider pairs such that $X^{(f)}(t^\ast) = X^{(b)}(t^\ast)$. However, in
order to decrease the variance of the estimator, it may once again be
advantageous to relax this condition by introducing a \emph{kernel}.

By a kernel, we understand a function $\kappa: \zset^d \to \rset$ satisfying
\begin{equation*}
  \sum_{x \in \zset^d} \kappa(x) = 1.
\end{equation*}
Moreover, we call $\kappa$  a kernel of order $r \ge 0$ if, in addition,
\begin{equation*}
  \sum_{x \in \zset^d} x^\alpha \kappa(x) = 0
\end{equation*}
for any multi-index $\alpha$ with $1\le|\alpha| \le r$,
 $\alpha := \alpha_1+\cdots+\alpha_d$, and $x^{\alpha}:= x_1^{\alpha_1}\cdots x_d^{\alpha_d},\, \alpha\in \{0,1,2,\ldots\}$. For instance, any
non-negative symmetric kernel has order $r=1$ in this sense.
%\red{I am not sure we need or want that in the end.}

Having fixed one such kernel $\kappa$, we define a whole family of kernels
$\kappa_\epsilon$, indexed by the \emph{bandwidth} $\epsilon \ge 0$, by
\begin{equation*}
  \kappa_\epsilon(x) = C_\epsilon \kappa\left( \frac{x}{\epsilon} \right)
\end{equation*}
with the constant $C_\epsilon$ being defined by the normalization condition
$\sum_{x \in \zset^d} \kappa_\epsilon(x) = 1$. {Here, we implicitly assume the
kernel, $\kappa$, to be extended to $\rset^d$, for instance in a piecewise constant way.}
 As we necessarily have
$\kappa(x) \to 0$ as $|x| \to \infty$, it is easy to see that we have the
special case
\begin{equation*}
  \kappa_0(x) =
  \begin{cases}
    1, & x = 0,\\
    0, & x \neq 0.
  \end{cases}
\end{equation*}

\begin{rem}
  \label{rem:kernel-bandwidth}
  % Note that $\epsilon \mapsto C_\epsilon$ is piecewise constant. Moreover, the
  % same is true for the function-valued map $\epsilon \mapsto
  % \kappa_\epsilon$. In particular, 
  The Kronecker kernel $\kappa_0$ can also be realized as $\kappa_0 = \kappa_{\epsilon_0}$ for some
  $\epsilon_0 > 0$, which will depend on the base kernel $\kappa$, provided
  that the base kernel $\kappa$ has finite support.
\end{rem}

\begin{theorem}
  \label{thr:representation}
  Let $\Phi$ be a continuous real-valued functional  on the space of piecewise
  constant functions defined on $[s,t]$ and taking values in $\zset^d$
  (w.r.t.~uniform topology) such that
  both $\mathcal{H}$ and the right hand side of~(\ref{eq:themain}) is finite for any
  $\epsilon$. 
  %\red{Precise conditions?}
  With $\kappa_\epsilon$, $X^{(f)}$ and $X^{(b)}$ as above, we have
  \begin{equation}
    \label{eq:themain}
    \mathcal{H}(x,y) = \lim_{\epsilon \to 0} \frac{\expt{ \Phi\left( X^{(f)} \circ
          X^{(b)}, [s,t] \right) \kappa_\epsilon(X^{(f)}(t^\ast) -
        X^{(b)}(t^\ast)) \Psi\left( X^{(b)}, [t^\ast,t] \right)
      }}{\expt{ \kappa_\epsilon(X^{(f)}(t^\ast) -
        X^{(b)}(t^\ast)) \Psi\left( X^{(b)}, [t^\ast,t] \right) }}\COMMA
  \end{equation}
  where $X^{(f)} \circ X^{(b)}$ denotes the \emph{concatenation} of the paths
  $X^{(f)}$ and $X^{(b)}$ in the sense defined by
  \begin{equation*}
    X^{(f)} \circ X^{(b)} (u) \equiv
    \begin{cases}
      X^{(f)}(u), & s \le u \le t^\ast, \\
      X^{(b)}(u), & t^\ast < u \le t\COMMA
    \end{cases}
  \end{equation*}
  and
  \begin{equation*}
    \Psi(Z, [a,b]) {:=} \exp\left( \int_a^b c\left( Z(u) \right) du\right).
  \end{equation*}
\end{theorem}

\begin{rem}
  In line with Remark~\ref{rem:kernel-bandwidth}, we note that we could easily
  have avoided taking limits in Theorem~\ref{thr:representation} by replacing
  $\kappa_\epsilon$ with $\kappa_0$ everywhere in~(\ref{eq:themain}). 
  At this stage we note that the Monte Carlo estimator based on~(\ref{eq:themain})
  with positive $\epsilon$ will have considerable smaller variance than the
  version with $\epsilon = 0$, potentially outweighing the increased bias.
\end{rem}

\begin{proof}[Sketch of proof of Theorem~\ref{thr:representation}]
%  {\color{green}
    For simplicity, we assume that the kernel $\kappa$ has finite support and     that the functional $\Phi$ is uniformly bounded.
%  }
  We will prove convergence of the numerator and the denominator
  in~(\ref{eq:themain}) separately. Let us, hence, prove the more general
  case first, i.e., the convergence
  \begin{multline}
    \label{eq:auxiliary-rep}
    \mathrm{h}(x,y) {:=} \mathcal{H}(x,y) \,p(s,x,t,y) = \\
    \lim_{\epsilon\to0} \expt{ \Phi\left(
        X^{(f)} \circ X^{(b)}, [s,t] \right) \kappa_\epsilon(X^{(f)}(t^\ast) -
      X^{(b)}(t^\ast)) \Psi\left( X^{(b)}, [t^\ast,t] \right) }.
  \end{multline}
  
  In the first step, we assume that $\Phi(Z, [s,t])$ only depends on the
  values of $Z$ on a fixed grid, say $s = t_0 < t_1 < \cdots < t_n = t$, i.e.,
  \begin{equation*}
    \Phi(Z, [s,t]) = g\left(Z(t_0), \ldots, Z(t_n) \right).
  \end{equation*}
  Then~(\ref{eq:auxiliary-rep}) is proved (with minor modifications) in
  {\cite{Bayer} (Theorem 3.4)}. Indeed, a closer look at that proof reveals
  that only Markovianity of $X$ is really used.

  Furthermore, note that any continuous functional $\Phi$ can be approximated
  by functionals $\Phi_n$ depending only on the values of the process on a
  (ever finer) finite grid $t_0, \ldots, t_n$. As, on the one side,
  \begin{multline*}
    \mathrm{h}(x,y) = \expt{\left. \Phi\left(X, [s,t] \right) \, \right| \, X(s) =
    x,\, X(t) = y }p(s,x,t,y) = \\
  \lim_{n\to\infty}
  \expt{\left. \Phi_n\left(X, [s,t] \right) \, \right| \, X(s) = x,\, X(t) =
    y} p(s,x,t,y)
  \end{multline*}
  and, on the other side,
  \begin{multline*}
    \lim_{\epsilon\to0} \lim_{n\to\infty} \expt{ \Phi_n\left(
        X^{(f)} \circ X^{(b)}, [s,t] \right) \kappa_\epsilon(X^{(f)}(t^\ast) -
      X^{(b)}(t^\ast)) \Psi\left( X^{(b)}, [t^\ast,t] \right) } =\\
    \lim_{\epsilon\to0} \expt{ \Phi\left( 
        X^{(f)} \circ X^{(b)}, [s,t] \right) \kappa_\epsilon(X^{(f)}(t^\ast) -
      X^{(b)}(t^\ast)) \Psi\left( X^{(b)}, [t^\ast,t] \right) }\PERIOD
  \end{multline*}
  We are left to prove that
  \begin{multline*}
    \lim_{\epsilon\to0} \lim_{n\to\infty} \expt{ \Phi_n\left(
        X^{(f)} \circ X^{(b)}, [s,t] \right) \kappa_\epsilon(X^{(f)}(t^\ast) -
      X^{(b)}(t^\ast)) \Psi\left( X^{(b)}, [t^\ast,t] \right) } =\\  
    \lim_{n\to\infty} \lim_{\epsilon\to0}  \expt{ \Phi_n\left(
        X^{(f)} \circ X^{(b)}, [s,t] \right) \kappa_\epsilon(X^{(f)}(t^\ast) -
      X^{(b)}(t^\ast)) \Psi\left( X^{(b)}, [t^\ast,t] \right) },
  \end{multline*}
  % which is trivial in the case when the kernel $\kappa$ has finite
  % support. 
  %\red{Full proof if wanted?}
  which follows as $\kappa_0 = \kappa_{\epsilon_0}$ for some $\epsilon_0 >
  0$. In fact, it even follows in the general case by
    dominated convergence.
  %{\color{green} (In fact, it even follows in the general case by    dominated convergence.)} 

  Finally, the proof of convergence of the denominator is a special case of
  the proof for the numerator, and therefore, the convergence of the fraction follows from the
  continuity of $(a,b) \mapsto a/b$ for $b > 0$.
\end{proof}

\section{The EM Algorithm for SRN}
% flatex input: [theEM.tex]
\label{sec:EM}
In this section, we present the EM algorithm for SRN, which is the main step for computing the parameter estimation. First, we explain  the EM algorithm in general, and then, we derive the log-likelihood function for a fixed realization of the process, $X$. Finally, we present the EM algorithm for SRN.

\subsection{The EM Algorithm}\label{met:EM}
{
The EM algorithm \cite{Dempster77,Casella, WatanabeYamaguchi, McLachlanEM}
its named from its two steps: expectation and maximization.
It is an iterative algorithm that, given an initial guess and a stopping rule, provides an approximation for a local maximum or saddle point of the likelihood function, $\text{lik}(\theta \SEP \mathcal{D})$.
It is a data augmentation technique in the sense that the maximization of the likelihood $\text{lik}(\theta \SEP \mathcal{D})$ is performed by treating the data $\mathcal{D}$ as a part of a larger data set, $(\mathcal{D},\tilde {\mathcal{D}})$, where the complete-likelihood, $\text{lik}^c(\theta \SEP \mathcal{D},\tilde {\mathcal{D}})$, is amenable to maximization.
Given an initial guess $\theta^{(0)}$, the EM algorithm maps $\theta^{(p)}$ into $\theta^{(p+1)}$ by the
\bigskip
\begin{enumerate}
\item expectation step: $Q_{\theta^{(p)}}(\theta \SEP\mathcal{D}) := \mathrm{E}_{\theta^{(p)}}\left[{\log(\text{lik}^c(\theta\SEP \mathcal{D},\tilde {\mathcal{D}}))\SEP \mathcal{D}}\right]$, and the 
\item maximization step: $\theta^{(p+1)} := \arg \max_{\theta} Q_{ \theta^{(p)}}(\theta \SEP\mathcal{D})$. 
\end{enumerate} 
\bigskip
Here, $\mathrm{E}_{\theta^{(p)}}\left[ \cdot \SEP \mathcal{D}\right]$, denotes the expectation associated with the distribution of $\tilde {\mathcal{D}}$ under the parameter choice  
$\theta^{(p)}$, conditional on the data, $\mathcal{D}$.
In many applications, the expectation step is computationally infeasible and $Q_{\theta^{(p)}}(\theta \SEP \mathcal{D})$ should be approximated by some estimate, 
\begin{align*}
\hat Q_{\theta^{(p)}}(\theta \SEP \mathcal{D}) &:= \hat{\mathrm{E}}_{\theta^{(p)}} \left[\log(\text{lik}^c(\theta\SEP \mathcal{D},\tilde {\mathcal{D}}))\SEP \mathcal{D}\right]
\PERIOD
\end{align*}

\begin{rem}[The Monte Carlo EM]\label{rem:MCEM}
If we know how to sample a sequence of $M$ independent variates $\seqof{\tilde {\mathcal{D}}_i}{i=1}{M} \sim \tilde {\mathcal{D}} \SEP \mathcal{D}$, with parameter $\theta^{(p)}$, then we can define the following Monte Carlo estimator of $Q_{\theta^{(p)}}(\theta \SEP\mathcal{D})$, 
\end{rem}
\begin{align*}
\hat Q_{\theta^{(p)}}(\theta \SEP\mathcal{D}) &:= 
\frac{1}{M} \sum_{i=1}^M 
\log(\text{lik}^c(\theta\SEP \mathcal{D},\tilde {\mathcal{D}}_i))%\SEP \theta^{(p)},\mathcal{D}
%\avg{\log(\text{lik}^c(\theta\SEP \mathcal{D},\tilde {\mathcal{D}}))\SEP \theta^{(p)},\mathcal{D}}{M}
\PERIOD
\end{align*}

In Section \ref{sec:FREM}, we describe how to simulate exact and approximate samples of $\tilde {\mathcal{D}} \SEP \mathcal{D}$.
}

\subsection{The Log-Likelihood Function for Continuously Observed Paths}\label{sec:contobservedpaths}
The goal of this section is to derive an expression for the likelihood of a particular path, $(X(t,\omega_0))_{t\in[0,T]}$, of the process $X$, where $\omega_0 \in \Omega$ is a fixed realization. 
An important assumption in this work is that the propensity functions $a_j$
can be written as $a_j(x) = c_j g_j(x)$ for $j{=}1,\ldots,J$ and
$x\in\latt$ {where $g_j$ are known functionals and $c_j$ are
  considered the unknown parameters.} Define $\theta{:=}(c_1,\ldots,c_J)$.
%This is known as mass action kinetics \cite{Holmes}.
Let us denote the jump times of $(X(t,\omega_0))_{t\in[0,T]}$ in $(0,T)$ 
by $\xi_1, \xi_2, \ldots, \xi_{N-1}$. 
Define $\xi_0:=0$, $\xi_N:=T$ and $\Delta \xi_i = \xi_{i+1} - \xi_{i}$ for $i=0,1,\ldots,N-1$.

Let us assume that the system is in the state $x_0$ at time $0$.
We have that $\xi_1$ is the time to the first reaction, or equivalently, the time that the system spend at $x_0$ (sojourn time or holding time at state $x_0$). 
Let us denote by $\nu_{\xi_1}$ the reaction that takes place at $\xi_1$, and therefore, the system at time $\xi_1$ is in the state $x_1:= x_0 + \nu_{\xi_1}$. 
From the SSA algorithm, it is easy to see that the probability density corresponding to this transition is the product $a_{\nu_{\xi_1}}(x_0) \exp{(-a_0(x_0) \Delta \xi_0)}$. 

By the Markov property we can see that the density of one path $\seqof{(\xi_i,x_i)}{i=0}{N-1}$ is given by 
\begin{equation}\label{likpath}
\prod_{i=1}^{N-1} a_{\nu_{\xi_i}}(x_{i-1})\exp{(-a_0(x_{i-1})\Delta \xi_{i-1})}\times \exp{(-a_0(x_{N-1})\Delta \xi_{N-1})}\PERIOD
\end{equation}
The last factor in \eqref{likpath} is due to the fact that we know that the system will remain in the state $x_{N-1}$ in the time interval $[\xi_{N-1},T)$.

Rearranging the factors in \eqref{likpath}, we obtain 
\begin{equation}\label{eprod}
\exp{\LP -\sum_{i=0}^{N-1} a_0(x_i) \Delta \xi_{i}\RP} \prod_{i=1}^{N-1} a_{\nu_{\xi_i}}(x_{i-1})\PERIOD
\end{equation}
Now, taking logarithms in \eqref{eprod} we have
\begin{equation*}
 -\sum_{i=0}^{N-1} a_0(x_i) \Delta \xi_{i}  +  \sum_{i=1}^{N-1} \log(a_{\nu_{\xi_i}}(x_{i-1}))\COMMA
\end{equation*}
which by the definition of $a_0$ can be written as
\begin{equation*}
 -\sum_{i=0}^{N-1} \sum_{j=1}^J a_j(x_i) \Delta\xi_{i}  +  \sum_{i=1}^{N-1} \log(c_{\nu_{\xi_i}}g_{\nu_{\xi_i}}(x_{i-1}))\PERIOD
\end{equation*}
Interchanging the order in the summation and denoting the number of times that the reaction $\nu_j$ occurred in the interval $[0,T]$ by $R_{j,[0,T]}$, we have
\begin{equation}\label{lastuseless}
\sum_{j=1}^J\LP -c_j \sum_{i=0}^{N-1} g_j(x_i) \Delta \xi_{i} + \log(c_j) R_{j,[0,T]} \RP + \sum_{i=1}^{N-1} \log(g_{\nu_{\xi_i}}(x_{i-1}))\PERIOD 
\end{equation}
Observing that the last term in \eqref{lastuseless} does not depend on
$\theta$, the complete log-likelihood of the path $(X(t,\omega_0))_{t\in[0,T]}$
is {up to constant terms} given by
\begin{equation}
\ell^c(\theta) := \sum_{j=1}^J  \log(c_j) R_{j,[0,T]} - c_j F_{j,[0,T]},\, \text{ with } \theta {=} (c_1,\ldots,c_J)\COMMA
\end{equation}
where $F_{j,[0,T]}:= g_j(x_0)\Delta \xi_0 +\cdots+ g_j(x_{N-1})\Delta
\xi_{N-1} = \int_0^T
g_j(X(s))\,ds$. {The last equality is due to $g_j$ being piecewise constant in the partition $\{\xi_0,\xi_1,\ldots,\xi_N\}$.}  

Now let us assume that we have a collection of intervals, $\seqof{I_k =[s_k,t_k]}{k=1}{K}\subset [0,T]$, where we have continuously observed the process $(X(t,\cdot))_{t\in I_k}$ at each $I_k$. 
We define the log-likelihood function as:
\begin{equation*}
\ell^c(\theta) := \sum_{j=1}^J \LP \log(c_j)\sum_{k=1}^K R_{j,I_k} - c_j \sum_{k=1}^K F_{j,I_k}\RP\PERIOD
\end{equation*}
%\todo{Why $\omega_k$ instead of just $\omega$? AM: this is because we do not care about paths once we organize the data in intervals. If we put just $\omega$ it may seem like there is only one observed path, and we could have many, or even different pieces of different paths.}

\begin{rem} \label{rem:nonrandom} 
%Note that that $R_{j,I_k}(\omega)$ and $F_{j,I_k}(\omega)$ are particular realizations of the random variables,
%$R_{j,I_k}$ and $F_{j,I_k}$, respectively.
Note that $R_{j,I_k}$ and $F_{j,I_k}$ are random variables, which are functions of the full paths of $X$ but not of the discretely observed paths. Hence, they are random given the data $\mathcal{D}$ as defined in \eqref{def:data}.
\end{rem}
\vspace{1mm}

\subsection{The EM Algorithm for SRNs}

%Let us assume that we have observed the process $X$ only at the extremes of the intervals $I_k$. 
%In such case, our data, $\mathcal{D}$, is the collection  
%$\seqof{I_k{=}[s_k,t_k],X(s_k),X(t_k)}{k=1}{K}$. 
%Consider the data, $\mathcal{D}$, as defined in \eqref{def:data}.

According to the Section \ref{met:EM}, for a particular value of the parameter $\theta$, say $\theta^{(p)}$, we define 
\begin{align*}\label{loglikc}
Q_{\theta^{(p)}}(c_1,\ldots,c_J \SEP \mathcal{D}) := \sum_{j=1}^J \LP \log(c_j)\sum_{k=1}^K 
\mathrm{E}_{\theta^{(p)}}\left[{R_{j,I_k}\SEP \mathcal{D}}\right] - c_j \sum_{k=1}^K \mathrm{E}_{\theta^{(p)}}\left[{F_{j,I_k}\SEP \mathcal{D}}\right]\RP\COMMA
\end{align*}
%\todo{Some parenthesis are wrong here!}
where 
$\mathrm{E}_{\theta^{(p)}}\left[{R_{j,I_k}\SEP \mathcal{D}}\right] =\mathrm{E}_{\theta^{(p)}}\left[{R_{j,I_k}\SEP X(s_k){=}x(s_k),X(t_k){=}x(t_k)}\right]$ (by the Markov property), and analogously for $F_{j,I_k}$.

Now consider the partial derivatives of $Q_{\theta^{(p)}}(c_1,\ldots,c_J \SEP \mathcal{D})$ with respect to $c_j$ 
\begin{align*}
\partial_{c_j} Q_{\theta^{(p)}}(c_1,\ldots,c_J \SEP \mathcal{D})=  
\frac{1}{c_j}\sum_{k=1}^K \mathrm{E}_{\theta^{(p)}}\left[{R_{j,I_k}\SEP \mathcal{D}}\right]
  -\sum_{k=1}^K \mathrm{E}_{\theta^{(p)}}\left[{F_{j,I_k}\SEP\mathcal{D}}\right]\PERIOD
\end{align*}
Therefore, $\nabla Q_{\theta^{(p)}}(c_1,\ldots,c_J \SEP\mathcal{D}) = 0$ is obtained at $\theta^*= \LP c^*_1,\ldots,c^*_J \RP$ such that 
\begin{equation}
c^*_j = \frac{\sum_{k=1}^K \mathrm{E}_{\theta^{(p)}}\left[{ R_{j,I_k}\SEP
      \mathcal{D}}\right]}{\sum_{k=1}^K
  \mathrm{E}_{\theta^{(p)}}\left[{F_{j,I_k}\SEP \mathcal{D}}\right]}, \ j{=}1,
\ldots, J.
\end{equation}
This is clearly the global maximization point of the function $Q_{\theta^{(p)}}(\cdot \SEP \mathcal{D})$.

%Now, the idea is to approximate $\expt{R_{j,I_k}\SEP \theta, \mathcal{D}}$ and $ \expt{F_{g_j,I_k}\SEP \theta,\mathcal{D}}$ by the Monte Carlo method. To this end, for each interval $I_k$, we should sample $M_k$ bridges of the process $X$ and compute $\avg{R_{j,I_k}\SEP \theta, \mathcal{D}}{M_k}$ and 
%$\avg{F_{j,I_k}\SEP \theta, \mathcal{D}}{M_k}$
%
%
%By approximating the expectations by empirical averages, we have
%\begin{equation}
%\tilde c^*_j = \frac{\sum_{k=1}^K\avg{ R_{j,I_k}\SEP \theta^*, \mathcal{D}}{M_k}}{\sum_{k=1}^K\avg{F_{j,I_k}\SEP \theta^*,\mathcal{D}}{M_k}}\COMMA
%\end{equation}
%where each $M_k$ has to be chosen sufficiently large such that the coefficient of variation of the estimates fall bellow a certain prescribed threshold. Observe that $ R_{j,I_k}$ and 
%$F_{j,I_k}$ are not independent random variables since they are computed using the same SRN-bridge.

The EM algorithm for this particular problem generates a {deterministic} sequence 
$\seqof{\theta^{(p)}}{p=1}{+\infty}$ that starts from a deterministic initial guess
$\theta^{(0)}$ provided by phase I (see Section \ref{sec:phaseI}) and evolves by 
\begin{equation}\label{eq:EMiteration1}
 c^{(p+1)}_j =  \frac{\sum_{k=1}^K \mathrm{E}_{\theta^{(p)}}\left[{ R_{j,I_k}\SEP \mathcal{D}}\right]}{\sum_{k=1}^K \mathrm{E}_{\theta^{(p)}}\left[{F_{j,I_k}\SEP \mathcal{D}}\right]}\COMMA
\end{equation}
where $\theta^{(p)} = \LP c_1^{(p)},\ldots, c_J^{(p)} \RP$.

% flatex input end: [theEM.tex]

%\clearpage

%\clearpage
\section{Forward-Reverse Monte Carlo EM Algorithm for SRNs}
% flatex input: [FREM.tex]
\label{sec:FREM}
In this section, we present a two-phase algorithm for estimating the parameter $\theta$. 
Phase I is deterministic while  phase II is stochastic.
We consider the data, $\mathcal{D}$, as given by \eqref{def:data}. The main goal of this section is to provide a Monte Carlo version of formula \eqref{eq:EMiteration1}.

\subsection{Phase I: using Approximating ODEs}\label{sec:phaseI}
The objective of phase I is to address the key problem of finding a suitable initial point ${\theta}^{(0)}_{I\!I}$ to reduce the variance (or the computational work) of phase II, 
thereby increasing (in some cases dramatically) the number of SRN bridges from the sampled forward-reverse trajectories for all time intervals.

Let us now describe phase I. From the user-selected seed, $\theta^{(0)}_{I}$,
we solve the following deterministic optimization problem using some appropriate numerical iterative method:
\begin{align}\label{eq:seedI}
{\theta}^{(0)}_{I\!I} := \operatorname*{arg\,min}_{\theta\geq 0} 
\sum_{k} w_k\, 
\norm{\tilde Z^{(f)}(t_k^*;\theta)- \tilde Z^{(b)}(t_k^*;\theta)}^2
%\text{ starting from } \theta^{(0)}_{I}
\PERIOD
\end{align} 
%\todo{Why not optimize sum over $k$, i.e., a single $\lambda_i$?. AM: We already tried this way, but we sometimes obtained bad seed compared with the current approach.}
Here, $\tilde Z^{(f)}$ is the ODE approximation defined by \eqref{eq:ODE} in the interval $[s_k,t_k^*]$, to the SRN defined by the reaction channels, $\seqof{(\nu_j,a_j)}{j=1}{J}$, and the initial condition $x(s_k)$; 
$\tilde Z^{(r)}$ is the ODE approximation in the interval $[t_k^*,t_k]$ to the SRN defined by the reaction channels, 
$\seqof{(-\nu_j,\tilde a_j)}{j=1}{J}$, and by the initial condition $x(t_k)$. 
Let us recall that in Section \ref{sec:reverse}, $\tilde a_j(x)$ was  defined as $a_j(x{-}\nu_j)$. 
We define $\tilde Z^{(b)}(u,\theta){:=}\tilde Z^{(r)}(t_k^*{+}t_k{-}u,\theta)$ for $u\in[t_k^*,t_k]$. Furthermore, $w_k {:=} (t_k{-}s_k)^{-1}$ and $\norm{\cdot}$ is the Euclidean norm in $\rset^d$.
The rationale behind this particular choice of the weight factors is based on the mitigation of the effect of very large time intervals, where the evolution of the process, $X$, may be more uncertain. A better (but more costly) measure would be the inverse of the maximal variance of the SRN bridge.

\begin{rem}[An alternative definition of ${\theta}^{(0)}_{I\!I}$]\label{rem:alternative} In some cases, convergence issues arise when solving the problem \eqref{eq:seedI}. We found it useful to solve a set of simpler problems whose answers can be combined to provide a reasonable seed for phase II: 
more precisely, we solve $K$ deterministic optimization problems, one for each time interval $[s_k,t_k]$:
\begin{align*}%\label{eq:optproblem}
\lambda_{k} := \operatorname*{arg\,min}_{\theta\geq 0} \norm{\tilde Z^{(f)}(t_k^*;\theta)-\tilde Z^{(b)}(t_k^*;\theta)}\COMMA
\end{align*} 
all of which were solved iteratively with the same seed, $\theta^{(0)}_{I}$. Then, we define
\begin{align}\label{eq:seedIalt}
{\theta}^{(0)}_{I\!I}:= \frac{\sum_k w_k \lambda_{k} }{\sum_k w_k}\PERIOD
\end{align}
\end{rem}

\subsection{Phase II: the Monte Carlo EM}
\begin{figure}[h!]
\centering
%\begin{minipage}{0.59\textwidth}
\begin{minipage}{0.49\textwidth}
\includegraphics[scale=0.40]{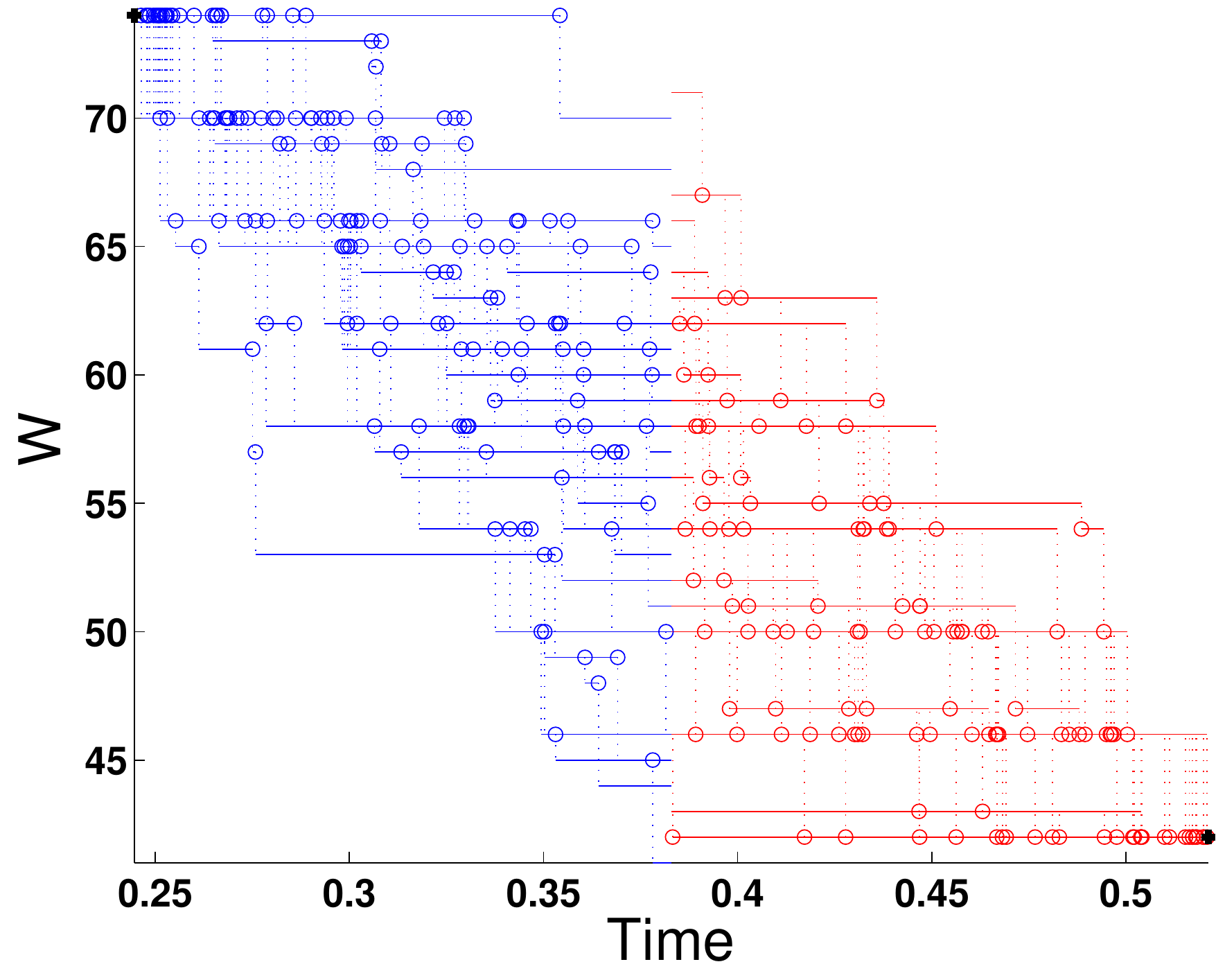}
%\myinclude{Figures/Wear_Cilindri_T1_FR_paths.pdf}
\end{minipage}
\hfill
%\begin{minipage}{0.4\textwidth}
\begin{minipage}{0.49\textwidth}
\includegraphics[scale=0.3]{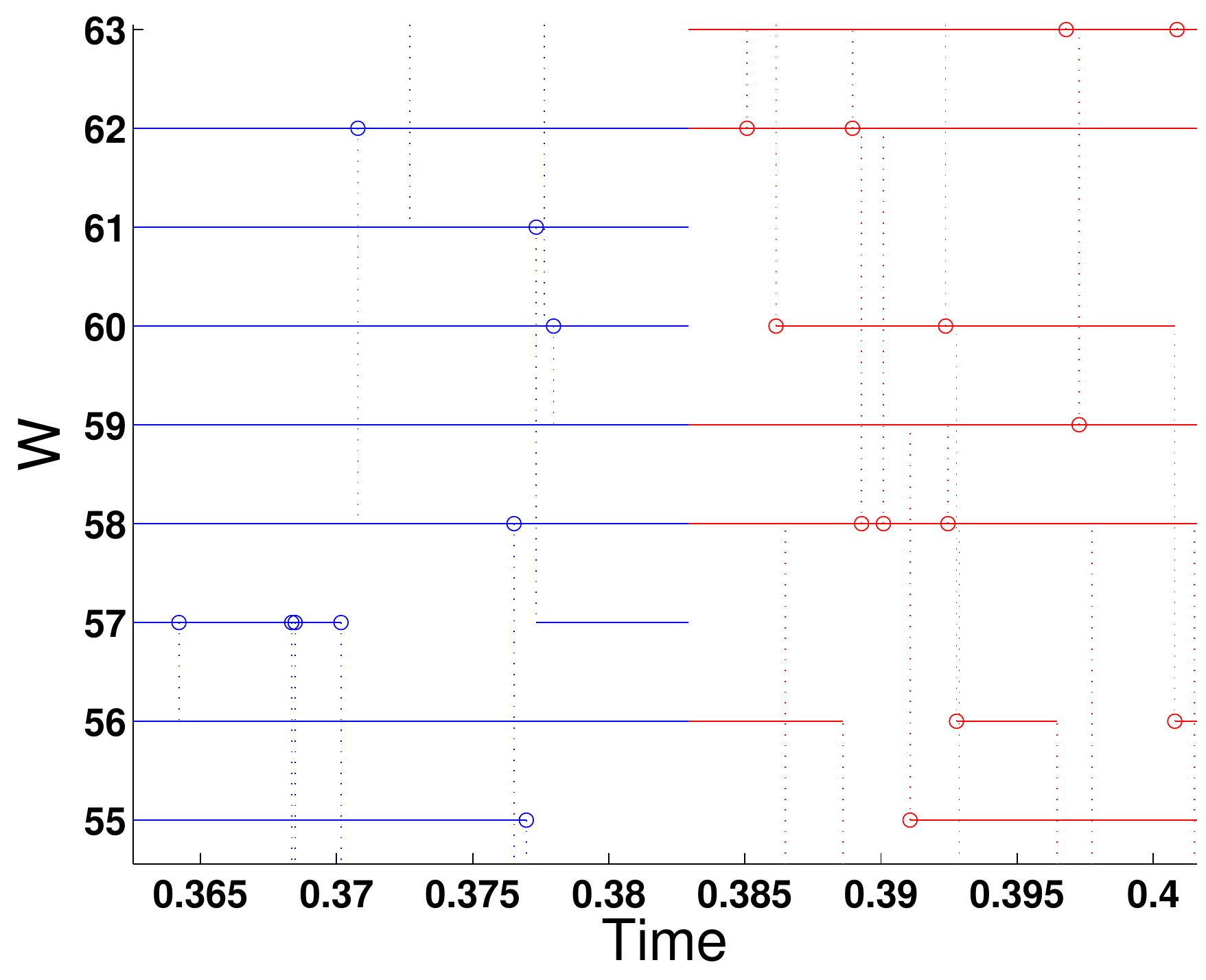}
%\myinclude{Figures/Wear_Cilindri_T1_FR_paths_zoom.pdf}
\end{minipage}
\caption{Left: Illustration of the forward-reverse path simulation in Phase II. The plot corresponds to a given interval for the wear data, presented in Section \ref{ex:wear}. The observed values are marked with a black circle (beginning and end of the interval). On the y-axis we plot the thickness process $X(t)$, derived from the wear process of the cylinder liner. Observe that every forward path that ends up at a certain value will be joined with every backward path that ends up at the same value when using the Kronecker kernel. For example, this happens at value 58, where several forward paths end and several backward paths start. Right: Zoom near value 58.}
\label{fig:frpaths}
\end{figure}

%This phase implements the Monte Carlo EM Algorithm for SRNs. 
%It starts with the simulation of SRN-bridges and the Monte Carlo approximation of the expected values given by \eqref{eq:themain}.
%We also address the use of Epanechnikov kernels.
In our statistical estimation approach, the Monte Carlo EM Algorithm uses data (pseudo-data) generated by those forward and backward simulated paths that result in exact or approximate SRN bridges. In Figure \ref{fig:frpaths}, we illustrate this idea for the wear example data presented in Section \ref{ex:wear}. Phase II implements the Monte Carlo EM algorithm for SRNs. 

\subsubsection{Simulating Forward and Backward Paths}\label{sec:clouds}
This phase starts with the simulation of forward and backward paths at each time interval $I_k$, for $k{=}1,...,K$. More specifically, given an estimation of the true parameter $\theta$, say, $\hat{\theta} = (\hat{c}_1, \hat{c}_2,\ldots, \hat{c}_J)$, the fist step is to simulate $M_k$ forward paths with reaction channels $\seqof{\nu_j,\hat{c}_j g_j(x)}{j=1}{J}$ in $[s_k,t_k^*]$, all of them starting at $s_k$ from $x(s_k)$ (see Section \ref{sec:Mk} for details about the selection of $M_k$).  Then, we  simulate $M_k$ backward paths with reaction channels $\seqof{-\nu_j,\hat{c}_j g_j(x-\nu_j)}{j=1}{J}$ in $[t_k^*, t_k]$, all starting at $t_k$ from $x(t_k)$.
%\begin{enumerate}
%\item Simulate $M_k$ forward paths with reaction channels $\seqof{\nu_j,\hat{c}_j g_j(x)}{j=1}{J}$ in $[s_k,t_k^*]$, starting at $s_k$ from $x(s_k)$. 
%\item Simulate $M_k$ backward paths with reaction channels $\seqof{-\nu_j,\hat{c}_j g_j(x-\nu_j)}{j=1}{J}$ in $[t_k^*, t_k]$, 
%starting at $t_k$ from $x(t_k)$.
%\end{enumerate}
%\vspace{1mm}
Let $\seqof{\tilde X^{(f)}(t_k^*,\tilde{\omega}_m)}{m=1}{M_k}$ and $\seqof{\tilde X^{(b)}(t_k^*,\tilde{\omega}_{m'})}{m'=1}{M_k}$ denote the values of the simulated forward and backward paths at the time $t_k^*$, respectively. If the intersection of these two sets of points is nonempty, then, there exists at least one $m$ and one $m'$ such that the forward and backward paths  can be linked as one SRN path that connects  $x(s_k)$ and $x(t_k)$ data values. 
%By construction, this SRN-bridge has the correct probability law and there is no need for any acceptance-rejection scheme.

When the number of simulated paths $M_k$ is large enough, and an appropriate guess of the parameter ${\theta}$ is used to generate those paths, then, due to the discrete nature of our state space $\latt$, we expect to generate a sufficiently large number of \emph{exact SRN bridges} to perform statistical inference.
However, at early stages of the Monte Carlo EM algorithm, our approximations to the unknown parameter ${\theta}$ are not expected to provide a large number of exact SRN bridges. 
In such a case, we can use kernels to relax  the notion of an exact SRN bridge (see Section \ref{sec:SRNformula}).
Notice that in the case of exact  SRN bridges, we are implicitly using a Kronecker  kernel in the formula
\eqref{eq:themain},
that is, $\kappa$ takes the value $1$ when $\tilde X^{(f)}(t_k^*,\tilde{\omega}_m) = \tilde X^{(b)}(t_k^*,\tilde{\omega}_{m'})$ and $0$ otherwise. 
We can relax this condition to obtain \emph{approximate SRN bridges}. 

To make computationally efficient use of kernels, we sometimes transform the endpoints of the forward and backward paths generated in the interval $I_k$, 
\begin{align}\label{eq:cloudX}
\mathcal{X}_k := (&
\tilde X^{(f)}(t_k^*,\tilde{\omega}_1),
\tilde X^{(f)}(t_k^*,\tilde{\omega}_2),\ldots,
\tilde X^{(f)}(t_k^*,\tilde{\omega}_{M_k}),\\ \nonumber
& 
\tilde X^{(b)}(t_k^*,\tilde{\omega}_{M_k+1}), 
\tilde X^{(b)}(t_k^*,\tilde{\omega}_{M_k+2}),\ldots, 
\tilde X^{(b)}(t_k^*,\tilde{\omega}_{2M_k})
)\COMMA
\end{align}

into 
\begin{align}\label{eq:cloudY}
H(\mathcal{X}_k) := (&
\tilde Y^{(f)}(t_k^*,\tilde{\omega}_1),
\tilde Y^{(f)}(t_k^*,\tilde{\omega}_2),\ldots,
\tilde Y^{(f)}(t_k^*,\tilde{\omega}_{M_k}),\\ \nonumber
& 
\tilde Y^{(b)}(t_k^*,\tilde{\omega}_{M_k+1}), 
\tilde Y^{(b)}(t_k^*,\tilde{\omega}_{M_k+2}),\ldots, 
\tilde Y^{(b)}(t_k^*,\tilde{\omega}_{2M_k})
)\COMMA
\end{align}
 by  a linear transformation $H$ with the aim of eliminating possibly high correlations in the  components of  $\mathcal{X}_k$. 
The original cloud of points $\mathcal{X}_k$ formed by  extremes of the forward and backward paths is then transformed into 
$H(\mathcal{X}_k)$, which hopefully has a covariance matrix close to a multiple of the $d$-dimensional identity matrix $\alpha I_d$.  Ideally, the coefficient $\alpha$ should be chosen in such way that each $d$-dimensional unitary cube centered at $\tilde Y^{(f)}(t_k^*,\tilde{\omega}_m)$ contains on average one element of $\cup_{m'}\{\tilde Y^{(b)}(t_k^*,\tilde{\omega}_{m'})\}$.
Note that this transformation changes (generally slightly) the variances of our estimators (see Section \ref{sec:trasf} for details about the selection of $\alpha$ and $H$). 
%\todo{For me, this is a bit too fast: if I understand correctly, the idea of %$H$ is to properly choose the bandwidth of the kernel. This needs to be %explained! AM: Absolutely. This is the objective of the section 5.1. I %planned to write down a careful explanations, using some plots.}

In our numerical examples, we use the  Epanechnikov kernel
\begin{equation}\label{eq:epa}
\kappa(\eta) := \LP \frac{3}{4} \RP^d \,\prod_{i=1}^{d} (1-\eta_i^2) \indicator{\abs{\eta_i}\leq1}\COMMA
\end{equation}
where  $\eta$ is defined as 
\begin{align}\label{eq:eta}
\eta \equiv \eta_k(m,m') := \tilde Y^{(f)}(t_k^*,\tilde{\omega}_m)-\tilde Y^{(b)}(t_k^*,\tilde{\omega}_{m'})\PERIOD
\end{align}

{This choice is motivated by the way in which we compute $\eta_k(m,m')$ avoiding whenever possible to make $M_k^2$ calculations. The support of $\kappa$ is perfectly adapted to our strategy of dividing $\rset^d$ into unitary cubes with vertices in $\zset^d$.}

%In Section \ref{sec:details}, we provide details about the linear transformation $H$, the value of $\alpha$.

\subsubsection{Kernel-weighted Averages for the Monte Carlo EM}
As we previously mentioned, the only available data in the interval $I_k$ correspond to the observed values of the process, $X$, at its extremes. 
Therefore, the expected values  $\mathrm{E}_{\theta^{(p)}}\left[{ R_{j,I_k}\SEP \mathcal{D}}\right]$  and $\mathrm{E}_{\theta^{(p)}}\left[{ F_{j,I_k}\SEP \mathcal{D}}\right]$ in the formula \eqref{eq:EMiteration1} must be approximated by SRN-bridge simulation.
{To this end, we generate a set of $M_k$ forward paths in the interval $I_k$ using $\hat{\theta}_{I\!I}^{(p)}$ as the current guess for the unknown parameter $\theta^{(p)}$. Having generated those paths, we record $R^{(f)}_{j,I_k}(\tilde{\omega}_m)$ and  $F^{(f)}_{j,I_k}(\tilde{\omega}_m)$ for all $j=1,2,\ldots,J$ and  $m=1,2,\ldots,M_k$
as defined in Section \ref{sec:contobservedpaths}.
Analogously, we record $R^{(b)}_{j,I_k}(\tilde{\omega}_{m'})$ and  $F^{(b)}_{j,I_k}(\tilde{\omega}_{m' })$ for all $j=1,2,\ldots,J$ and  $m'=1,2,\ldots,M_k$.}

%In this section, we motivate a kernel-based method for approximating 
%$\mathrm{E}_{\theta^{(p)}}\left[{ R_{j,I_k}\SEP \mathcal{D}}\right]$ and $\mathrm{E}_{\theta^{(p)}}\left[{ F_{j,I_k}\SEP \mathcal{D}}\right]$. In Section XX, we indicate how to sample exact and approximate SRN-bridges, where the approximate SRN-bridges are weighted by kernel, $\kappa(\eta)$. 

Consider the following $\kappa$-weighted averages, where $\kappa=\kappa_{\epsilon}$ for an appropriate choice of bandwidth $\epsilon$ that approximate  
$\mathrm{E}_{\theta^{(p)}}\left[{ R_{j,I_k}\SEP \mathcal{D}}\right]$ and 
$\mathrm{E}_{\theta^{(p)}}\left[{ F_{j,I_k}\SEP \mathcal{D}}\right]$, respectively:
\begin{align}\label{eq:averagesIk}
\avgsub{ R_{j,I_k}\SEP \mathcal{D}}{\kappa}{\hat{\theta}_{I\!I}^{(p)}}&:=\frac{\sum_{m,m'} \LP R^{(f)}_{j,I_k}(\tilde{\omega}_m) + R^{(b)}_{j,I_k}(\tilde{\omega}_{m'})\RP\kappa(\eta_k(m,m'))\psi_k(m')}{\sum_{m,m'} \kappa(\eta_k(m,m'))\psi_k(m')} \text{ and}\\
\nonumber \avgsub{ F_{j,I_k}\SEP \mathcal{D}}{\kappa}{\hat{\theta}_{I\!I}^{(p)}}&:=\frac{\sum_{m,m'} \LP F^{(f)}_{j,I_k}(\tilde{\omega}_m) + F^{(b)}_{j,I_k}(\tilde{\omega}_{m'})\RP\kappa(\eta_k(m,m'))\psi_k(m')}{\sum_{m,m'} \kappa(\eta_k(m,m'))\psi_k(m')}\COMMA
\end{align}
where $\eta_{\kappa}(m,m')$ has been defined in \eqref{eq:eta} and $m,m'=1,2,\ldots,M_k$ and 
$\psi_k(m') := \exp\LP \int_{t^*_k}^{t_k} c_j(\tilde X^{(b)}(s,\tilde{\omega}_{m'}))ds \RP$ (according to Theorem \ref{thr:representation}).
%This averages will be used in \eqref{eq:EMiteration1}. 
Observe that we generate $M_k$ forward and reverse paths in the interval $I_k$ but we do not directly control  the number of exact or approximate SRN bridges that are formed. The number $M_k$ is chosen using a coefficient of variation criterion, as explained in Section \ref{sec:Mk}.
%such that either the number of SRN-bridges is of order $O(M_k)$ or we reach a computational budget $M_b$, which is $200$ in our numerical experiments. 
In Section \ref{sec:complexity}, we indicate an algorithm to reduce the computational complexity of computing those $\kappa$-weighted averages from $O(M_k^2)$ to $O(M_k \log (M_k) )$.

Finally, the Monte Carlo EM algorithm for this particular problem generates a {stochastic} sequence 
$\seqof{\hat{\theta}_{I\!I}^{(p)}}{p=1}{+\infty}$ staring from the initial guess
${\theta}^{(0)}_{I\!I}$ provided by phase I \eqref{eq:seedI} and evolving by 
\begin{equation}\label{eq:MCEMiteration}
 \hat{c}^{(p+1)}_{} =  \frac{\sum_{k=1}^K \avgsub{ R_{j,I_k}\SEP \mathcal{D}}{\kappa}{\hat{\theta}_{I\!I}^{(p)}}}{\sum_{k=1}^K \avgsub{ F_{j,I_k}\SEP \mathcal{D}}{\kappa}{\hat{\theta}_{I\!I}^{(p)}}}\COMMA
\end{equation}
where $\hat{\theta}_{I\!I}^{(p)} = \LP \hat{c}_{1}^{(p)},\ldots,\hat{c}_{J}^{(p)} \RP$. In Section \ref{sec:stopping}, a stopping criterion based on techniques widely used in Monte Carlo Markov chains is applied.

% flatex input end: [FREM.tex]

%\clearpage

%\clearpage
\section{Computational Details}
% flatex input: [compdetails.tex]
\label{sec:compdetails}
%\blue{
This section is intended to show computational details omitted in Section \ref{sec:FREM}. 
Here, we explain why and how we transform the clouds $\mathcal{X}_k$ consisting of endpoints of forward and reverse paths in the time interval $I_k$ at the time $t_k^*$, for $k{=}1,...,K$. 
Then, we explain how to chose the number of simulated forward and backward paths, $M_k$, in the time interval $I_k$ to obtain accurate estimates of the expected values of $R_{j,I_k}$ and $F_{j,I_k}$ for $j=1,2,\ldots,J$. 
Next, we show how to reduce the computational cost of computing approximate SRN bridges from $O(M_k^2)$ to $O(M_k \log (M_k))$ using a strategy introduced by Bayer and Schoenmakers \cite{BayerMC}. 
Finally, we indicate how to choose the initial seeds for  phase I and a stopping criteria for  phase II.

\subsection{On the Selection of the Number of Simulated Forward-Backward Paths}
\label{sec:Mk}
The selection strategy of the number of sampled forward-backward paths, $M_k$, for interval $I_k$, is determined by the following sampling scheme:
\begin{enumerate}
\item First sample $M$ forward-reverse paths (in the numerical examples we use $M{=}100$).
\item If the number of joined forward-reverse paths using a delta kernel is less than a certain threshold, $\gamma$, we transform the data as described in Section \ref{sec:trasf}. This data transformation allows us to use the Epanechnikov kernel \eqref{eq:epa}. In this way, we are likely to obtain a larger number of joined paths.
\item We then compute the coefficient of variation of the sample mean of the sum of the number of times that each reaction $j$ occurred in the interval $I_k$, $R^{(f)}_{j,I_k} {+} R^{(b)}_{j,I_k}$ and $F^{(f)}_{j,I_k} {+} F^{(b)}_{j,I_k}$, for $j{=}1,...,J$. Here $F^{(f)}_{j,I_k}=  \int_{I_k}
g_j(X^{(f)}(s))\,ds$ and the coefficient of variation of the sample mean of the sum $F^{(b)}_{j,I_k}=  \int_{I_k}
g_j(X^{(b)}(s))\,ds$. Further details can be found in Section \ref{sec:contobservedpaths}.
The coefficient of variation ($cv$) of a random variable is defined as the ratio of its standard deviation $\sigma$ over its mean $\mu$,
$cv := \frac{\sigma}{\abs{\mu}}$.
In this case, for the reaction channel $j$ in the interval $I_k$, we have:
$$cv_{\bar R}(I_k,j) = L_k^{-1/2}\, \frac{\sdev{R^{(f)}_{j,I_k}(\tilde{\omega}_m) {+} R^{(b)}_{j,I_k}(\tilde{\omega}_{m})}{L_k}}{{\avg{R^{(f)}_{j,I_k}(\tilde{\omega}_m) {+} R^{(b)}_{j,I_k}(\tilde{\omega}_{m})}{L_k}}}$$
and
$$cv_{\bar F}(I_k,j) = L_k^{-1/2}\, \frac{\sdev{F^{(f)}_{j,I_k}(\tilde{\omega}_m) {+} F^{(b)}_{j,I_k}(\tilde{\omega}_{m})}{L_k}}{{\avg{F^{(f)}_{j,I_k}(\tilde{\omega}_m) {+} F^{(b)}_{j,I_k}(\tilde{\omega}_{m})}{L_k}}} \COMMA$$

where $\sdev{Y}{L}{:=} \avg{Y^2}{L}-\avg{Y}{L}^2$ is the sample standard deviation of the random variable $Y$ over an ensemble of size $L$ and $\avg{Y}{L}{:=}\frac{1}{L}\sum_{m=1}^L Y(\omega_m)$ is its sample average. Here $L_k$ denotes the number of joined paths in the interval $k$, which is bounded by $M_k^2$. 
In the case that $L_k$ is small, we compute a bootstrapped coefficient of variation. 

The idea is that by controlling both coefficients of variation, we can control the variation of the $p$-th iteration estimation $\hat \theta_{II}^{(p)}$. Our numerical experiments confirm this fact. 
\item If each coefficient of variation is less than a certain threshold then the sampling for interval $I_k$ finishes, where $M_k$ is the total number of sampled paths, and  accepting the quantities in step 3. and also the quantities $\kappa(\eta_k(m,m'))\psi_k(m')$, $m,m'=1,...,L$ as defined in Section \ref{sec:contobservedpaths}. Otherwise, we sample additional forward-reverse paths (increasing the number of sampled paths at each iteration $M$) and go to step 2.
\end{enumerate}
This selection procedure is implemented in  Algorithm \ref{alg:fr_path}.
\subsection{On the Complexity of the Path Joining Algorithm}
\label{sec:complexity}
In this section, we describe the computational complexity of  Algorithm \ref{alg:fr_join}  for joining paths in  phase II, and show that this complexity is $O(M \log (M))$ on average.
%We have two cases: when $\kappa$ is  Epanechnikov's kernel, and when $\kappa$ is Kronecker's kernel.
%If the Epanechnikov kernel is used, since it has compact support,
%is supported in a ball of radius $R > 0$, and we use a bandwidth proportional to $M^{1/d}$, we have that  Algorithm \ref{alg:fr_join}
%the join algorithm can be implemented with a computational cost of order $O(M)$.
%In case ii), the complexity of joining $M$ paths is $O(M^2)$, because the assumption \textbf{A)} (described below) cannot be fulfilled in practice.

Let us describe the idea.
First, fix  a time interval $I_k$ and a  reaction channel $j$. We use
 the following double sum as an example,%$$\sum_{m=1}^{M} \sum_{m'=1}^{M} \LP R^{(f)}_{j,I_k}(\tilde{\omega}_m) + R^{(b)}_{j,I_k}(\tilde{\omega}_{m'})\RP\kappa_{m,m'} \COMMA$$ 
\begin{equation*}
\sum_{m=1}^{M} \sum_{m'=1}^{M} \LP R^{(f)}_{j,I_k}(\tilde{\omega}_m) + R^{(b)}_{j,I_k}(\tilde{\omega}_{m'})\RP\kappa_{m,m'}\PERIOD
\end{equation*}
A double sum like this one appears in the numerator of \eqref{eq:averagesIk}.  
%This double sum is one of the two components of the MCEM estimator \eqref{eq:EMiteration1}, the other component is an analogous double sum. 
%Here $\kappa_{m,m'}:=\kappa_{m,m'}(X_m^{(f)}(t^*),X_{m'}^{(b)}(t^*))$ represents the kernel weight factor, and depends on the end points of the $M$ samples. 
Instead of computing a double loop which always takes $O(M^2)$ steps (and many of those steps contribute 0 to the sum), we take the following alternative approach: 
let 
$\times_{i=1}^d [A_i,B_i]$ be the smallest hyperrectangle  of sides $[A_i,B_i]$, $i=1,...,d$ that contains the cloud 
$\mathcal{Y}$, defined in \eqref{eq:cloudY}. Let us also assume that 
$A_i,B_i$, $i=1,...,d$ are integers.
%(or cloud) of end points or its linear \blue{P:(homeomorphic, canonical?)} transformation  (as explained in \ref{transformation}) 
%can be fit on the set $\times_{i=1}^d [A_i,B_i]$, that is, a hyperrectangle of sides $[A_i,B_i]$, $i=1,...,d$. 
The length $B_i-A_i$ depends on how sparse the cloud is in its $i$-th dimension. Given the cloud, it is easy to check that the values $A_i,B_i$, $i=1,...,d$ can be computed in $O(M)$ operations. Now, we subdivide the hyperrectangle into sub-boxes of size-length  1, with sides parallel to the coordinate axis. %Let $C$ be the total number of such sub-boxes. 

Since we have a finite number of those sub-boxes, we can associate an index for each one in such a way that it is possible to directly retrieve each one using a suitable data structure (for example  an efficient sparse matrix or a hash table). The average access cost of such structure is constant with respect of $M$. For each sub-box, we associate a list of \emph{forward} points that ended up in that sub-box. It is also direct to see that the construction of such a structure takes a computational cost of $M$ steps on average.
%If we assume that ``\textbf{A)} the average list size of each sub-box is $O(1)$'',{(which can be proved under certain conditions),} 
Then, instead of evaluating the double sum which has $O(M^2)$ steps, we evaluate only the non zero terms. This is because, when a kernel $\kappa$ is used, $\kappa(x,y) \neq 0$ if and only if $x$ and $y$ are situated in neighboring sub-boxes. 
That is,
\begin{align*}
\sum_{m=1}^{M} \sum_{m'=1}^{M}& \LP R^{(f)}_{j,I_k}(\tilde{\omega}_m) +  R^{(b)}_{j,I_k}(\tilde{\omega}_{m'})\RP \kappa_{m,m'}  \\
 &=\sum_{m'=1}^{M} \sum_{i=1}^{3^d} \sum_{l=1}^{n(b_i)} \LP R^{(f)}_{j,I_k}(\tilde{\omega}_{\ell(l)}) + R^{(b)}_{j,I_k}(\tilde{\omega}_{m'})\RP\kappa_{\ell(l),m'} \COMMA 
\end{align*}
where $n(b_i)$ is the total quantity of \emph{reverse} end points associated with the $i$-th neighbor of the sub-box to which the \emph{forward} end-point, $\tilde Y^{(f)}(t_k^*,\tilde{\omega}_{m})$, belongs, whereas $\ell(l)$ indexes one of those reverse end points. %Here $I_k$ is a fixed interval and $j$ is a fixed reaction channel.
Note that the constant of this complexity depends exponentially on the dimension ($3^d$). %This reasoning is based on \cite{bayer}.

The cost that dominates the triple sum on the right-hand side is the expected maximum number of reverse points that can be found in a sub-box. This size can be proved to be $O(\log (M))$, which makes the whole joining algorithm of order $O(M \log (M))$. For additional details we refer to \cite{Bayer}.

%When Epanechnikov's kernel is used with a kernel-width of size 1, $\kappa_e(x,y) \neq 0$ if and only if $x$ and $y$ are situated in neighboring sub-boxes. In this case, the assumption \textbf{A)} is realistic, and, on average, the number of paths that ended in a given neighborhood is of order $O(1)$. 

%If Kronecker's  kernel is used, then the range of the end points belongs to $\zset$ and then assumption \textbf{A)} may not be realistic. In this case, the average order tends to be larger than $M$. This is because we only use Kronecker's  kernel in a regime in which the values are clustered around a few neighbors.
\subsection{A Linear Transformation for the Epanechnikov Kernel}
\label{sec:trasf}
\begin{figure}[h!]
\centering
\begin{minipage}{0.49\textwidth}
\includegraphics[width=\textwidth]{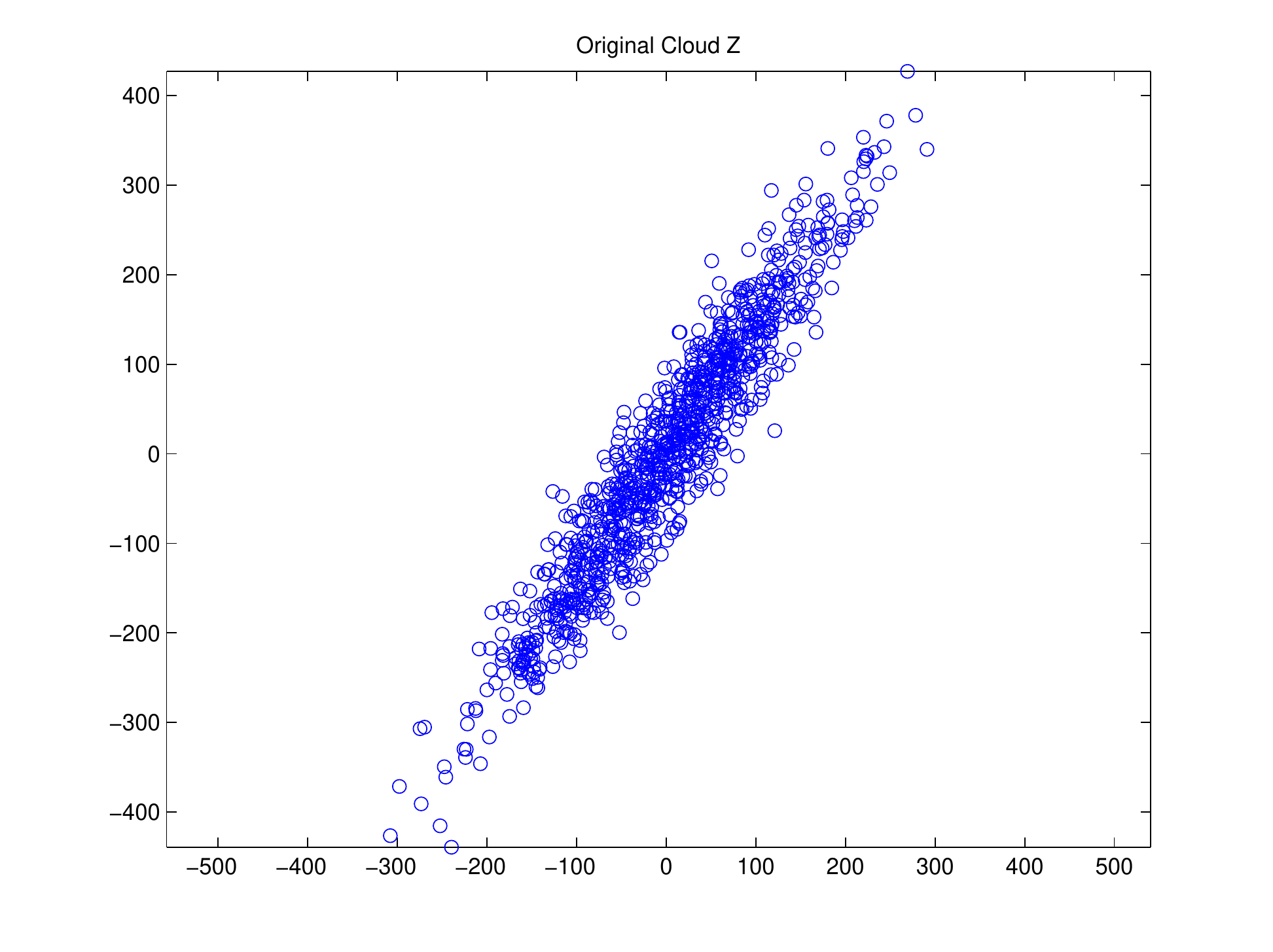}
\end{minipage}
\hfill
\begin{minipage}{0.49\textwidth}
\includegraphics[width=\textwidth]{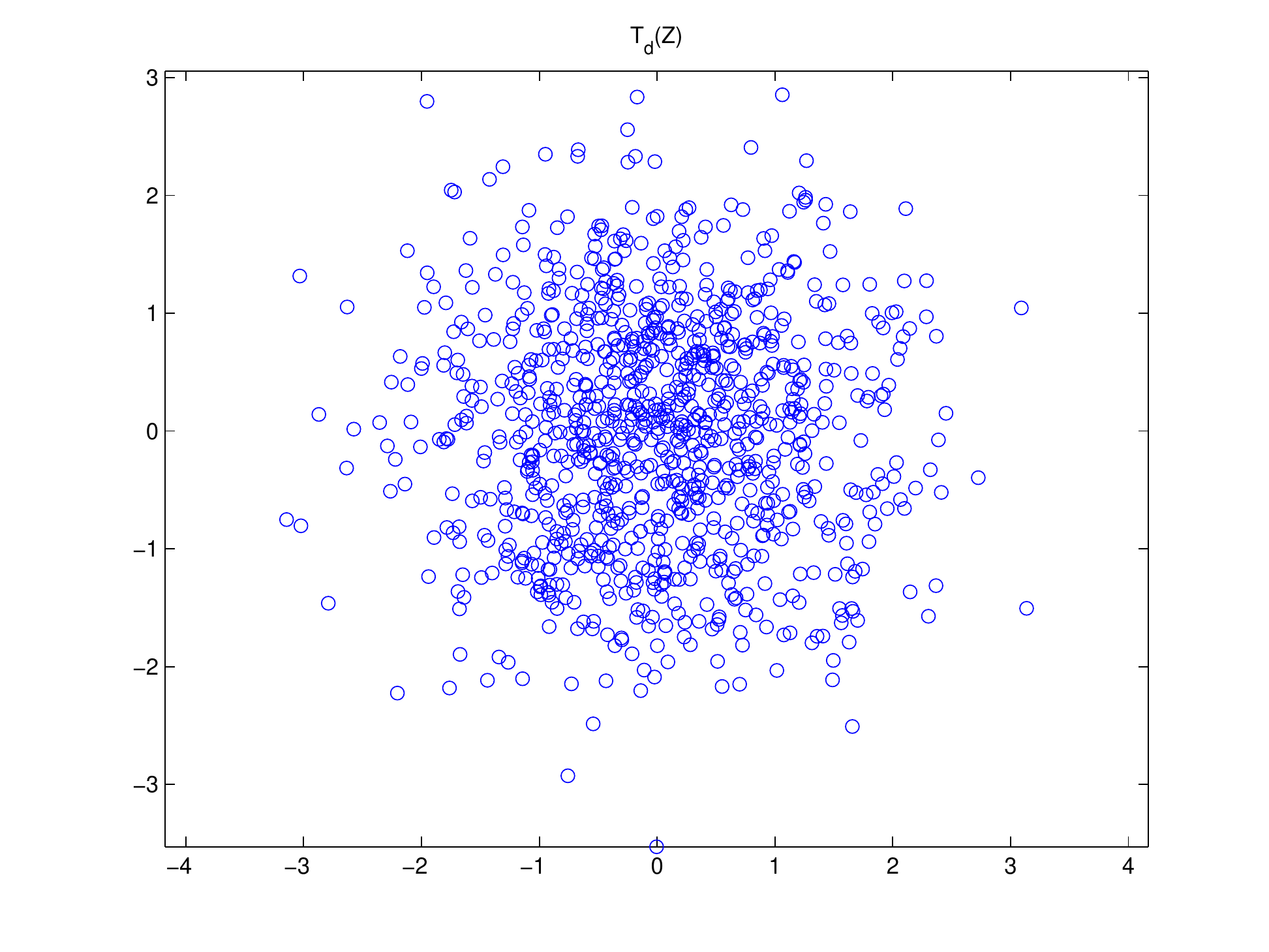}
\end{minipage}
\caption{Left: A bivariate Gaussian cloud, $\mathcal{Z}$. Right: Its corresponding decorrelated and scaled version $T(\mathcal{Z})$.}
\label{fig:CloudZ}
\end{figure}

%\red{In this section we present a data transformation related to the window size of the Epanechnikov kernel.}
Our numerical experiments show that clouds formed by the endpoints of simulated paths,  $\mathcal{X}$, usually have a shape similar to the \emph{cloud $\mathcal{Z}$} shown in the left panel of Figure \ref{fig:CloudZ}.

It turns out that partitioning the space into $d$-dimensional cubes with sides parallel to the coordinate axis is not idle for selecting  kernel domains and consequently for finding SRN bridges. 
It is more natural way to divide the space into a system of parallelepipeds with sides parallel to the principal directions of   cloud $\mathcal{Z}$ having sides proportional to the lengths of its corresponding semi-axes to use  as  supports for our kernels.

Another way of proceeding (somehow related but not totally equivalent)  is to transform the original cloud $\mathcal{Z}$ to obtaining another cloud  $T(\mathcal{Z})$ with a near-spherical shape. Then, scale it to have on average one point of the cloud in each $d$-dimensional cube (with sides parallel to the coordinate axis). In this new cloud, $H(\mathcal{Z})$, we can naturally find neighbors using the algorithm described in Section \ref{sec:complexity} below and we have the Epanechnikov kernel to assign weights. 
This is why in Section \ref{sec:FREM}  we wanted to transform the data $\mathcal{X}_k$ into an  isotropic cloud, such that, every unitary cube centered in 
$\tilde Y^{(f)}(t_k^*,\tilde{\omega}_m')$ 
contains, on average, one point of the cloud 
$\cup_{m'} \tilde Y^{(b)}(t_k^*,\tilde{\omega}_m')$.

We will now describe the details of the aforementioned transformations.

First, we show a customary procedure in statistics to motivate the  transformation.
Let $\Sigma := \text{cov}(\mathcal{Z})$ be the sample covariance matrix computed from a cloud of points $\mathcal{Z}$. %$\mathcal{X}_k$ defined in \eqref{eq:cloudX}. 
To obtain a de correlated version of $\mathcal{Z}$, the  linear transformation $T(z) = \Sigma ^{-1/2} \, z$ is widely used in statistics.  
For example, consider a cloud $\mathcal{Z}$ of points obtained by sampling $10^3$ independent highly correlated bi-variate Gaussian random variables. The corresponding cloud $T(\mathcal{Z})$, depicted in the right panel of Figure \ref{fig:CloudZ}, shows the aspect of a sphere with a radius $3$ units.
\begin{figure}[h!]
\centering
%\begin{minipage}{0.49\textwidth}
%\includegraphics[width=\textwidth]{./Figures/CloudX.pdf}
%\end{minipage}
%\hfill
%\begin{minipage}{0.49\textwidth}
\includegraphics[width=0.49\textwidth]{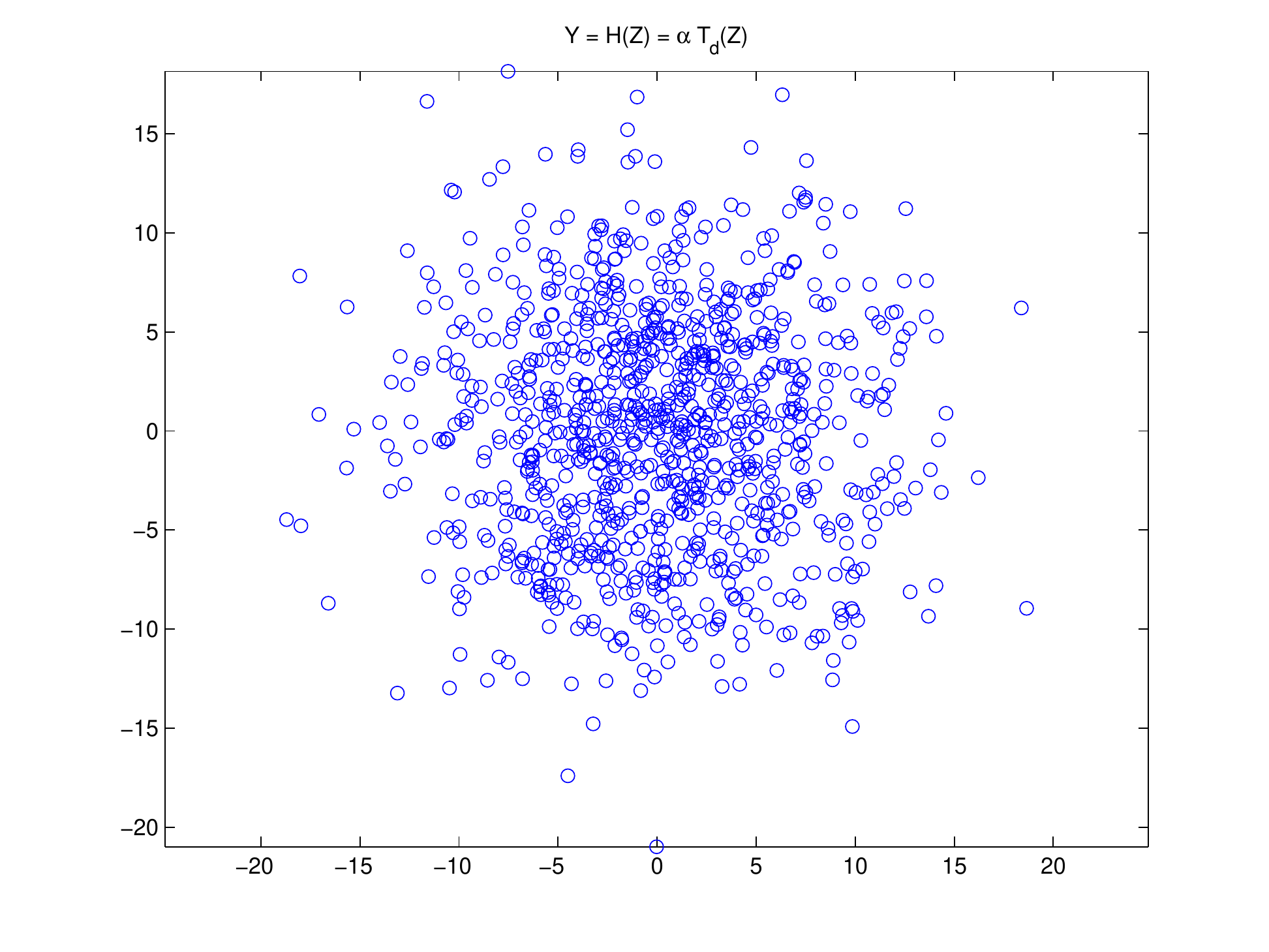}
%\end{minipage}
\caption{Cloud $H(\mathcal{Z})$.}
\label{fig:CloudZ3}
\end{figure}
The next step is to obtain a radius $\alpha$ such that the volume of a $d$-dimensional sphere of radius $3\alpha$ equals to the volume of $M$ unitary $d$-dimensional cubes.
From the equation $M = (3\alpha)^d\, V_d$,
%\begin{align}
%M &= (3\alpha)^d\, V_d\\ \nonumber
%  &= (3\alpha)^d\, \frac{\pi^{d/2}}{\Gamma(d/2+1)}
%\end{align} 
we obtain $\alpha = \frac{1}{3} \,(\frac{M}{V_d})^{1/d}$, where $V_d = \frac{\pi^{d/2}}{\Gamma(d/2+1)}$ is the volume of the unitary sphere in $\rset^d$. Therefore, the linear transformation $H$ is defined by $H(x):= \alpha T(x)$.  
The result of this transformation is depicted in Figure \ref{fig:CloudZ3} in our Gaussian example.

In  general, we do not expect to have a Gaussian-like distribution for $\mathcal{X}_k$, however it seems to be a good approximation in our numerical examples. At this point, it is worth mentioning that in examples with several dimensions (species), the number of approximate SRN bridges we get by using the transformation may be of the order of  $M^2$. This indicates that the bandwidth is too large, and consequently the bias introduced in the estimation may be large. In these cases, we expand $\alpha$ by a factor of $1.5$, for example, until $O(M)$ approximate bridges are formed. Generally, one or two expansions are enough. 

A motivation for the Gaussian approximation is that, for short time intervals and in certain regimes of activity of the system, specially where the total propensity, $a_0$, is high enough, a Langevin approximation of our SRN provides an Ornstein-Uhlenbeck process, which can potentially be close in distribution to our SRN (see \cite{ourInf}). 

\begin{rem}
According to the transformation $H$, the kernel used in our case is approximately equal to
\begin{equation*}
\kappa_H(z) := \frac{1}{\det(H)}\,\kappa\LP H^{-1}(z)\RP\COMMA
\end{equation*}
where $\kappa$ is the Epanechnikov kernel defined in \eqref{eq:epa}, since it corresponds with the continuous case and not with the lattice case.
\end{rem}

\begin{rem}
We can even consider a perturbated version of $T$, say $T_c$, by adding a multiple of the diagonal matrix formed by the diagonal elements of $\Sigma$, \ie, $T_c = (\Sigma + c\,\text{diag}(\Sigma))^{-1/2}$, where $c$ is a positive constant of order $O(1)$. The linear transformation $T_c$ can be considered as a regularization of $T$ that does not change the scale of the transformation $T$.
\end{rem}

\subsection{On the Stopping Criterion}
\label{sec:stopping}

A well-known fact about the EM algorithm is that, given a starting point, it converges to a saddle point or a local maximum of the likelihood function. Unless we know beforehand that the likelihood function has a unique global maximum, we cannot be sure that the output of the EM Algorithm is the MLE we are looking for.
The same phenomenon occurs in the case of the Monte Carlo EM algorithm, and for that reason Casella and Robert in \cite{Casella} recommend  generating a set of $N$ (usually $N$ around five) parallel-independent Monte Carlo EM sequences starting from a set of over dispersed initial guesses. Usually, we do not know even the scale of the coordinates of our unknown parameter $\theta = (c_1,c_2,\ldots,c_d)$. For that reason, we recommend running only  phase I of our algorithm over a set of uniformly distributed random samples drawn from a $d$-dimensional hyperrectangle $\prod_{i=1}^d (0,C_i]$, where $C_i$ is a reasonable, case dependent, upper bound for each reaction rate parameter $c_i$.
We observed in our numerical experiments that the result of this procedure is a number of points laying on a low dimensional manifold. 
Once this manifold is identified, $N$ different initial guesses are taken as over dispersed seeds for phase II.  

Note that the stochastic iterative scheme given by formula \eqref{eq:MCEMiteration} may be easily adapted to produce $N$ parallel stochastic sequences, $\seqof{\hat{\theta}_{I\!I,i}^{(p)}}{p=1}{+\infty}$, where, for each $i=1,2,\ldots,N$, the distribution of the random variable $\hat{\theta}_{I\!I,i}^{(p+1)}$ depends on its history of realizations, $\seqof{\hat{\theta}_{I\!I,i}^{(k)}}{k=1}{p}$, only through its previous value, $\hat{\theta}_{I\!I,i}^{(p)}$. In this sense, the $N$ sequences, $\seqof{\hat{\theta}_{I\!I,i}^{(p)}}{p=1}{+\infty}$, are  MCMC sequences \cite{Norris, Casella}. 

There is a number of convergence assessment techniques or convergence diagnostic tools in the MCMC literature; in this article, we adopt the $\hat R$ criterion by Gelman and Rubin \cite{Rhat, BayesianData}, which monitors the convergence of $N$ parallel random sequences $\seqof{\psi_i^{(p)}}{p=1}{+\infty}$, where $i=1,2,\ldots,N$. 

Compute:
\begin{align*}
B_p &:= \frac{1}{N-1} \sum_{i=1}^N \LP \bar{\psi}_{p, i} - \dbar{\psi}_p\RP^2, \text{ where } \bar{\psi}_{p, i} := \frac 1 p \sum_{k=1}^p \psi_i^{(k)} \text{ and } \dbar{\psi}_p := \frac 1 N \sum_{i=1}^N \bar{\psi}_{p, i},  \text{ and}\\
W_p &:= \frac{1}{N} \sum_{i=1}^N s^2_{p,i}, \text{ where }  s^2_{p,i} := \frac{1}{p-1} \sum_{k=1}^p \LP \psi_i^{(k)} - \bar{\psi}_{p, i}\RP^2\PERIOD
\end{align*}
Then define
\begin{align}
V_p &:= \frac{p-1}{p} W_p + B_p \text{ and } \hat{R}_p := \sqrt{\frac{V_p}{W_p}} \PERIOD
\end{align}
%\begin{align}
%\nonumber \bar{\psi}_{\cdot, i} &:= \frac 1 p \sum_{k=1}^p \psi_i^{(k)}\\
%\nonumber \dbar{\psi} &:= \frac 1 N \sum_{i=1}^N \bar{\psi}_{\cdot, i}\\
%B &:= \frac{p}{N-1} \sum_{i=1}^N \LP \bar{\psi}_{\cdot, i} - \dbar{\psi}\RP^2\\
%\nonumber s^2_i &:= \frac{1}{p-1} \sum_{k=1}^p \LP \psi_i^{(k)} - \bar{\psi}_{\cdot, i}\RP^2\\
%W &:= \frac{1}{N} \sum_{i=1}^N s^2_i\\
%\nonumber V &:= \frac{p-1}{p} W + \frac{1}{p} B\\
%\hat R &:= \sqrt{\frac{V}{W}}
%\end{align}
$B$ and $W$ are known as between and within variances, respectively. 
It is expected that $\hat R$ (potential scale reduction) declines to $1$ as $p\to+\infty$. In our numerical experiments we use $1.4$ as a threshold. %Quoting Gelman and Shirley in  Chapter 6 of \cite{BrooksGelmanJonesMeng201105} ``At convergence, the chains will have mixed, so that the distribution of the simulations between and within chains will be identical, and the ratio $\hat R$ should equal 1. If $\hat R$ is greater than 1, this implies that the chains have not fully mixed and that further simulation might increase the precision of inferences''.

Observe that if for all $p$, the values $\bar{\psi}_{p, i}$ are grouped in a very small cluster, \ie, $\bar{\psi}_{p, i} \approx \dbar{\psi}_p$ and therefore we have essentially only one Markov chain, then $B_p$ is close to zero and $\hat{R}_p \approx \sqrt{\frac{p-1}{p}} \to 1$ as $p \to +\infty$ independently of the behavior of the chain. To avoid this undesirable situation, we propose to observe also the behavior of the moving averages of order $L$, that is, 
\begin{align}\label{eq:philags}
\tilde{\psi}_{p}:= \frac{1}{N} \sum_{i=1}^N \LP \tilde{\psi}_{p,i}-\tilde{\psi}_{p-1,i}\RP^2 \text{ where } \tilde{\psi}_{p,i} := \frac{1}{L}\sum_{\ell=0}^{L-1} {\psi}^{(p-\ell)}_{i} \PERIOD
\end{align}
We stop when $\tilde{\psi}_{p}$ is sufficiently small. %In our numerical experiments we use $L=3$.

Once we stop to iterate after $p^*$ iterations, the individual outputs \begin{equation*}
\hat{\theta}_{I\!I,1}^{(p^*)}, \hat{\theta}_{I\!I,2}^{(p^*)},\ldots, \hat{\theta}_{I\!I,N}^{(p^*)}
\end{equation*} 
form a small cluster. 
Although we cannot be certain that this cluster is near the MLE, we do have at least some confidence. Therefore, we can use the mean of this small cluster as a MLE estimation of our unknown parameter, $\theta$. 
Otherwise, if we have two or more clusters or over dispersed results, we should perform a more careful analysis.
\begin{rem}
The $\hat R$ stopping criterion only works if the over dispersed seeds obtained in phase I lie in the basin of attraction of one local maximum of the likelihood function. Otherwise $\hat R$ may not decrease to 1, even worse, it may go to $+\infty$. For that reason, it is recommendable to monitor the evolution of $\hat R$. In our numerical examples we have that $\hat R$ is decreasing and we stop the algorithm using $\hat R_0=1.4$ as a threshold.
\end{rem}

%}

% flatex input end: [compdetails.tex]

%\clearpage

%\clearpage
\section{Numerical Examples}
% flatex input: [num_examples]
\label{sec:numex}

In this section, we present numerical results that show the performance of our FREM algorithm.
In phase I, we use the alternative definition of  ${\theta}_{I\!I,i}^{(0)}$ described in Remark \ref{rem:alternative}. For  phase II, we run $N=4$ parallel sequences using $1.4$ as a threshold for $\hat R$ (described in Section \ref{sec:stopping}). 
The moving average order used in all numerical examples is $L=3$ (see formula \ref{eq:philags}), and the associated tolerance is $0.05$.
As a point estimator of $\theta$, we provide the cluster average of the sequence 
$\hat{\theta}_{I\!I,1}^{(p^*)}, \hat{\theta}_{I\!I,2}^{(p^*)},\ldots, \hat{\theta}_{I\!I,N}^{(p^*)}$. 

For each example, we report i) the number of iterations of  phase II, $p^*$; ii) a table containing a) the initial points, ${\theta}_{I,i}^{(0)}$, b) the outputs of the phase I, ${\theta}_{I\!I,i}^{(0)}$, and c)  the outputs of  phase II, $\hat{\theta}_{I\!I,i}^{(p^*)}$; and iii) a Figure with all those values.

For the examples wehre we generate synthetic data, we provide the seed parameter $\theta_G$ we used to generate the observations. It is important to stress that the distance from our point estimator to $\theta_G$ depends of the number of generated observations.

\newcommand{\BD}{num_examples/main_new_generic_2014_10_12_7_53_46_BD/runs/Birth-death_2014_10_12_8_5_18}
\newcommand{\DEC}{num_examples/main_new_generic_2014_10_12_20_6_29_DEC2/runs/Decay_two_reactions_2014_10_12_20_6_29}
\newcommand{\SIR}
{num_examples/main_new_generic_2015_1_6_11_51_43_SIR2D_new/runs/SIR2D_2015_1_6_11_51_43}
\newcommand{\WEAR}
{num_examples/main_new_generic_2014_11_7_22_16_54_WEAR/runs/Wear_Cilindri_2014_11_7_22_16_55/}

\newcommand{\DECENS}{num_examples/ensemble/Decay_two_reactions_2015_1_19_1_13_11}
\newcommand{\WEARENS}{num_examples/ensemble/Wear_Cilindri_2015_1_20_12_7_22}
\newcommand{\BDENS}{num_examples/ensemble/Birth-death_2015_1_19_18_6_55}
\newcommand{\SIRENS}{num_examples/ensemble/SIR2D_2015_1_18_17_54_56}

\newcommand\rsp{\rule[10pt]{0pt}{0pt}}

\subsection{The Decay Process}
We start with a simple decay model with only one species and two reaction channels.
%\begin{align*}
%X \xrightarrow{c_1} \emptyset,& \ \ X \xrightarrow{c_2} \emptyset
%\end{align*}
Its stoichiometric matrix and propensity function are:
\begin{align*}
\nu^T = \left( 
 \begin{array}{cccc}     -1   \\
     -4  
 \end{array} 
 \right) \mbox{   and   }  a(X) = \left( \begin{array}{l}  c_1 X \\ c_2 X \cdot \indicator{X \geq 4} \end{array} \right), \mbox{   respectively} \PERIOD
\end{align*}

%We set $X_0{=}100$, $T{=}1$ and we use a uniform grid of measurements of size $\Delta t {=}\frac{1}{16}$.  We generate a single data trajectory, using the parameters $c_1{=}3.78$ and $c_2{=}7.2$.
%
%%We start from the points: 
%The FREM estimation is $(3.681488;7.495539)$, and the MLE is $(3.708467;7.555186)$.
%   
%The FREM algorithm took 3 iterations to converge (minimum imposed), for a R threshold of $1.4$.

\begin{figure}[h!]
\centering
\includegraphics[scale=0.4]{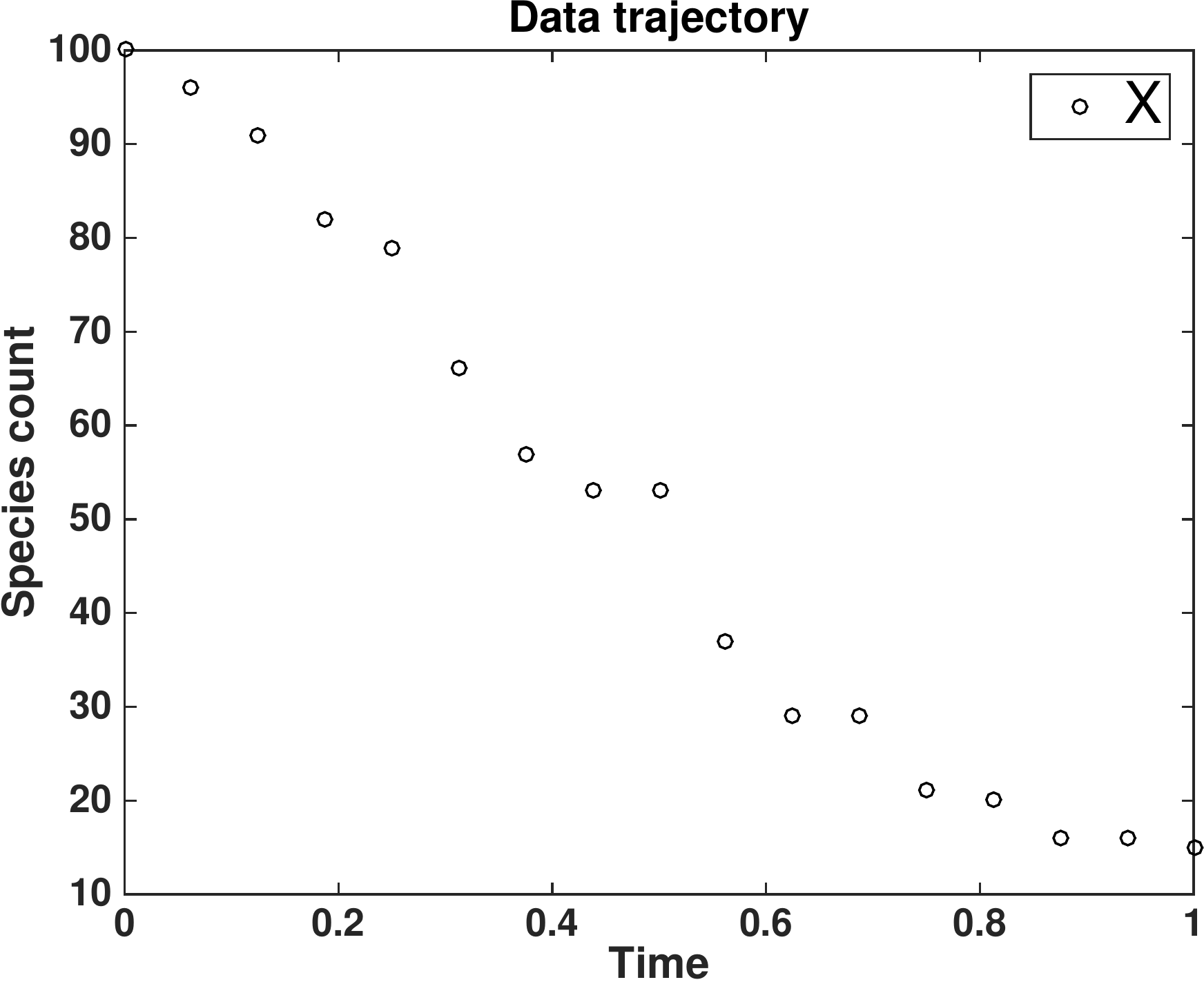}
\caption{Data trajectory for the decay example. This is obtained by observing the values of an SSA path at uniform time intervals of size
$\Delta t {=}1/16$.}
\label{fig:dataDec}
\end{figure}

We set $X_0{=}100$, $T{=}1$ and consider synthetic data observed in uniform time intervals of size $\Delta t {=}\frac{1}{16}$. This determines a set of $17$ observations generated from a single path using the parameter $\theta_G {=}(3.78,7.20)$. The data trajectory is shown in Figure \ref{fig:dataDec}.

For this example, we use $N{=}4$ FREM sequences starting at $\theta_{I,1}^{(0)}{=}(1, 5)$, $\theta_{I,2}^{(0)}{=}(6, 5)$, $\theta_{I,3}^{(0)}{=}(1, 9)$, and $\theta_{I,4}^{(0)}{=}(6, 9)$. In this  and  the following examples, for each interval we run a minimum of $M=100$ forward-reverse sample paths and we set a coefficient of variation threshold of $0.1$ (see Section \ref{sec:Mk}).

\begin{figure}[h!]
\centering
\begin{minipage}{0.48\textwidth}
\includegraphics[width=\textwidth]{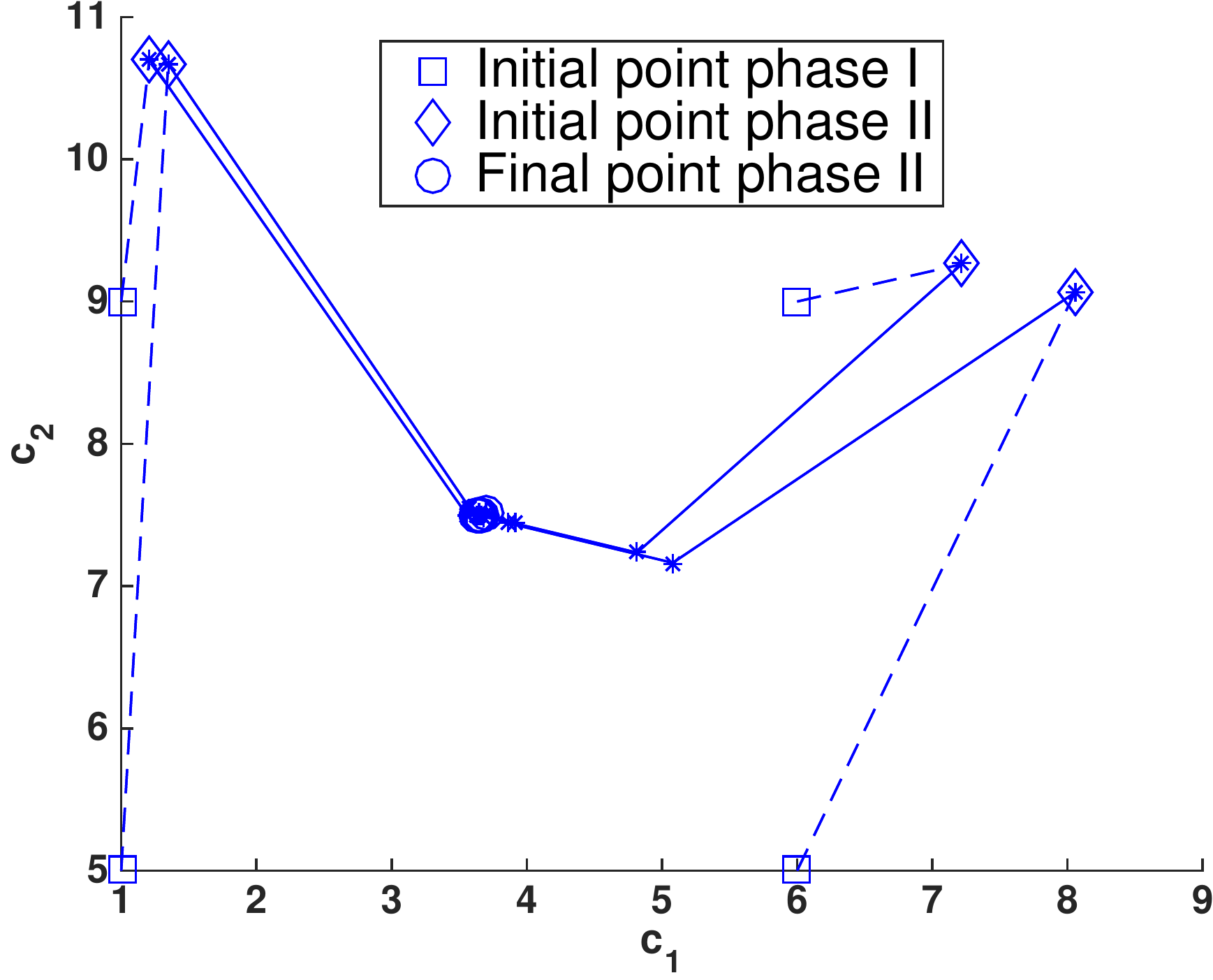}
\end{minipage}
%\hfill
\begin{minipage}{0.51\textwidth}
\includegraphics[width=\textwidth]{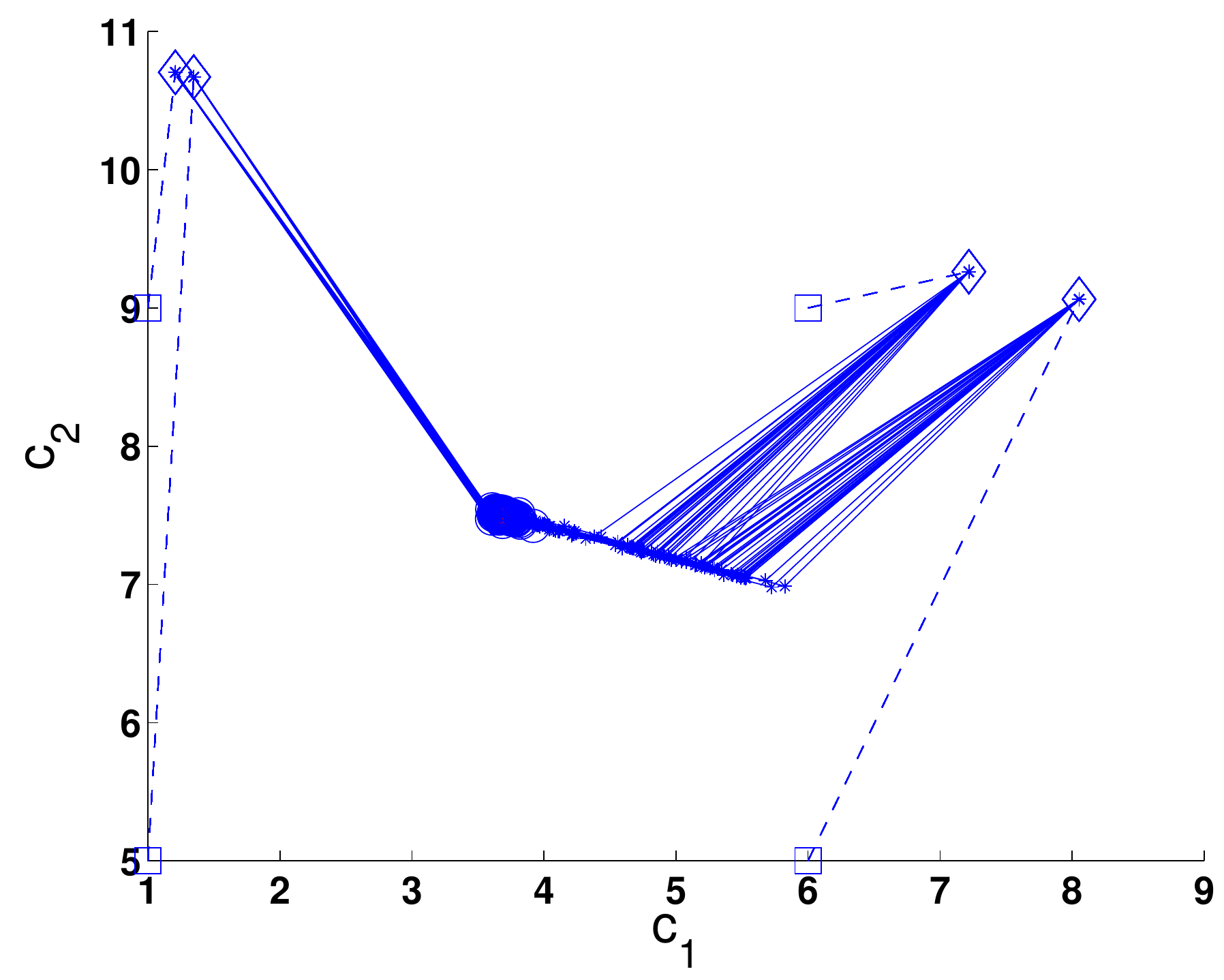}
\end{minipage}
\caption{Left: One FREM estimation (phase I and phase II) for the decay example. The $N$ final values of this particular run are shown as circles. Right: We show 30 independent runs of the FREM algorithm.}
\label{fig:dec2}
\end{figure}

We illustrate one run of the FREM algorithm in the left panel of Figure \ref{fig:dec2} and in Table \ref{tab:dec}. For that run, the cluster average is 
$\hat{\theta} {=} (3.68, 7.50)$, and it took $p^*{=}3$ iterations to converge for a $\hat{R}$ threshold equal to 1.4.
We take $\hat{\theta}$ as a MLE point estimation of the unknown parameters. 
%Details can be found in Table \ref{tab:dec} and Figure \ref{fig:dec2}.

\begin{table}[h!]
\centering
\begin{tabular}{cccc}
$i$ & $\square {=} {\theta}_{I,i}^{(0)}$ & $\Diamond {=} {\theta}_{I\!I,i}^{(0)}$ & $\bigcirc {=} \hat{\theta}_{I\!I,i}^{(p^*)}$ \\ 
\hline  \rsp
1 &(1, 5) & (1.35, 10.67) & (3.65, 7.52) \\ 
%\hline 
2 &(6, 5) & (7.85, 9.11) & (3.80, 7.46) \\ 
%\hline 
3 & (1, 9) & (1.20, 10.71)  & (3.63, 7.50) \\ 
%\hline 
4 &(6, 9) & (7.06, 9.30) & (3.65, 7.50) \\ 
%\hline 
\end{tabular} 
\caption{Values computed by one run of the FREM Algorithm for the decay example, corresponding to the left panel of Figure \ref{fig:dec2}.}
\label{tab:dec}
\end{table}

We computed an ensemble of 30 independent runs (and obtained 30 cluster averages). The result is shown in the right panel of Figure \ref{fig:dec2}. We observe that the variability of the cluster average is indeed very small, indicating the robustness of the method and that 1.4 is a reasonable choice as a threshold for $\hat R$. Details are shown in Table \ref{tab:dec_ens}.
\begin{table}[h!]
\centering
\begin{tabular}{c|cccc}
%$n$ & \multicolumn{4}{c}{$\hat{c}_n$} \\
%\cline{2-5} \rsp
& Average & Average CI at $95\%$ & Min Value & Max Value  \\
\hline  \rsp
$\hat{c}_1$ &  3.69   &(3.681, 3.699) & 3.66   & 3.77    \\ 
$\hat{c}_2$ & 7.50  &(7.495, 7.505) & 7.48 & 7.51
\end{tabular} 
%\caption{Values computed for an ensemble of 30 independent runs of the FREM Algorithm for the decay example. In each run we obtain a cluster average. For each parameter, we show its ensemble average, a $95\%$ confidence interval for the ensemble average, and finally the minimum and maximum values recorded in the ensemble.}
\caption{Values computed for an ensemble of 30 independent runs of the FREM algorithm for the decay example. 
In each run, we obtain a cluster average, $\hat{\theta}^{(i)}$, as an MLE point estimate. Define $\mathcal{C} {:=} \seqof{\hat{\theta}^{(i)}}{i=1}{30}$.
For each unknown coefficient $c_j$ in $\theta$, we show i) the average of    $\mathcal{C}$, ii) a $95\%$ confidence interval for the mean of $\mathcal{C}$, and iii) the minimum and maximum values of $\mathcal{C}$. }
\label{tab:dec_ens}
\end{table}

\begin{rem}
Recall that the distance between the value $\theta_G$ used to generate synthetic data and the estimation $\hat{\theta}$ is meaningless for small data sets. The relevant distance in this estimation problem is the one we obtain from our FREM algorithm $\hat{\theta}$ and the 
$\hat{\theta}_{\text{MLE}}$ based on maximizing the true likelihood function however, the later is not available in most cases.
\end{rem}

\subsection{Wear in Cylinder Liners}
\label{ex:wear}
\begin{figure}[htp!]
\centering
\includegraphics[scale=0.5]{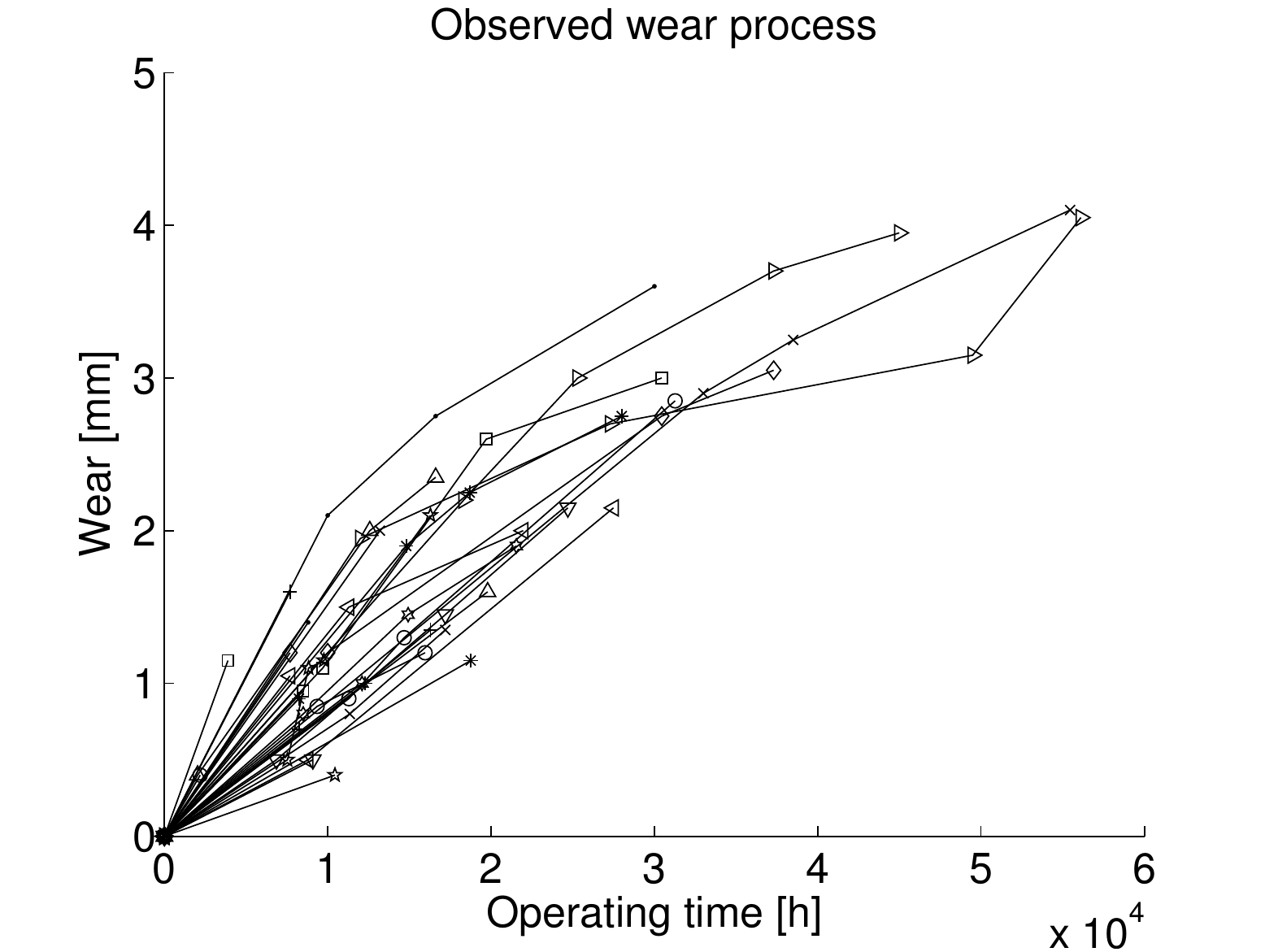}
\caption{Data set from {\rm \cite{gio2011}}. Data refer to cylinder liners used in ships of the Grimaldi Group. 
%\red{plot confidence band}
}
\label{fig:data}
\end{figure}

We now test our FREM algorithm by using real data. 
%This example is similar to the previous one from the modeling point of view, but in this case, we have real data.
%This example, which 
The data set $\mathbf{w} = \{w_i\}_{i=1}^n$, taken from \cite{gio2011}, consists of wear levels observed on $n= 32$ cylinder liners of eight-cylinder SULZER engines as measured by a caliper with a precision of 
$\Delta = 0.05$ mm. Data are presented in Figure \ref{fig:data}.

% Figure 1

The finite resolution of the caliper allows us to represent the set of possible measurements using a finite lattice.
Let $X(t)$ be the \textit{thickness process} derived from the wear of the cylinder liners up to time $t$, i.e., $X(t) = X_0 - W(t)$, where $W$ is the wear process and $X_0$ is the initial thickness. The final time of some observations is close to $T{=}60,000$ hours.

We model $X(t)$ as a decay processes with two reaction channels and  $\Delta = 0.05$, since a simple decay process is not enough to explain the data. The two considered intensity-jump pairs are $(a_1(x),\nu_1) =  (c_1x, -\Delta)$ and $(a_2(x),\nu_2) = (c_2x, -4\Delta)$. 
%where $k$ is a positive integer to be determined, and 
Here, $c_1$ and $c_2$ are coefficients with dimension $(\text{mm}\cdot \text{hour})^{-1}$. 
%Therefore, the probability of observing a thickness decrement in a small time interval $(t,t+dt)$ is
%\begin{eqnarray}\label{mod:tworeactions}
%{P(X(t+dt) = X(t){-}\Delta \bigm| X(t)=x) = \indicator{x \geq \Delta} c_1x \,dt ,}\\
%{P(X(t+dt) = X(t){-}k\Delta \bigm| X(t)=x) = \indicator{x \geq k\Delta} c_2x \,dt \PERIOD} \nonumber 
%\end{eqnarray}

The linear propensity functions, the value $X_0{=}5$ mm and the initial values for phase I:  $\theta_{I,1}^{(0)}{=}(1,1)$, $\theta_{I,2}^{(0)}{=}(10,1)$, $\theta_{I,3}^{(0)}{=}(1,10)$ and $\theta_{I,4}^{(0)}{=}(10,10)$, are motivated by previous studies of the same data set  (see \cite{ourInf} for details).
   
%The FREM estimation is $(8.9126;5.7400)$, and the MLE is $(10.1078;5.5277)$.
%      
%The FREM algorithm took 93 iterations to converge, for a R threshold of $1.4$.
%For this example we ran $N{=}4$ FREM sequences.

In our computations, we re scaled the original problem by setting 
$\Delta{=}1$ and $T{=}1$.

\begin{figure}[h!]
\centering
\begin{minipage}{0.48\textwidth}
\includegraphics[width=\textwidth]{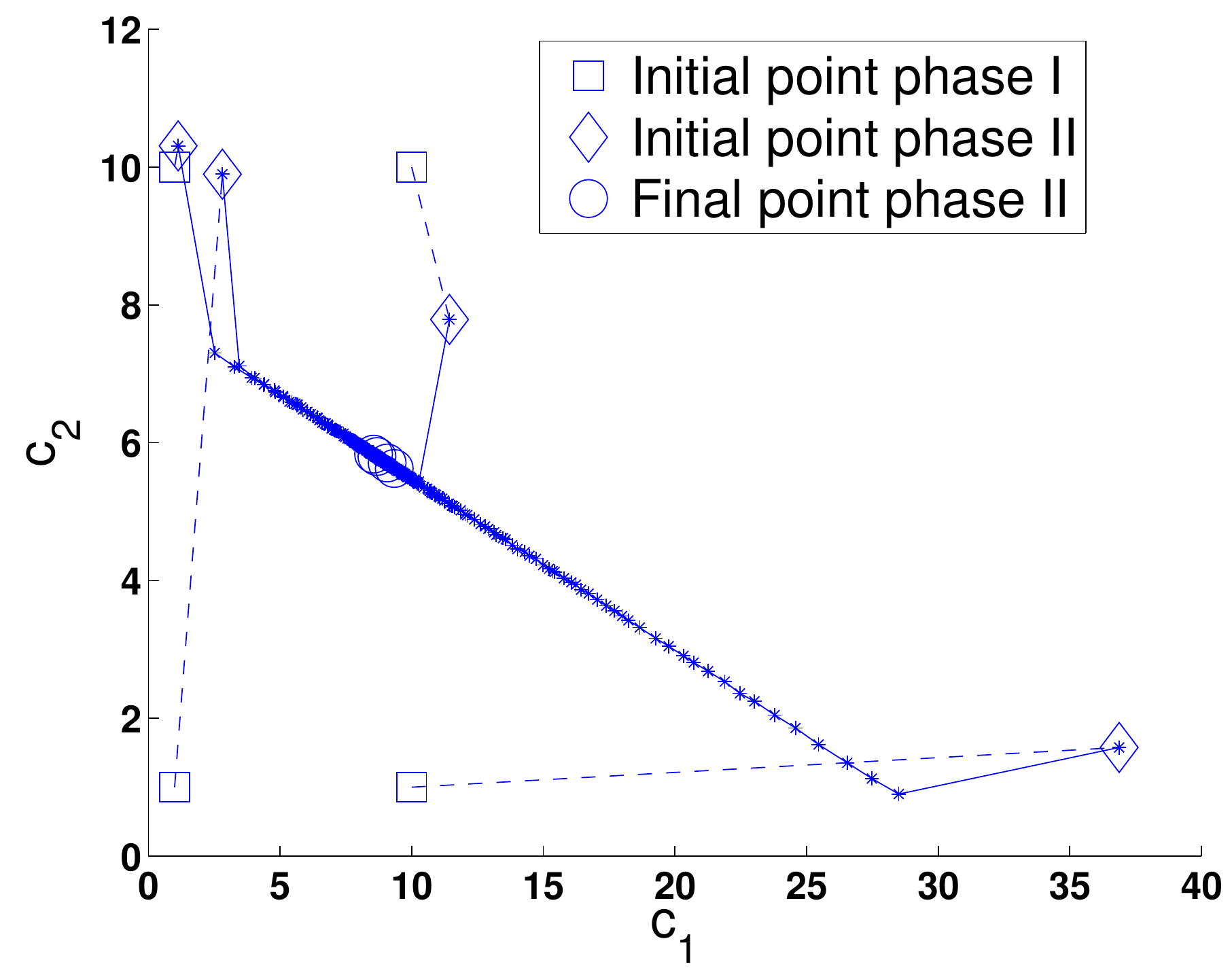}
\end{minipage}
\begin{minipage}{0.51\textwidth}
\includegraphics[width=\textwidth]{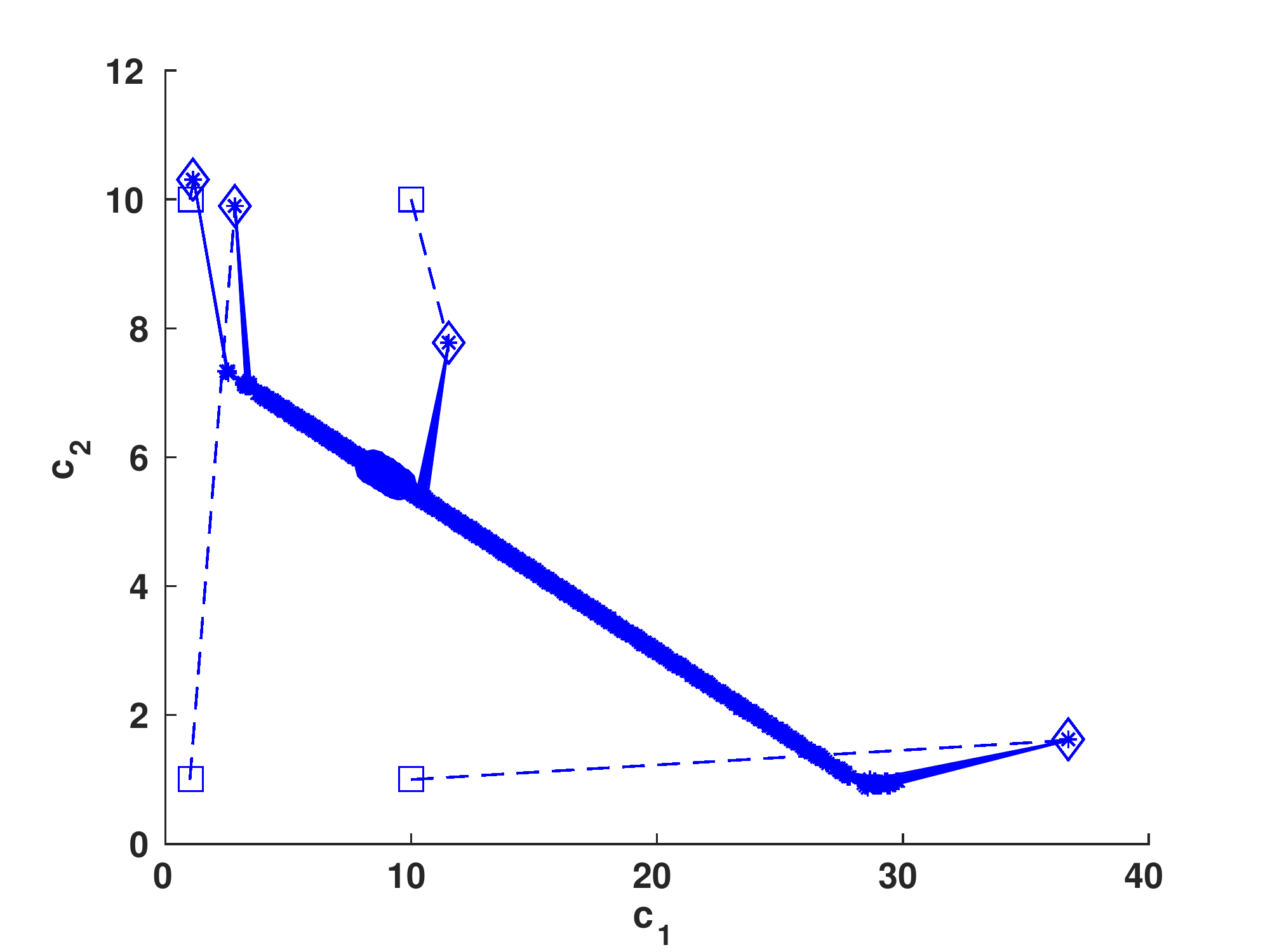}
\end{minipage}
\caption{Left: FREM estimation (phase I and phase II) for the wear example. The $N$ final values of this particular run are shown as circles. Right: We show 30 independent runs of the FREM algorithm.}
\label{fig:dec2b}
\end{figure}

We illustrate one run of our FREM algorithm in the left panel of Figure \ref{fig:dec2b} and in Table \ref{tab:dec2b}. For that run, we obtained a cluster average of 
$\hat{\theta} {=} (8.91 , 5.74)$, which corresponds to $\hat{{\theta}}_o {=} (1.5 \cdot 10^{-4} , 0.97 \cdot 10^{-4})$ in the non scaled model. 
The algorithm converged after $p^*{=}93$ iterations using 1.4 as a threshold for $\hat{R}$. 
We take that cluster average as an MLE point estimation of the unknown parameters.
%Details can be found in Table \ref{tab:dec2b} and Figure \ref{fig:dec2b}.

\begin{table}[h!]
\centering
\begin{tabular}{cccc}
$i$ & $\square {=} {\theta}_{I,i}^{(0)}$ & $\Diamond {=} {\theta}_{I\!I,i}^{(0)}$ & $\bigcirc {=} \hat{\theta}_{I\!I,i}^{(p^*)}$ \\ 
\hline 
1 &(1, 1) &   (2.81, 9.90) & (8.56, 5.83) \\ 
%\hline 
2 &(10, 1) &  (36.88, 1.58) & (9.07, 5.71) \\ 
%\hline 
3 & (1, 10) & (1.13, 10.31)  & (8.68, 5.80) \\ 
%\hline 
4 &(10, 10) & (11.44, 7.79) & (9.34, 5.62) \\ 
%\hline 
\end{tabular} 
\caption{Values computed by one run of the FREM algorithm for the wear example corresponding to the left panel of Figure \ref{fig:dec2b}.}
\label{tab:dec2b}
\end{table}

We computed an ensemble of 30 independent runs (and obtained 30 cluster averages). The result is shown in the right panel of Figure \ref{fig:dec2b}. We observe that there is a small  variability in the estimates indicating the robustness of the method. Details are shown in Table \ref{tab:dec2b_ens}.
\begin{table}[h!]
\centering
\begin{tabular}{c|cccc}
%$n$ & \multicolumn{4}{c}{$\hat{c}_n$} \\
%\cline{2-5} \rsp
& Average & Average CI at $95\%$ & Min Value & Max Value \\
\hline  \rsp
$\hat{c}_1$ &  8.94       &(8.90,  8.98) & 8.71       & 9.22  \\ 
$\hat{c}_2$ & 5.73  &(5.72, 5.74) & 5.66 & 5.79
\end{tabular} 
\caption{Values computed for an ensemble of 30 independent runs of the FREM algorithm for the wear example. 
In each run, we obtain a cluster average, $\hat{\theta}^{(i)}$, as an MLE point estimate. Define $\mathcal{C} {:=} \seqof{\hat{\theta}^{(i)}}{i=1}{30}$.
For each unknown coefficient $c_j$ in $\theta$, we show i) the average of    $\mathcal{C}$, ii) a $95\%$ confidence interval for the mean of $\mathcal{C}$, and iii) the minimum and maximum values of $\mathcal{C}$. }
\label{tab:dec2b_ens}
\end{table}

\begin{figure}[h!]
\centering
\begin{minipage}{0.49\textwidth}
\includegraphics[width=\textwidth]{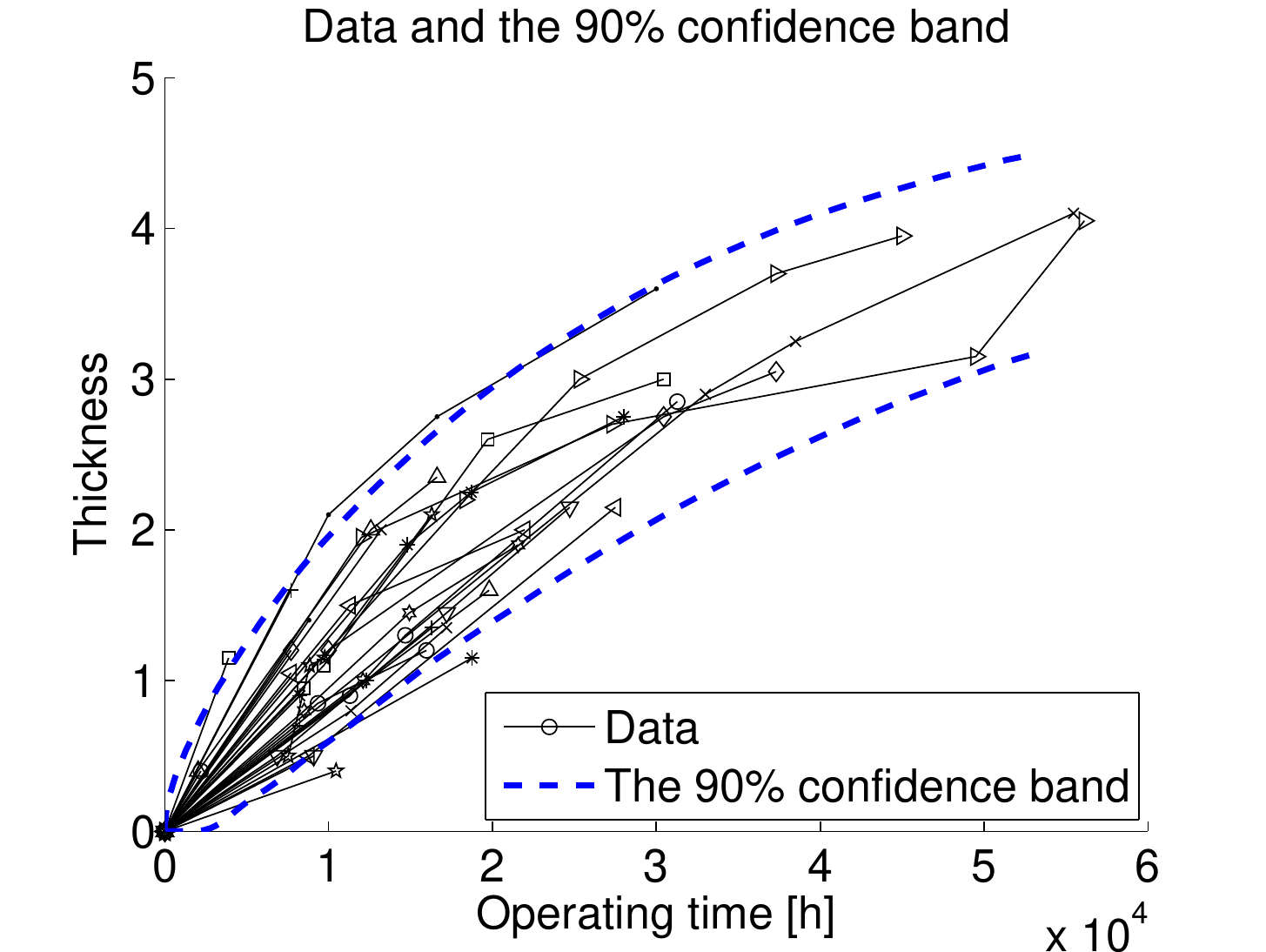}
\end{minipage}
\hfill
\begin{minipage}{0.49\textwidth}
\includegraphics[width=\textwidth]{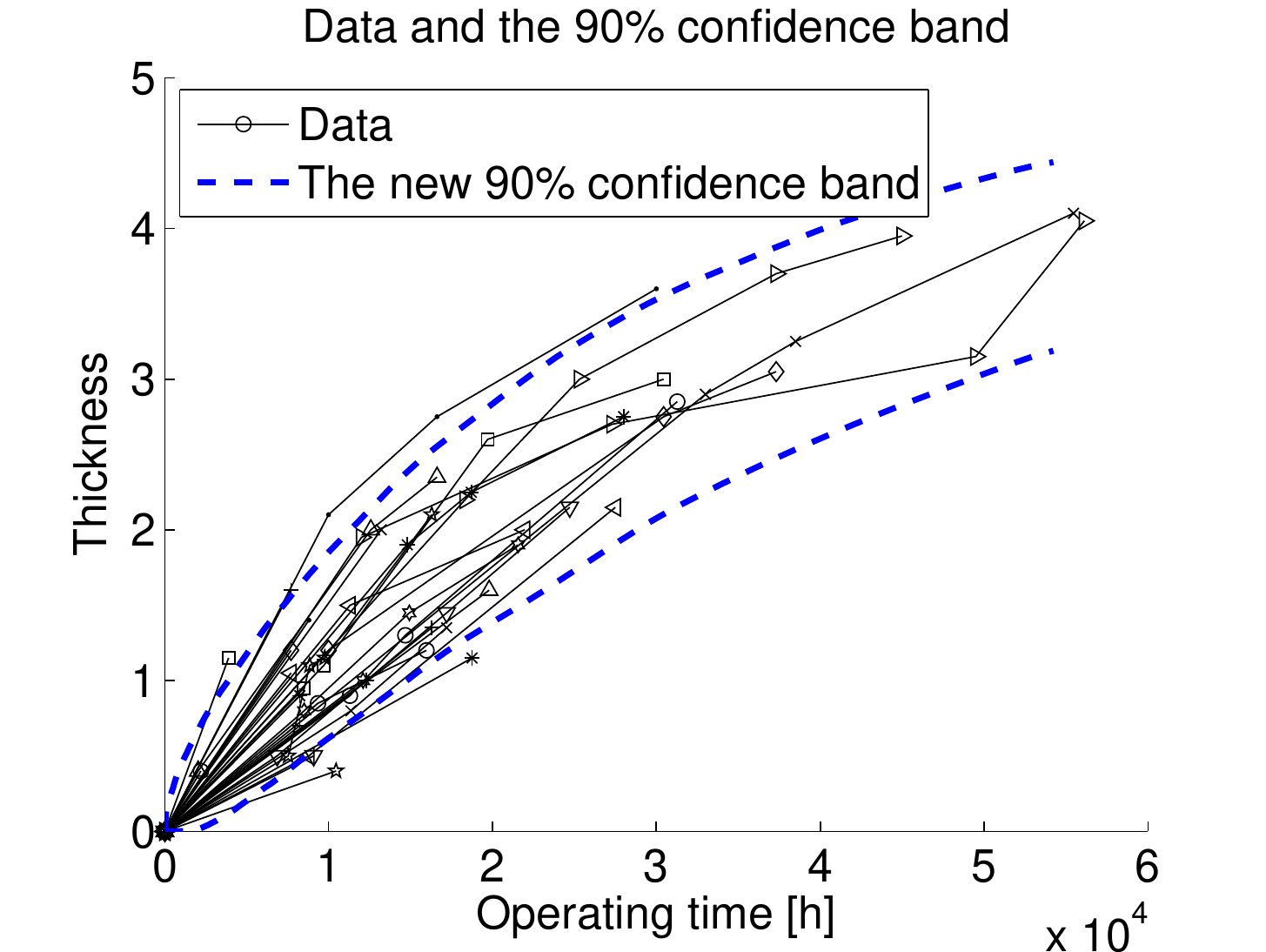}
\end{minipage}
\caption{Left: confidence band with the parameter $\tilde{\theta}$ obtained in \cite{ourInf} for the wear example. Right: the confidence band obtained with the FREM algorithm.}
\label{fig:twoCI}
\end{figure}

\begin{rem}
In this particular example, the data set was obtained using a caliper with finite precision. Therefore, our likelihood should also incorporate  the distribution of the measurement errors, which may be assumed Gaussian,  independent, and identically distributed with mean zero and variance equals to the caliper's precision. We omitted this step in our analysis for the sake of simplicity and brevity.
\end{rem}

\begin{rem}
Comparing our FREM estimate, $\hat{\hat{\theta}} {=} (1.5 \cdot 10^{-4} , 0.97 \cdot 10^{-4})$, with the value obtained in \cite{ourInf} for the same data set and the same model, ${\tilde{\theta}} {=} (0.63 \cdot 10^{-4} , 1.2 \cdot 10^{-4})$, we obtained the same scale in the coefficients and a quite similar confidence band, see Figure \ref{fig:twoCI}.
\end{rem}

\subsection{Birth-Death Process}\label{ex:bd}
This model has one species and two reaction channels:
\begin{align*}
\emptyset \xrightarrow{c_1} X,& \ \ X \xrightarrow{c_2} \emptyset
\end{align*}
described  by the stoichiometric matrix and the propensity function
\begin{align*}
\nu^T = \left( 
 \begin{array}{r}     
      1   \\
     -1  
 \end{array} 
 \right) \mbox{   and   }  a(X) = \left( \begin{array}{l}  c_1 \\ c_2 \,X \end{array} \right), \text{ respectively}\PERIOD
\end{align*}
Since we are not continuously observing the paths of $X$, 
an increment of size $k$ in the number of particles in a time interval 
$[t_1,t_2]$ may be the consequence of any combination of $n{+}k$ firings of channel 1 and $n$ firings of channel 2 in that interval.
This fact makes  the estimation of $c_1$ and $c_2$ nontrivial.

\begin{figure}[h!]
\centering
\includegraphics[scale=0.4]{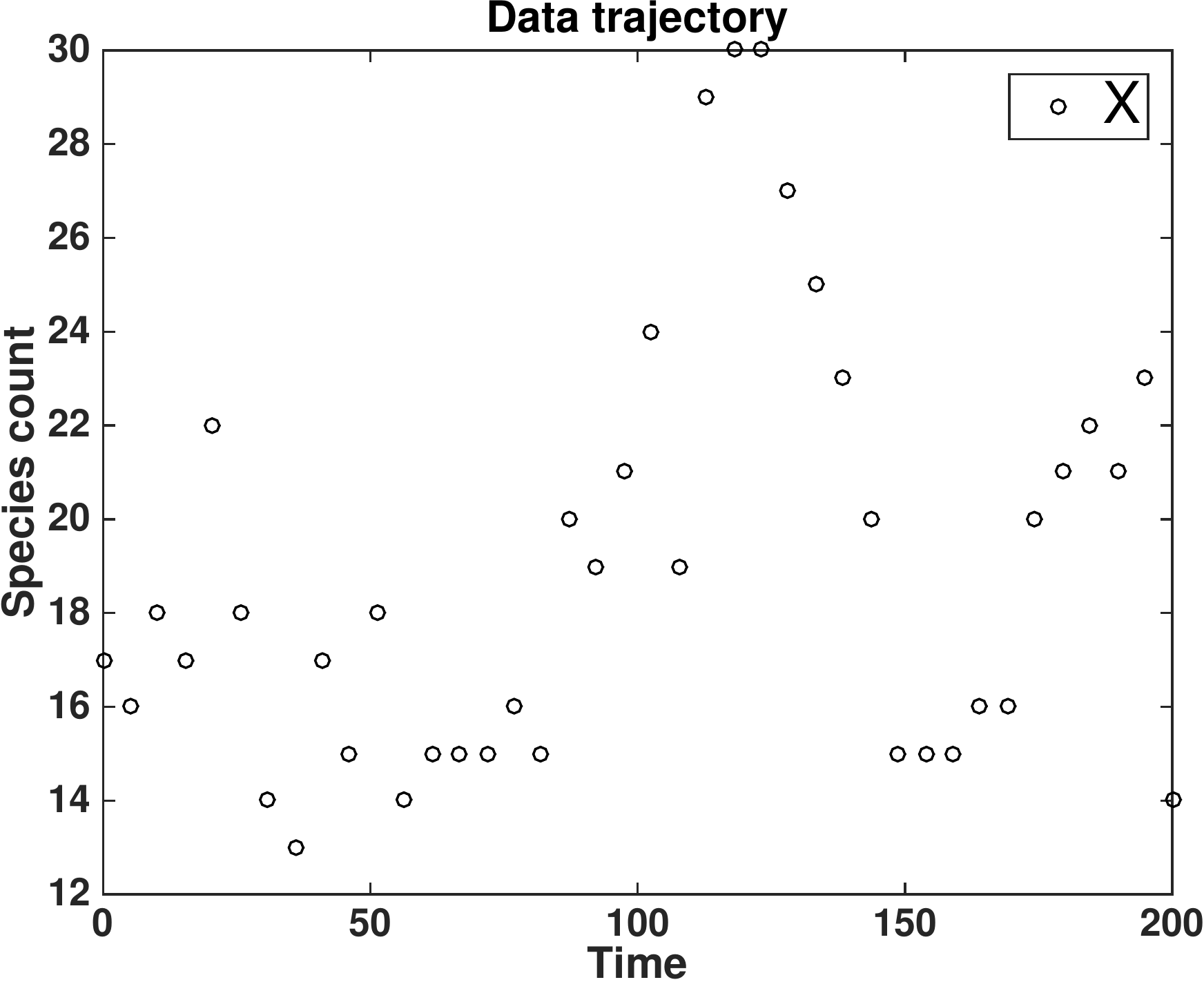}
\caption{Data trajectory for the Birth-death example. This is obtained by observing the values of an SSA path at uniform time intervals of size
$\Delta t {=}5$.}
\label{fig:dataBD}
\end{figure}

We set $X_0{=}17$, $T{=}200$ and consider synthetic data observed in uniform time intervals of size $\Delta t {=}5$. This determines a set of $41$ observations generated from a single path using the parameter $\theta_G {=}(1, 0.06)$.  The data trajectory is shown in Figure \ref{fig:dataBD}.

For this example, we ran $N{=}4$ FREM sequences starting at $\theta_{I,1}^{(0)}{=}(0.5,0.04)$, $\theta_{I,2}^{(0)}{=}(0.5,0.08)$, $\theta_{I,3}^{(0)}{=}(1.5,0.04)$, and $\theta_{I,4}^{(0)}{=}(1.5,0.08)$. Those points where chosen after a previous exploration with phase I.

\begin{figure}[h!]
\centering
\begin{minipage}{0.49\textwidth}
\includegraphics[width=\textwidth]{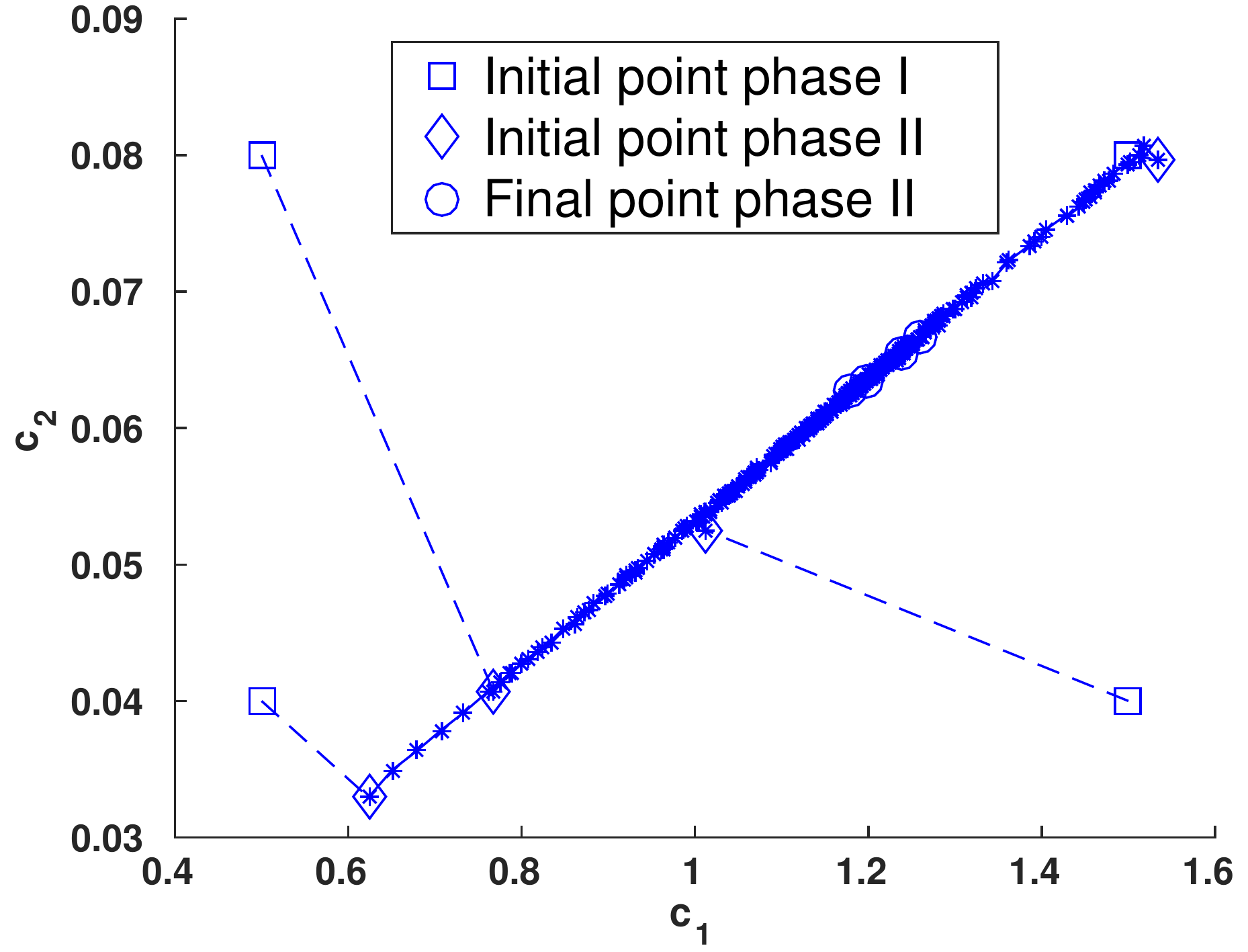}
\end{minipage}
%\hfill
\begin{minipage}{0.49\textwidth}
\includegraphics[width=\textwidth]{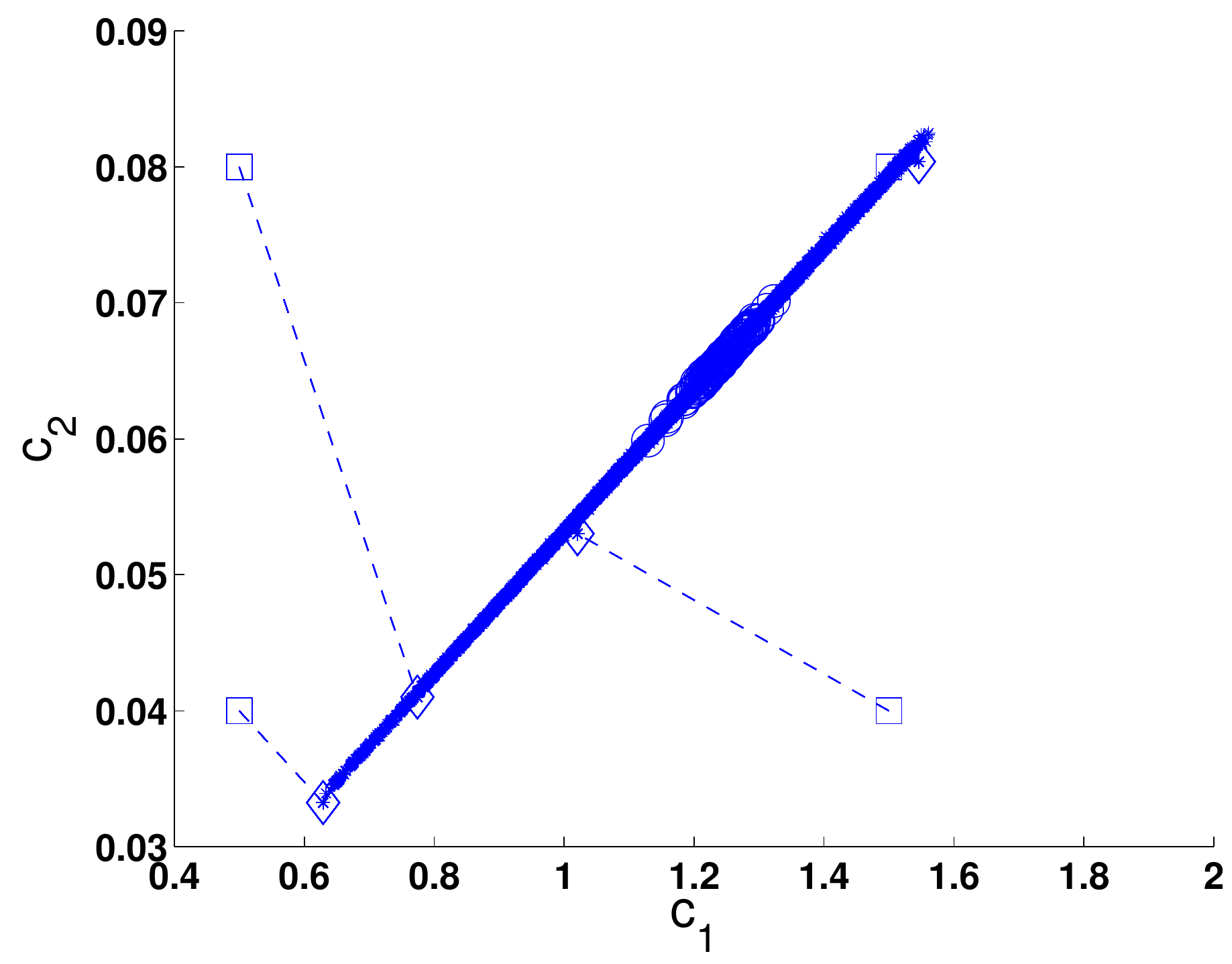}
\end{minipage}
\caption{Left: FREM estimation (phase I and phase II) for the birth-death  example. The $N$ final values of this particular run are shown as circles. Right: We show 30 independent runs of the FREM algorithm.}
\label{fig:bd}
\end{figure}

We illustrate one run of our FREM algorithm in the left panel of Figure \ref{fig:bd} and Table \ref{tab:bd}. For that run, we obtained a cluster average of $\hat{\theta} {=} (1.22,0.065)$. 
%FREM $(1.221058;0.064756)$ and the MLE is $(1.218448;0.064599)$.
The FREM algorithm took  $p^*{=}95$ iterations to converge using a threshold of 1.4 for  $\hat{R}$.
We take that cluster average as a MLE estimation of the unknown parameters. %Details can be found in Table \ref{tab:bd} and Figure \ref{fig:bd}.

\begin{table}[h!]
\centering
\begin{tabular}{cccc}
$i$ & $\square = {\theta}_{I,i}^{(0)}$ & $\Diamond = {\theta}_{I\!I,i}^{(0)}$ & $\bigcirc = \hat{\theta}_{I\!I,i}^{(p^*)}$ \\ 
\hline 
1 &(0.5, 0.04) & (6.24e-01,   3.29e-02) & (1.24e+00,   6.55e-02) \\ 
%\hline 
2 &(0.5, 0.08) & (7.68e-01,   4.07e-02) & (1.29e+00,   6.67e-02) \\ 
%\hline 
3 & (1.5, 0.04) & (1.01e+00,   5.25e-02)  & (1.18e+00   6.27e-02) \\ 
%\hline 
4 &(1.5, 0.08) & (1.53e+00,   7.97e-02) & (1.20e+00,   6.34e-02) \\ 
%\hline 
\end{tabular} 
\caption{Values computed by one run of the FREM Algorithm for the birth-death example corresponding to the left panel of Figure \ref{fig:bd}.}
\label{tab:bd}
\end{table}

We compute an ensemble of 30 independent runs (and obtained 30 cluster averages), and the result is shown in the right panel of Figure \ref{fig:bd}. We observe a moderate variability in the estimates. This may indicate that the $\hat{R}$ threshold needs to be decreased and consequently more iterations of the algorithm may be needed. Details are shown in Table \ref{tab:bd_ens}.
\begin{table}[h!]
\centering
\begin{tabular}{c|cccc}
%$n$ & \multicolumn{4}{c}{$\hat{c}_n$} \\
%\cline{2-5} \rsp
& Average & Average CI at $95\%$ & Min Value & Max Value \\
\hline  \rsp
$\hat{c}_1$ &  1.243    &(1.237, 1.249) & 1.213       &  1.284    \\ 
$\hat{c}_2$ &   0.0659  &(0.0655, 0.0663) & 0.0643 &     0.0681
\end{tabular} 
%\caption{Ensemble values computed by the FREM Algorithm for the birth-death example. In each run we obtain a cluster average. For each parameter, we show its ensemble average, a $95\%$ confidence interval for the ensemble average, and finally the minimum and maximum values recorded in the ensemble.}
\caption{Values computed for an ensemble of 30 independent runs of the FREM algorithm for the birth-death example. 
In each run, we obtain a cluster average, $\hat{\theta}^{(i)}$, as an MLE point estimate. Define $\mathcal{C} {:=} \seqof{\hat{\theta}^{(i)}}{i=1}{30}$.
For each unknown coefficient $c_j$ in $\theta$, we show i) the average of    $\mathcal{C}$, ii) a $95\%$ confidence interval for the mean of $\mathcal{C}$, and iii) the minimum and maximum values of $\mathcal{C}$. }
\label{tab:bd_ens}
\end{table}

\subsection{SIR Epidemic Model}
In this section we consider the SIR epidemic model, where $X(t)=(S(t),I(t),R(t))$ (susceptible-infected-removed individuals) and the total population is constant, $S{+}I{+}R=N$ (see \cite{SIR}). The importance of this example lies in the fact that has a nonlinear propensity function and it has two dimensions.

\begin{figure}[h!]
\centering
\includegraphics[scale=0.4]{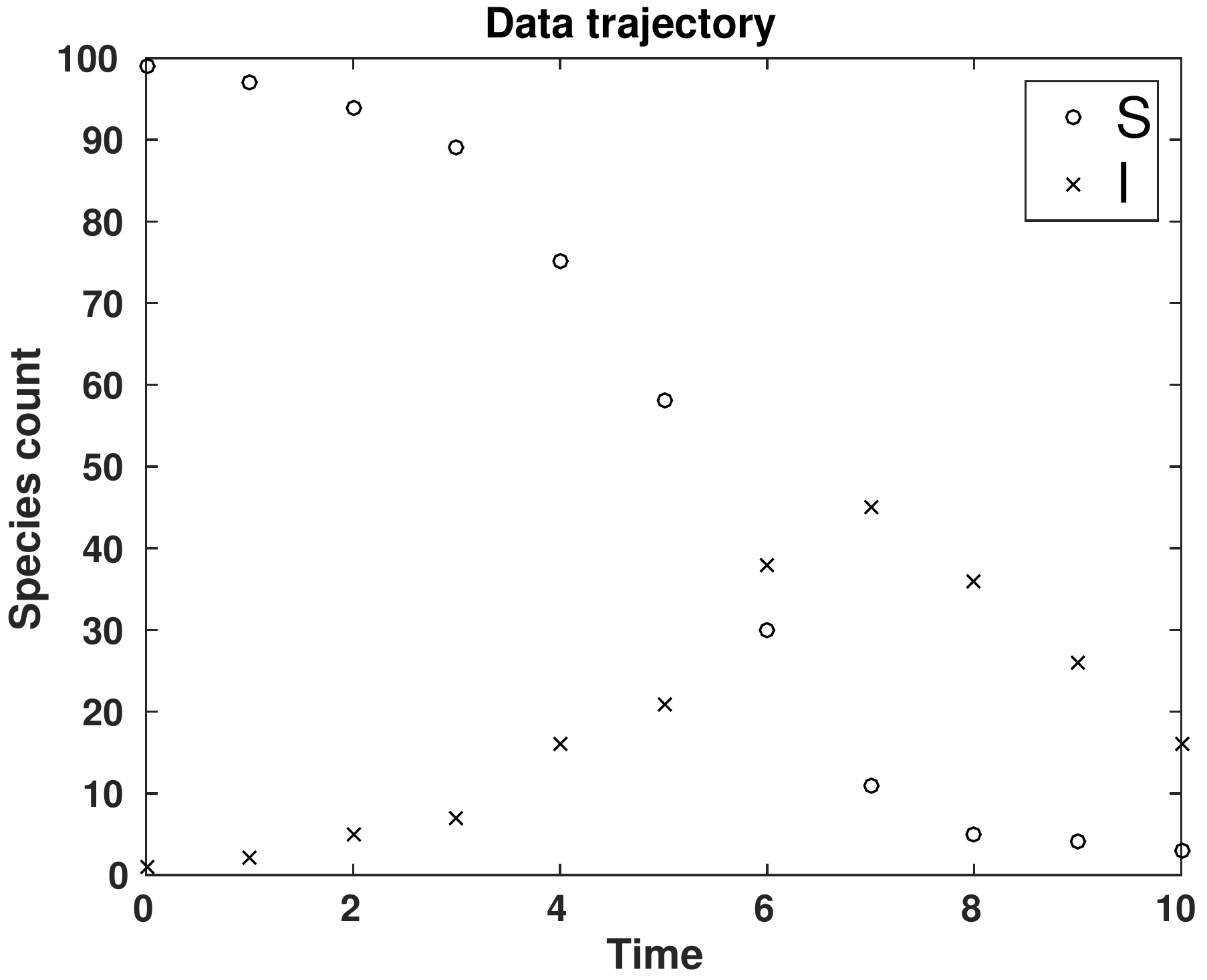}
\caption{Data trajectory for the SIR example. This is obtained by observing the values of an SSA path at uniform time intervals of size
$\Delta t {=}1$.  
}
\label{fig:dataSIR}
\end{figure}

This model has two reaction channels
\begin{align*}
S {+} I \xrightarrow{\beta} 2I, \ \ I \xrightarrow{\gamma} R 
\end{align*}
described by the stoichiometric matrix and the propensity function
\begin{align*}
\nu^T = \left( \begin{array}{rr}  -1 & 0 \\ 1 & -1 \\ 0 & 1 \end{array} \right) \mbox{   and   }  a(X) = \left( \begin{array}{l}  \beta \,S I \\ \gamma \,I \end{array} \right)\PERIOD
\end{align*}
We set $X_0{=}(300,5)$, $T{=}10$ and consider synthetic data generated using the parameters 
$\theta_G {=} (1.66, 0.44)$ by observing $X$ at uniform time intervals of size  $\Delta t {=}1$. The data trajectory is shown in Figure \ref{fig:dataSIR}. 

%%We start from the points: 
%The FREM estimation is $(1.856324;0.428335)$, and the MLE is $(1.855722;0.430352)$.
%      
%The FREM algorithm took 3 iterations to converge (minimum imposed), for a R threshold of $1.4$.
For this example we ran $N{=}4$ FREM sequences starting at $\theta_{I,1}^{(0)}{=}(0.40, 0.05)$, $\theta_{I,2}^{(0)}{=}(0.40, 1.00)$, $\theta_{I,3}^{(0)}{=}(3.00, 0.05)$, and $\theta_{I,4}^{(0)}{=}(3.00, 1.00)$. Those points where chosen after some previous exploration with  phase I.

\begin{figure}[h!]
\centering
\begin{minipage}{0.49\textwidth}
\includegraphics[width=\textwidth]{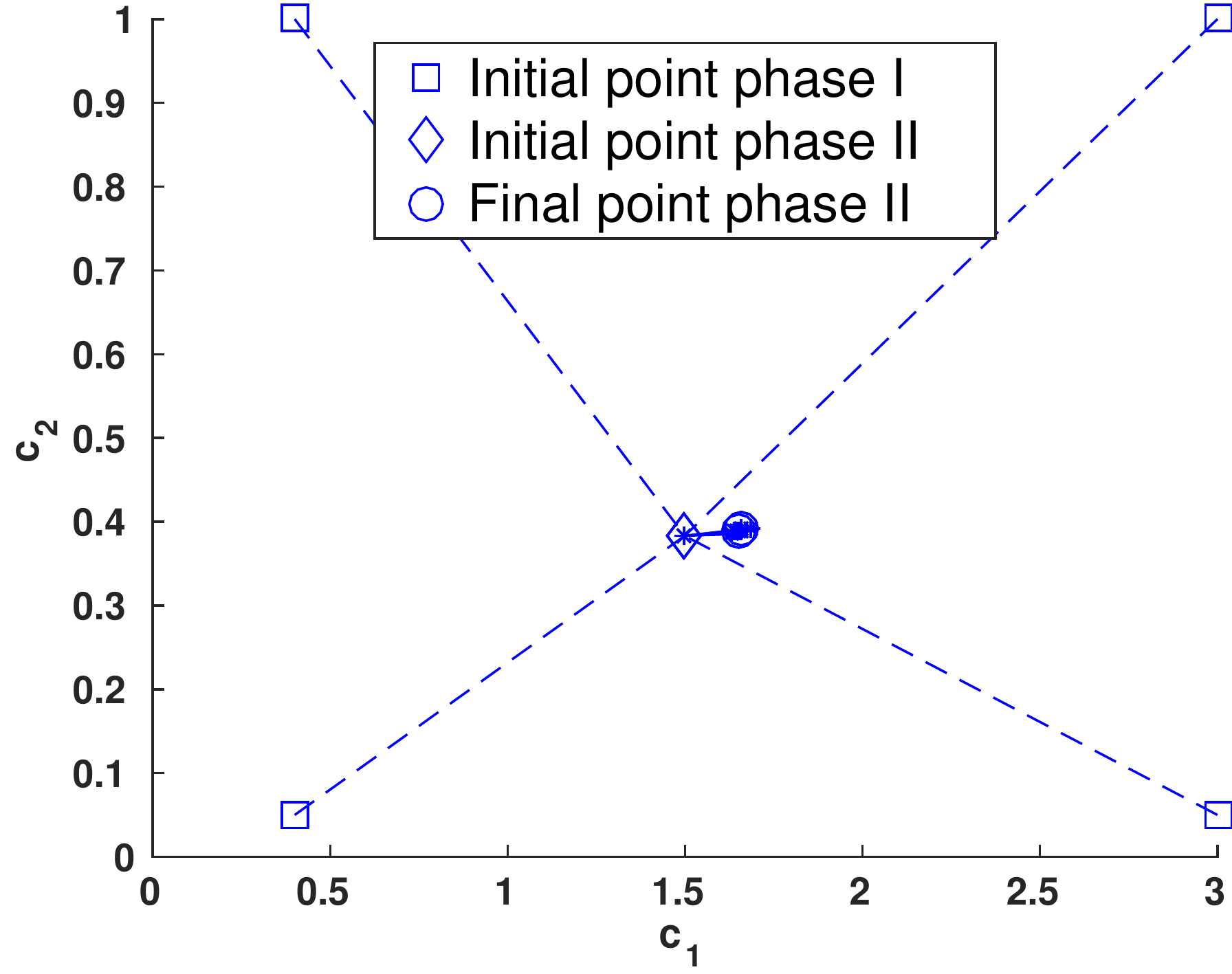}
\end{minipage}
%\hfill
\begin{minipage}{0.49\textwidth}
\includegraphics[width=\textwidth]{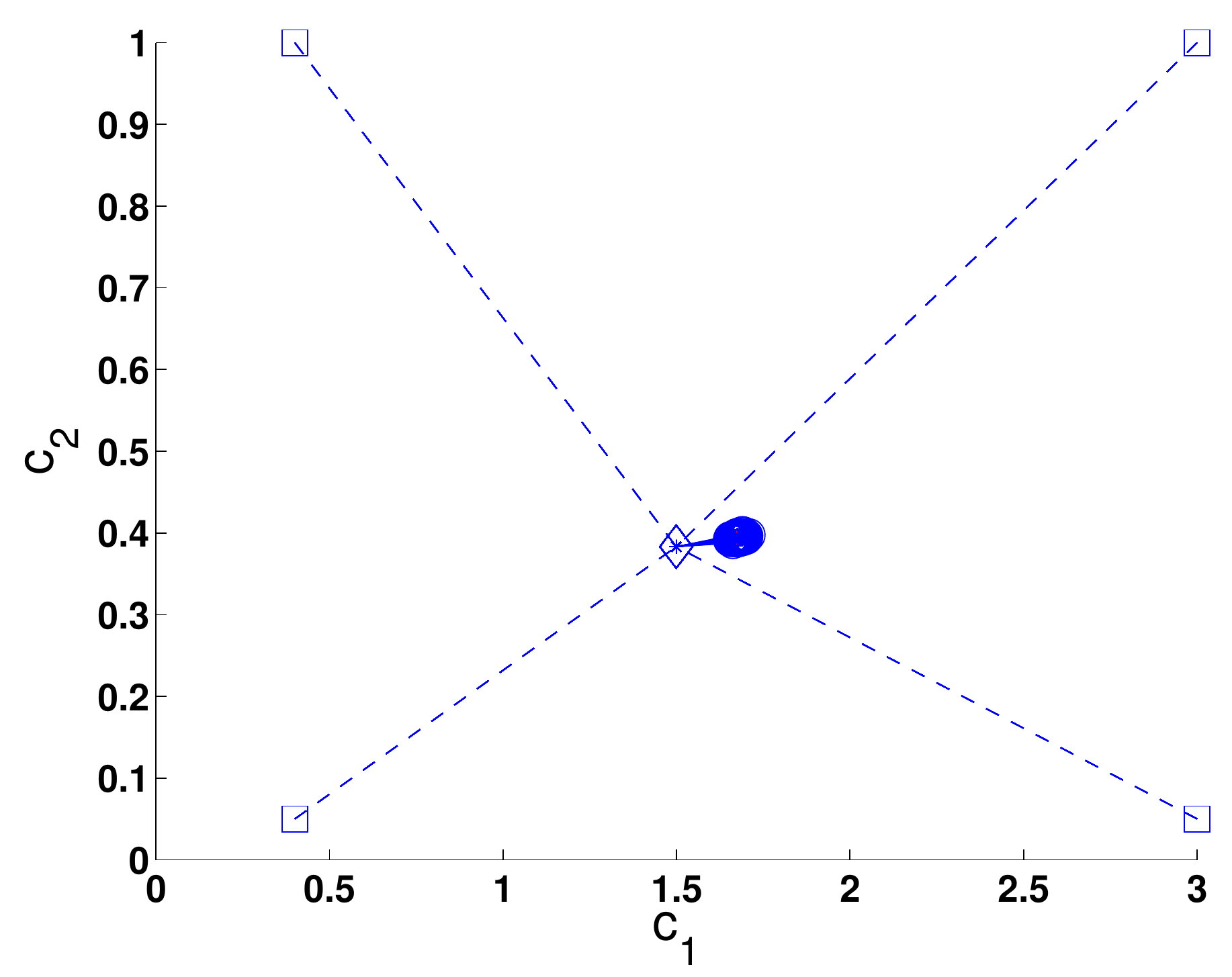}
\end{minipage}
\caption{Left: FREM estimation (phase I and phase II) for the SIR example. The $N$ final values of this particular run are shown as circles. 
In this particular case, where the results of phase I collapses to a single point, $N=4$ FREM sequences seem to be unnecessary, but we note that the $\hat R$ criterion needs at least 2 sequences.
Right: We show 30 independent runs of the FREM algorithm.}
\label{fig:sir}
\end{figure}

We illustrate one run of the FREM algorithm in the left panel of Figure \ref{fig:sir}. Our MLE point estimation is obtained as the cluster average of the values shown in Table \ref{tab:sir}, that is 
$\hat{\theta} {=} (1.65, 0.39)$. The FREM algorithm took $p^*{=}3$ iterations to converge, using 1.4 as a threshold for  $\hat{R}$. 
%We take that cluster average as a MLE estimation of the unknown parameters. 

\begin{table}[h!]
\centering
\begin{tabular}{cccc}
$i$ & $\square = {\theta}_{I,i}^{(0)}$ & $\Diamond = {\theta}_{I\!I,i}^{(0)}$ & $\bigcirc = \hat{\theta}_{I\!I,i}^{(p^*)}$ \\ 
\hline 
1 &(0.40, 0.05) &   (1.50,   0.38) & (1.65,   0.39) \\ 
%\hline 
2 &(0.40, 1.00) &  (1.50,   0.38) & (1.65,   0.39) \\ 
%\hline 
3 & (3.00, 0.05) & (1.50,   0.38)  & (1.66,   0.39) \\ 
%\hline 
4 &(3.00, 1.00) & (1.50,   0.38) & (1.66,   0.39) \\ 
%\hline 
\end{tabular} 
\caption{Values computed by one run of the FREM Algorithm for the SIR example corresponding to the left panel of Figure \ref{fig:sir}.}
\label{tab:sir}
\end{table}
%   2.9597e+00   6.6266e-01   1.8559e+00   4.2880e-01
%   2.9598e+00   6.6279e-01   1.8585e+00   4.2896e-01
%   2.9600e+00   6.6264e-01   1.8557e+00   4.2792e-01
%   2.9609e+00   6.6292e-01   1.8553e+00   4.2766e-01
%\subsubsection{SIR with small number of initial infected individuals}
%We set $X_0{=}(300,1,0)$, $T{=}10$ and we use a uniform grid of measurements of size $\Delta t {=}$.
%In this case we have two equilibrium states, which reflects the fact that the epidemic outbreak can occur or not, because at time zero there is only 1 infected individual.
%
%%\begin{figure}[h!]
%%\centering
%%\includegraphics[scale=0.43]{plots/numexamples/\SIR/SIR_T10_trajectories.pdf}
%%\caption{10 SSA trajectories, its average, and the mean field for the SIR model.}
%%\label{fig:rev-trajs}
%%\end{figure}
%
%\subsection{SIR with demography}
%

We computed an ensemble of 30 independent runs (and obtained 30 cluster averages);  results are shown in the right panel of Figure \ref{fig:sir}. We observe a very small variability in our estimates; details are shown in Table \ref{tab:sir_ens}.
\begin{table}[h!]
\centering
\begin{tabular}{c|cccc}
%$n$ & \multicolumn{4}{c}{$\hat{c}_n$} \\
%\cline{2-5} \rsp
& Average & Average CI at $95\%$ & Min Value & Max Value \\
\hline  \rsp
$\hat{c}_1$ &  1.6784    &(1.6764, 1.6804) & 1.6648   & 1.6891 \\ 
$\hat{c}_2$ & 0.3942  &(0.3939, 0.3945) & 0.3920 & 0.3956
\end{tabular} 
%\caption{Ensemble values computed by the FREM Algorithm for the SIR example. In each run we obtain a cluster average. For each parameter, we show its ensemble average, a $95\%$ confidence interval for the ensemble average, and finally the minimum and maximum values recorded in the ensemble.}
\caption{Values computed for an ensemble of 30 independent runs of the FREM algorithm for the SIR example. 
In each run, we obtain a cluster average, $\hat{\theta}^{(i)}$, as an MLE point estimate. Define $\mathcal{C} {:=} \seqof{\hat{\theta}^{(i)}}{i=1}{30}$.
For each unknown coefficient $c_j$ in $\theta$, we show i) the average of    $\mathcal{C}$, ii) a $95\%$ confidence interval for the mean of $\mathcal{C}$, and iii) the minimum and maximum values of $\mathcal{C}$. }
\label{tab:sir_ens}
\end{table}

\subsection{Auto-Regulatory Gene Network}
\begin{figure}[h!]
\centering
\includegraphics[scale=0.4]{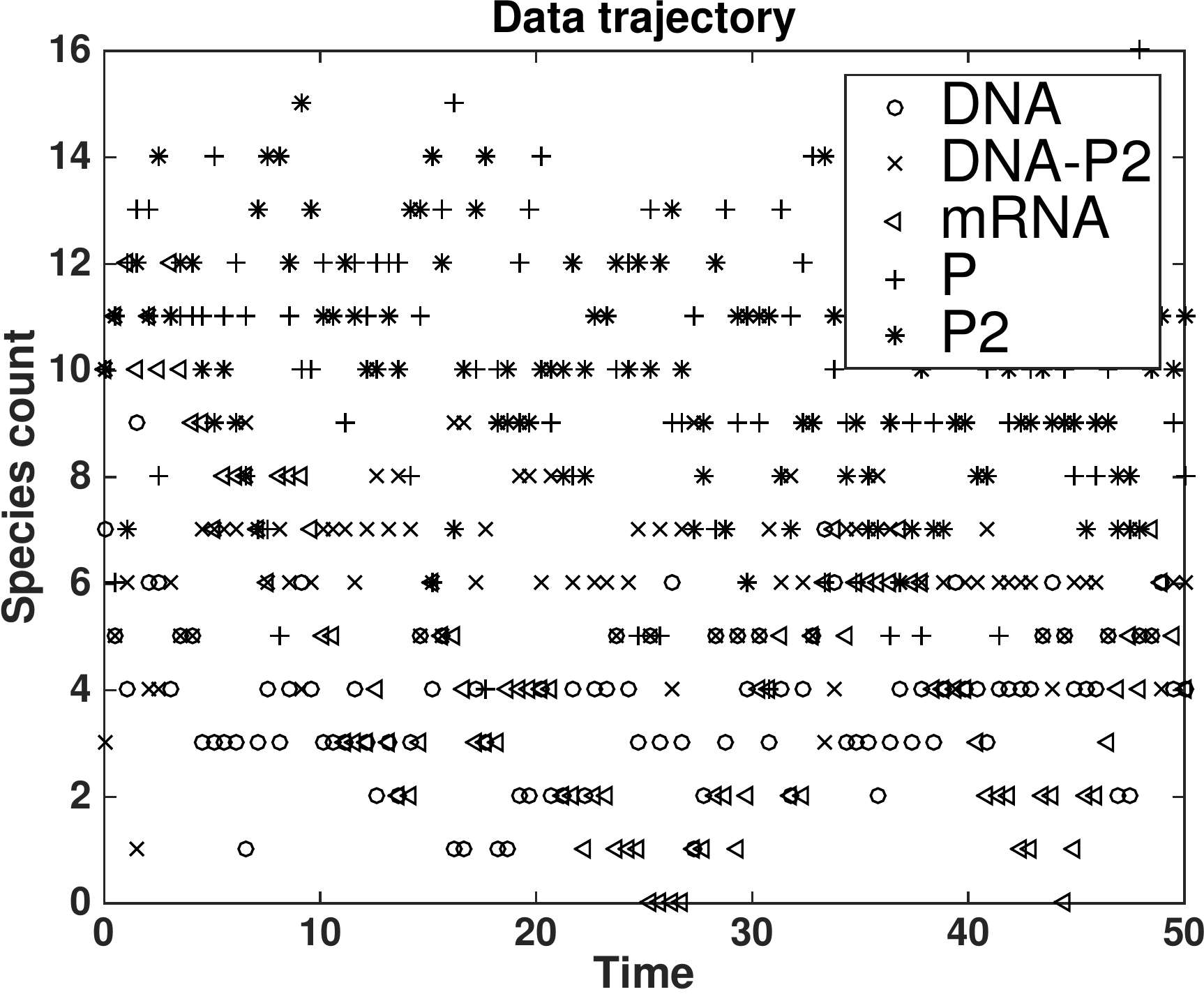}
\caption{Data trajectory for the auto-regulatory gene network example obtained by observing the values of an SSA path at uniform time intervals of size
$\Delta t {=}\frac{1}{2}$.  
}
\label{fig:dataARG}
\end{figure}

The following model, taken from \cite{daigle2012accelerated}, has eight reaction channels and five species,
\begin{align*}
DNA + P_2 &\xrightarrow{c_1} DNA{-}P_2, \ \ &
DNA{-}P_2 &\xrightarrow{c_2}  DNA + P_2\\
DNA &\xrightarrow{c_3} DNA + mRNA, \ \ &
mRNA  &\xrightarrow{c_4} \emptyset\\
P+P &\xrightarrow{c_5} P_2, \ \ &
P_2 &\xrightarrow{c_6} P+P\\
mRNA &\xrightarrow{c_7} mRNA +P, \ \ &
P&\xrightarrow{c_8}\emptyset
\end{align*}
 is described respectively by the stoichiometric matrix and the propensity function
\begin{align*}
\nu^T = \left( 
 \begin{array}{rrrrr}     
   -1   &  1  &    0  &  0  &  -1 \\
   1    &  -1  &    0 &  0   &  1 \\
    0   &   0  &  1   & 0   &   0\\
    0   &   0 &  -1  &  0 &  0\\
    0  &   0  &   0  & -2 & 1\\
    0 &   0 &  0 &  2 & -1\\
    0 &  0 & 0 & 1 & 0 \\
    0 & 0 & 0 & -1 & 0
 \end{array} 
 \right) \mbox{   and   }  a(X) = \left( \begin{array}{l}  c_1\,DNA{-}P_2 \\ c_2\,DNA \cdot P_2 \\ c_3\,DNA\\ c_4\,mRNA \\ 
 c_5 \,P(P{-}1)\\ c_6\,P_2\\c_7\,mRNA\\c_8\,P  \end{array} \right)\PERIOD
\end{align*}
Quoting \cite{daigle2012accelerated}, ``$DNA$, $P$, $P_2$, and $mRNA$ represent $DNA$ promoters,
protein gene products, protein dimers, and messenger $RNA$ molecules, respectively.'' 
This model has been selected to test the robustness of our FREM algorithm to deal with several dimensions and several reactions.
Following cited works, we also set the initial state of the system at $$X_0 = (DNA,DNA{-}P_2,mRNA, P, P_2) = (7, 3, 10, 10, 10),$$ and run the system to the final time $T = 50$. Synthetic data is gathered by observing a single trajectory generated using 
$\theta_G = (0.1, 0.7, 0.35, 0.3, 0.1, 0.9,0.2, 0.1)$ at uniform time intervals of size  $\Delta t {=}\frac{1}{2}$. The data trajectory is shown in Figure \ref{fig:dataARG}. 
For this example we ran $N{=}2$ FREM sequences starting at
$\theta_{I,1}^{(0)}=0.1\,v$ and
$\theta_{I,2}^{(0)}=0.5\,v$, respectively, where $v$ is the vector of $\rset^8$ with all its components equal to one. 
%\begin{align*}
%\theta_{I,1}^{(0)}&{=}(0.1,0.1,0.1,0.1,0.1,0.1,0.1,0.1) \text{ and}\\
%\theta_{I,2}^{(0)}&{=}(0.5,0.5,0.5,0.5,0.5,0.5,0.5,0.5)
%\end{align*}
%Our FREM algorithm estimation gave us a cluster average of 
%$$\hat \theta =(0.107,
%    0.649,
%    0.337,
%    0.319,
%    0.087,
%    0.835,
%    0.053,
%    0.025),$$ which is in the range of the components of $\theta_G$ and seems to be a satisfactory estimation of $\theta$. 
    
The FREM algorithm took, on average, $p^*{=}169$ iterations to converge, taking 2 days in our workstation configuration: a 12 core Intel GLNXA64 architecture and MATLAB version R2014a.

We computed an ensemble of 10 independent runs and obtained 10 cluster averages. We observe very small variability. Details are shown in Table \ref{tab:arg_ens}.
\begin{table}[h!]
\centering
\begin{tabular}{c|cccc}
%$n$ & \multicolumn{4}{c}{$\hat{c}_n$} \\
%\cline{2-5} \rsp
& Average & Average CI at $95\%$ & Min Value & Max Value \\
\hline  \rsp
$\hat{c}_1$ &  0.1011   &(0.1001, 0.1021) &0.0984   & 0.1033    \\ 
$\hat{c}_2$ &0.6207  &(0.6135, 0.6279) & 0.6005 & 0.6328 \\
$\hat{c}_3$ &  0.3398   &(0.3380, 0.3416) & 0.3358   & 0.3441    \\ 
$\hat{c}_4$ &  0.3182   &(0.3166, 0.3198) & 0.3139   & 0.3213    \\ 
$\hat{c}_5$ &  0.0637   &(0.0622, 0.0652) & 0.0595   & 0.0687    \\ 
$\hat{c}_6$ &  0.5891   &(0.5742, 0.6040) & 0.5485   &  0.6357    \\ 
$\hat{c}_7$ &  0.1444   &(0.1426, 0.1462) & 0.1392   & 0.1483    \\ 
$\hat{c}_8$ &  0.0630   &(0.0623, 0.0637) & 0.0618   & 0.0652    
\end{tabular} 
\caption{Values computed for an ensemble of 10 independent runs of the FREM algorithm for the auto-regulatory gene network example. 
In each run, we obtain a cluster average, $\hat{\theta}^{(i)}$, as an MLE point estimate. Define $\mathcal{C} {:=} \seqof{\hat{\theta}^{(i)}}{i=1}{10}$.
For each unknown coefficient $c_j$ in $\theta$, we show i) the average of    $\mathcal{C}$, ii) a $95\%$ confidence interval for the mean of $\mathcal{C}$, and iii) the minimum and maximum values of $\mathcal{C}$. }
\label{tab:arg_ens}
\end{table}        
%
%\subsection{Gene Transcription and Translation \cite{Anderson2012}} This model has five reactions,
%\begin{align*}
%\emptyset \xrightarrow{c_1} R,& \ \ R \xrightarrow{c_2} R+P \\
%2P \xrightarrow{c_3}  D,& \ \ R \xrightarrow{c_4} \emptyset \\
%P \xrightarrow{c_5} \emptyset 
%\end{align*}
%described respectively by the stoichiometric matrix and the propensity function
%\begin{align*}
%\nu = \left( 
% \begin{array}{cccc}     1    & 0 &    0 \\
%     0    & 1 &    0 \\
%     0   & -2 &    1 \\
%    -1    & 0 &    0 \\
%     0   & -1 &    0 
% \end{array} 
% \right) \mbox{   and   }  a(X) = \left( \begin{array}{c}  c_1 \\ c_2 R \\ c_3 P(P{-}1)\\ c_4R \\ c_5 P  \end{array} \right)\COMMA
%\end{align*}
%where $X(t)=(R(t),P(t),D(t))$, and $c_1 {=} 25$, $c_2 {=}10^3$, $c_3{=}0.001$, $c_4{=}0.1$, and $c_5{=}1$. We set $X_0{=}(0,0,0)$, $T{=}1$ and we use a uniform grid of measurements of size $\Delta t {=}$.
%
\begin{rem}
Observe that in the examples where the stoichiometric vectors are linearly dependent, the results of the phase I, ${\theta}_{I\!I,i}^{(0)}$, $i=1,2,3,4$, lies in a hyperplane that reflects a certain amount of indifference in the coefficient estimations. This does not happen in the SIR example where all the estimations in phase I are essentially the same.
\end{rem}

% flatex input end: [num_examples]

%\clearpage

%\clearpage
\section{Conclusions}
% flatex input: [conclusions.tex]
\label{conclusions}
In this work, we addressed the problem of efficiently computing approximations of expectations of functionals of bridges in the context of stochastic reaction networks by extending the forward-reverse technique developed by Bayer and Schoenmakers in \cite{Bayer}.
We also showed how to apply this technique to the statistical problem of inferring the set of coefficients of the propensity functions. 
We presented a two-phase approach, namely the Forward-Reverse Expectation-Maximization (FREM) algorithm, in which the first phase, based on reaction-rate ODEs is deterministic and  is intended to provide a starting point that reduces the computational work of the second phase, namely, the Monte Carlo EM Algorithm. 
Our novel algorithm for generating bridges provides a clear advantage over shooting methods and  methods based on acceptance rejection techniques. 
Our work  is illustrated with numerical examples.
In the future, we plan to incorporate higher-order kernels and multilevel Monte Carlo methods in the FREM algorithm.

% flatex input end: [conclusions.tex]

%\clearpage

\section*{Acknowledgments}
% flatex input: [acknowledgements.tex]
\thx{
 %\subsection*{Acknowledgements}
A. Moraes, R. Tempone and P. Vilanova are members of the KAUST SRI Center for
Uncertainty Quantification  at the Computer, Electrical and Mathematical Sciences and Engineering Division at King Abdullah University of Science and Technology (KAUST).
}
%%% Local Variables: 
%%% mode: latex
%%% TeX-master: "main"
%%% End: 

% flatex input end: [acknowledgements.tex]

%\clearpage

\newpage
\appendix
% flatex input: [appendix.tex]
\section{Algorithms}\label{sec:algorithms}
\begin{algorithm}[h!]
\caption{The F-R (forward-reverse) path generation algorithm in  phase II, for a given time interval, $[s,t]$. Inputs: the initial sample size, $M_0$, the coefficient of variation threshold, $cv_0$, the initial time, $s$, the final time, $t$, %the intermediate time factor, $t^*_f$, 
the initial observed state, $x(s)$, and the final observed state, $x(t)$. %A positive constant $c$.
Outputs: a sequence of the number of times that a reaction channel fired in the given  time interval, $\seqof{\seqof{r_{j,l}}{j=1}{J}}{l=1}{L}$, a sequence of forward Euler values for the given time interval, $\seqof{\seqof{u_{j,l}}{j=1}{J}}{l=1}{L}$ and a sequence of kernel weights for the given time interval, $\seqof{\seqof{w_{j,l}}{j=1}{J}}{l=1}{L}$.
Notes:
%based on $\expt{\tau_{\ssa}(\BX)}=1/a_0(\BX)$ and $\tau_{Ch}(\BX,\delta)$, this algorithm adaptively selects which  method to use: SSA or TL.  
Here %$D(t)$ is the data process value at time $t$, 
$V_d$ is the volume of a $d$ dimensional unit sphere, $\tilde{X}^{(f)}_{\cdot,\cdot,n}$ is the sampled forward process value at time $t_n^{(f)}$, $\tilde{X}^{(b)}_{\cdot,\cdot,n'}$ is the sampled reverse process at time $t_{n'}^{(b)}$, $\kappa_\delta$ is the Kronecker delta kernel and $\kappa_e$ is the Epanechnikov kernel, $L$ is the number of joined F-R paths in the time interval $[s,t]$, where $0 \leq L\leq \tilde{M}^2$. Finally, $0<\gamma<1 $ and $C_L$ is an integer greater than 1 (in our examples we use 2). %We assume that $L$ is approximately $\tilde{M}$.
}
\label{alg:fr_path}
\begin{algorithmic}[1]
	%\REQUIRE $a_0 \leftarrow \sum_{j=1}^J a_j > 0$
	\STATE $\tilde{M} \leftarrow 1$
	\STATE $M \leftarrow M_0$
	    %\STATE $T_k \leftarrow t_k-s_k$
		%\STATE $t^* \leftarrow t^*_f\,T_k$
		\STATE $t^* \leftarrow \frac{1}{2}(t-s)$
		\WHILE {$cv \geq cv_0$} %{$runs \leq $ MAX\_RUNS \AND $cv \geq cv_0$}
			\FOR {$m=\tilde{M}$ \TO $\tilde{M}{+}M{-}1$} 
				\STATE $(\seqof{\tilde{X}_{\cdot,m,n}^{(f)},t^{(f)}_{m,n}}{n=1}{N(m)},\seqof{r^{(f)}_{j,m}}{j=1}{J}) \leftarrow$ FW path from $s$ to $t^*$ starting at $x(s)$
				\STATE $u_{j,m}^{(f)} \leftarrow (t^{(f)}_{m,n+1}-t^{(f)}_{m,n}) g_j(\tilde{X}_{\cdot,m,n}^{(f)})$
				\STATE $(\seqof{\tilde{X}_{\cdot,m,n'}^{(b)},t^{(b)}_{m,n'}}{n=1}{N'(m)},\seqof{r^{(b)}_{j,m}}{j=1}{J}) \leftarrow$ RV path from $t$ to $t^*$ starting at $x(t)$\!\!
				
				\STATE $u_{j,m}^{(b)} \leftarrow (t^{(b)}_{m,n'+1}-t^{(b)}_{m,n'}) g_j(\tilde{X}_{\cdot,m,n'+1}^{(b)})$ 
			\ENDFOR
			\STATE $\seqof{u_{\cdot,l},r_{\cdot,l},w_{\cdot,l}}{l=1}{L} \leftarrow $ join F-R paths $(\tilde{X}_{\cdot,\cdot}^{(f,b)}(t^*),\seqof{r^{(f,b)}_{j,\cdot}}{j=1}{J},\seqof{\alpha_{j,\cdot}^{(f,b)}}{j=1}{J},\kappa_\delta)$
			\STATE $\hphantom{\seqof{f,r,w}{l=1}{L} \leftarrow}$ Here, $\alpha_{j,l} = \alpha^{(f)}_{j,m}+\alpha^{(b)}_{j,m}$ s.t. $m\in \{1,2,...,\tilde{M}\}$ and
			\STATE  $\hphantom{\seqof{f,r,w}{l=1}{L} \leftarrow and} $ $\kappa_\delta(\tilde{X}^{(f)}_{\cdot,m}(t^*),\tilde{X}^{(b)}_{\cdot,m}(t^*))>0$. Similarly for $r_{j,l}$.
%			\STATE $\hphantom{\seqof{f,r,w}{l=1}{L} \leftarrow and}$ 
			\IF {$L < \lceil \gamma \tilde{M} \rceil$}
			\STATE $\Sigma \leftarrow$ covariance matrix of $(\tilde{X}_{\cdot,m}^{(f)}(t^*),\tilde{X}_{\cdot,m}^{(b)}(t^*))$
			\STATE $\Sigma \leftarrow \Sigma + c\,diag(\Sigma)$, where $c$ is a positive constant.
			\IF {$\Sigma^{-1/2}$ \NOT singular}
				%\STATE $V_d \leftarrow $ (volume of a $d$ dimensional unit sphere)
				\STATE $H \leftarrow \frac{1}{3} \Sigma^{-1/2}(\frac{\tilde{M}}{V_d})^{1/d}$
				\STATE $\zeta \leftarrow 1$
				\REPEAT 
				\STATE $\tilde{Y}_{\cdot,m}^{(f)}(t^*) \leftarrow  \zeta H \tilde{X}_{\cdot,m}^{(f)}(t^*)$
				\STATE $\tilde{Y}_{\cdot,m}^{(b)}(t^*) \leftarrow  \zeta H \tilde{X}_{\cdot,m}^{(b)}(t^*)$
				\STATE $\seqof{u_{\cdot,l},r_{\cdot,l},w_{\cdot,l}}{l=1}{L} \leftarrow $ join F-R paths $(\tilde{Y}_{\cdot,\cdot}^{(f,b)}(t^*),\seqof{r^{(f,b)}_{j,\cdot}}{j=1}{J},\seqof{\alpha_{j,\cdot}^{(f,b)}}{j=1}{J},\kappa_e)$ \!\!
				
				\STATE $\zeta \leftarrow 1.5\zeta$
				\UNTIL{$L \leq C_L \tilde{M}$}
				%\STATE $\hphantom{\seqof{f,r,w}{l=1}{L} \leftarrow}$ Here, $\kappa_e (x,y) := \prod_{i=1}^d(1-(x_i{-}y_i)^2)\indicator{|x_i{-}y_i|<1}$
			\ENDIF
			\ENDIF
			\STATE compute the coefficient of variation of $\seqof{u_{\cdot,l}}{l=1}{L}$ and $\seqof{r_{\cdot,l}}{l=1}{L}$ (see section \ref{sec:Mk})
			\STATE $\tilde{M} \leftarrow \tilde{M} + M$
			\STATE $M \leftarrow 2M$
		\ENDWHILE
%	\IF{$K_1/a_0 < T_0-t$} 
%		\STATE  $\tau_{Ch}\leftarrow $  Algorithm \ref{alg:Ch} 
%		\IF{$\tau_{Ch} < K_2(\bar{X}(t),\delta)/a_0 $}
%			\RETURN $(\ssa, \tau_{\ssa})$
%		\ELSE
%			\RETURN $(TL, \tau_{Ch})$
%		\ENDIF
%	\ELSE
%		\RETURN $(\ssa, \tau_{\ssa})$
%	\ENDIF
\end{algorithmic}
\end{algorithm}

\begin{algorithm}[h!]
\caption{The F-R path-join algorithm for phase II. Inputs: a sequence of forward-backward samples for the time interval $[s,t]$ evaluated at the intermediate time, $t^*$, $\tilde{X}_{\cdot,\cdot}^{(f,b)}(t^*)$, a sequence of the number of times that a reaction channel fired in the \emph{forward} interval $[s,t^*]$ and in the \emph{reverse} interval $[t^*,t]$, $r^{(f,b)}_{\cdot,\cdot}$, the sequence of forward Euler values  for each reaction channel for the \emph{forward} interval $[s,t^*]$ and for the \emph{backward} interval $[t^*,t]$, $u^{(f,b)}_{\cdot,\cdot}$, and the kernel $\kappa$.
Outputs: the number of joined paths, $L$, a sequence of the number of times that a reaction channel fired in the interval $[s,t]$, $\seqof{\seqof{r_{j,l}}{j=1}{J}}{l=1}{L}$, the sequence of forward Euler values  for each reaction channel for the interval $[s,t]$, $\seqof{\seqof{u_{j,l}}{j=1}{J}}{l=1}{L}$ and the sequence of kernel weights for the interval $[s,t]$, $\seqof{\seqof{w_{j,l}}{j=1}{J}}{l=1}{L}$.
Notes: $S$ is a two dimensional sparse matrix of size $C \times \tilde{M}$. 
}
\label{alg:fr_join}
\begin{algorithmic}[1]
\STATE $L\leftarrow 0$
\FOR {$i=1$ \TO d}
\STATE $A_i \leftarrow \min_{m} \lfloor \tilde{X}_{i,m}^{(f,b)}(t^*) \rfloor$
\STATE $B_i \leftarrow  \max_{m }\lceil \tilde{X}_{i,m}^{(f,b)}(t^*) \rceil$
\STATE $E_i \leftarrow 1+B_i-A_i$
\ENDFOR
\FOR {$m=1$ \TO $\tilde{M}$}
\STATE $p_i \leftarrow 1+\lceil \tilde{X}_{i,m}^{(f)}(t^*) \rceil-A_i$
\STATE $c \leftarrow \text{convert}(p,E)$ (converts $d$ dimensional address to $\{1,...,C\}$)
\STATE $S_{c,n(c){+}1} \leftarrow m$, where $n(c)$ is the number of elements in row $c$ of $S$ 
%\STATE $\hphantom{S_{c,n(c){+}1}}$ ({\footnotesize If, on average, $n(c)$ is of order 1, then the order of the join algorithm is $\tilde{M}$.})
\STATE $n(c) \leftarrow n(c) + 1$ 
\ENDFOR
\FOR {$m=1$ \TO $\tilde{M}$}
\STATE $\seqof{b_k}{k=1}{3^d} \leftarrow $ get neighboring sub-boxes of $\tilde{X}_{\cdot,m}^{(b)}(t^*)$ s.t. $b_k \in \{1,...,C\}$
\FOR {$k=1$ \TO $3^d$}
\FOR {$j=1$ \TO $n(c_k)$}
\STATE $\ell \leftarrow S_{c_k,j}$
\STATE $v \leftarrow \kappa(\tilde{X}_{\cdot,\ell}^{(f)}(t^*),\tilde{X}^{(b)}_{\cdot,m}(t^*))$
\IF {$v>0$}
	\STATE $L \leftarrow L+1$
	\STATE $u_l \leftarrow u^{(f)}_{\ell} + u^{(b)}_{m}$
	\STATE $r_l \leftarrow r^{(f)}_{\ell} + r^{(b)}_{m}$
	\STATE $w_l \leftarrow v$
\ENDIF
\ENDFOR
\ENDFOR
\ENDFOR
\end{algorithmic}
\end{algorithm}

% flatex input end: [appendix.tex]

%\clearpage

\clearpage
\bibliographystyle{abbrv}
\bibliography{references}

\end{document}